\documentclass{article}

\usepackage{
 amsfonts
,amssymb
,amsmath
,mathabx
,yfonts
,fge
,pifont
}%

\usepackage{graphics}

\newtheorem{theorem}{Theorem}

\newtheorem{corollary}[theorem]{Corollary}
\newtheorem{conj}{Conjecture}{\bf}{\it}
\newtheorem{remark}{Remark}

\begin{document}

\title{
Special functions associated with automorphisms of the space of solutions
to special double confluent Heun equation
}
%

\author{S.I.~Tertychniy%
}

\date{}

\maketitle

\begin{abstract}
The family of quads of interrelated functions holomorphic
on the universal cover of the complex plane without zero
(for brevity, {\mbox{$pqr\!s$}}-functions),
revealing a number of remarkable properties, is introduced.
In particular, under certain conditions the transformations
 of the argument $z$ of {\mbox{$pqr\!s$}}-functions represented 
 by lifts of the replacements
$  z \leftleftharpoons -1/z $,
$  z \leftleftharpoons -z $,
and
$  z \leftleftharpoons 1/z $
are equivalent to linear transformations with known coefficients.
{\mbox{$P\!qr\!s$}}-functions arise in a natural way 
in constructing of certain linear operators
acting as automorphisms on the space of solutions 
to the special double confluent Heun equation (\mbox{sDCHE}).
Earlier such symmetries were known to exist only 
in the case of integer value of one of the constant parameters
when the  predecessors of {\mbox{$pqr\!s$}}-functions appear as polynomials.
In the present work, leaning on the generalized notion of {\mbox{$pqr\!s$}}-functions, 
discrete symmetries of the space of solutions to \mbox{sDCHE}
are extended to the general case, apart from some natural exceptions.
\end{abstract}


\medskip

\vbox{
}

%

\section{Definitions and basic properties }\label{s:010}

Let us consider the following system
of linear homogeneous first order
ODEs 
%
%
%
\begin{eqnarray}
 \label{eq:ppp'}
z^2 {p}'{{}}
&\!\!=\!\!&\;
\big(\mu + ({\ell}-1) z\big) {p}{{}} - {q}{{}}
+ 
z^2 {r}{{}}
,%
 \\
 \label{eq:qqq'}
{q}'{{}}
&\!\!=\!\!&\; 
\big(\lambda - ({\ell}+1) \mu z\big) {p}{{}} + \mu\, {q}{{}}
+
{s}{{}}
,%
 \\
 \label{eq:rrr'}
z^2 {r}'{{}}
&\!\!=\!\!&\;
-\big(\lambda + \mu^2\big) {p}{{}}
+ z \big(2 ({\ell}-1) - \mu z\big) {r}{{}}
- {s}{{}}
,%
 \\
 \label{eq:sss'}
z^2 {s}'{{}}
&\!\!=\!\!&\;
-
\big(\lambda + \mu^2\big) {q}{{}}
 + z^2 \big(\lambda - ({\ell}+1) \mu z\big) {r}{{}}
 + \big(({\ell}-1) z-\mu\big) {s}{{}}.
\end{eqnarray}
Here the symbols
${\ell}$, $\lambda$, $\mu$
denote
some complex constants.
The symbols ${p},{q},{r},{s}$ stand for 
holomorphic functions
of the complex variable
$z$.
For brevity,
we shall refer to
them
as {\em {\mbox{$pqr\!s$}}-functions}.

When resolved with respect to the derivatives
,  i.e.\ upon
division of  {{Eq.s}}{~}\eqref{eq:ppp'}-\eqref{eq:sss'}
by $z^2$, $1$, $z^2$, $z^2$, respectively,
the coefficients in their
right-hand sides
become rational
functions holomorphic
everywhere except for some of them at zero.
Hence 
all solutions
to the above system
are 
holomorphic in some vicinity of any point $z_0\not=0$.

One
may also regard 
{\mbox{$pqr\!s$}}-functions
as solutions of the Cauchy problem
for {{Eq.s}}~\eqref{eq:ppp'}-\eqref{eq:sss'}
with
arbitrary (but 
not totally
null)
initial data specified at arbitrary given $z=z_0\not=0$.
Obviously,
such local solution 
can be analytically continued to any other point of $\mathbb{C}$ except zero.
In particular,
all solutions  to {{Eq.s}}{~}\eqref{eq:ppp'}-\eqref{eq:sss'}
(i.e.\ {\mbox{$pqr\!s$}}-functions)
are
single-valued holomorphic functions
on any connected and simply connected subset
of $ \mathbb{C}^*=\mathbb{C}{\,\fgebackslash\,}\{0\}$
.

At the same time it has to be noted
that,
except for  
the very special conditions,
the natural (inextendible)
domain of holomorphicity of
{\mbox{$pqr\!s$}}-functions is
neither $\mathbb{C}^*$ nor any its subset 
but rather
{\em the universal cover} of $\mathbb{C}^*$.
%
This is the
Riemann surface
$\tilde{\mathbb{C}}\mathstrut^*$
diffeomorphic to
$\mathbb{C}$,
the covering projection
$\Pi:\tilde{\mathbb{C}}\mathstrut^*\simeq\mathbb{C}
\mapsto \mathbb{C}^*$
being realized by the natural exponential function.
However, in what follows,
we shall consider, unless otherwise specified, only a
part of  $\tilde{\mathbb{C}}\mathstrut^*$  (subdomain)
denoting it {{\mbox{$\mathstrut^\backprime{}\mathbb{C}^*$}}}.
It is representable 
by
the result  of the 
removing 
from $\mathbb{C}^*$  of the ray ${\mathbb R}_-$  of negative reals,
$^\backprime\mathbb{C}^* = \mathbb{C}^*{\,\fgebackslash\,}{\mathbb R}_- $.
When considered on {{\mbox{$\mathstrut^\backprime{}\mathbb{C}^*$}}}, 
any instance of {\mbox{$pqr\!s$}}-functions
combines the four
single-valued holomorphic functions
uniquely defined by their values (which may be arbitrary but all zero) at any
given point $ z_0\in {{\mbox{$\mathstrut^\backprime{}\mathbb{C}^*$}}}$.
The two their
single-side continuations to $ {\mathbb R}_- $ also
exist giving rise to real analytic functions in
the common domain $ {\mathbb R}_- $.
However,
as a rule, these one-side limits do not coincide pointwise.

{\mbox{$P\!qr\!s$}}-functions reveal 
a number of noteworthy properties.
The first of them is expressed
by the following statement.
                       \begin{theorem}\label{t:010}
Let {\mbox{$pqr\!s$}}-functions
obey at $z={\mathrm{i}}$ 
the constraint
\begin{equation}
\label{eq:050}
{q}({\mathrm{i}})-\mu{p}({\mathrm{i}})+{r}({\mathrm{i}})=0.
\end{equation}
Then the following equalities 
\begin{eqnarray}
\label{eq:pppA}
{p}(-1/z)
&\!\!=\!\!&
-e^{{\mathrm{i}}{\ell}\pi}
z^{2 \left(1 - {\ell}\right)}
{p}(z)
,
\\
\label{eq:qqqA}
{q}(-1/z)
&\!\!=\!\!&
e^{{\mathrm{i}}{\ell}\pi}
z^{-2{\ell}}
\left(
\mu\, {p}(z)
+z^2 {r}(z)\right)
,
\\
\label{eq:rrrA}
{r}(-1/z)
&\!\!=\!\!&
e^{{\mathrm{i}}{\ell}\pi}
z^{2 (1 - {\ell})} \left(\mu z^2 {p}(z) + {q}(z)\right)
,
\\
\label{eq:sssA}
{s}(-1/z)
&\!\!=\!\!&
-e^{{\mathrm{i}}{\ell}\pi}
z^{-2{\ell}} \left(
\mu\big(\mu z^2 {p}(z) + {q}(z)\big)
+ 
z^2 \left(\mu z^2 {r}(z) + {s}(z)\right)
\right);
\end{eqnarray}
hold true.  
Conversely, {{Eq.}}~\eqref{eq:050}
follows from {{Eq.s}}~\eqref{eq:pppA}-\eqref{eq:sssA}
evaluated at  $z={\mathrm{i}}$.
\end{theorem}
         \begin{remark}\label{r:010}
\hangindent=2ex
\rm
 {{Eq.}}~\eqref{eq:050} is obviously
implied 
by {{Eq.}}~\eqref{eq:qqqA} alone.
The remaining
three equations, when
evaluated at  $z={\mathrm{i}}$,
either turn out to be
 fulfilled identically or
follow from
{{Eq.}}~\eqref{eq:050} (and, thus, from {{Eq.}}~\eqref{eq:qqqA}).
\end{remark}
          \begin{remark}\label{r:025}
\hangindent=2ex
\rm
The constraint \eqref{eq:050}
does not affect  
the value of the function ${s}$ at the selected point
and, moreover,
${s}(\cdot)$
is present
only in
{{Eq.}}~\eqref{eq:sssA}
which might be considered as decoupled from the preceding
ones. However, there is an indirect
influence of the selection of ${s}$
(via the unrestricted setting of ${s}({\mathrm{i}})$)
to the other {\mbox{$pqr$}}-functions
in view of their ``unbreakable interrelation'' implied by
{{Eq.s}}{~}\eqref{eq:ppp'}-\eqref{eq:sss'}.
\end{remark}
          \begin{remark}\label{r:020}
\hangindent=2ex
\rm
The involutive
transformation 
which we shall here refer to as 
{\sl the transformation} A,
signified 
in the
left-hand sides of
{{Eq.s}}~\eqref{eq:pppA}-\eqref{eq:sssA}
by the replacement
    \begin{equation}%
    \label{eq:100}
    z \leftleftharpoons -1/z
   \end{equation}%
of the argument $z$ of the functions involved,
is here
tacitly regarded
as the map
keeping  
the particular argument $z={\mathrm{i}} $
unchanged.
This point is worth mentioning 
because
in the case we deal with,
i.e.~for functions possessing domains distinct of \mbox{$\mathbb{C}^*$},
 ``the reflected imaginary unit''
$-{\mathrm{i}}$ is {\em not\/}  a fixed point of the implied transformation of the arguments albeit
$-1/(-{\mathrm{i}})=(-{\mathrm{i}})$, formally.
Moreover, 
there is {\em another\/} transformation (let us denote it $\tilde{\mathrm{A}}$)
or, one might say,
another implementation  
of the rule \eqref{eq:100}
recognizing 
just $-{\mathrm{i}}$, but not $+{\mathrm{i}}$,
as a
fixed point in the domain of a
(this time)
$\tilde{\mathrm{A}}$-transformed functions.

Accordingly,
as long as we consider
$+{\mathrm{i}}$ as the fixed point 
of the transformation 
signified by the 
argument replacement \eqref{eq:100},
there exist connected and simply connected
open sets containing $+{\mathrm{i}}$ and contained in {\mbox{$\mathstrut^\backprime{}\mathbb{C}^*$}}{} 
such that
their
images through the transformation  ${\mathrm{A}}$
also contain $+{\mathrm{i}}$ and are contained in {\mbox{$\mathstrut^\backprime{}\mathbb{C}^*$}}. 
On them,
the asserted relations expressed 
by {{Eq.s}}~\eqref{eq:pppA}-\eqref{eq:sssA}
  are well defined.
It is here preferable to consider
the subdomain
{{\mbox{$\mathstrut^\backprime{}\mathbb{C}^*$}}}
only.
The extension of
{{Eq.s}}~\eqref{eq:pppA}-\eqref{eq:sssA}
to the whole domain of {\mbox{$pqr\!s$}}-functions
(the universal cover of $\mathbb{C}^*$) by means of
analytic continuation is obviously feasible
although
in general it might prove to be  
not representable {\it by the original formulas}.%

To clarify some specialties of the above interpretation, we consider the following example.
Let 
$z$ be continuously moving from ${\mathrm{i}}\in\,${{\mbox{$\mathstrut^\backprime{}\mathbb{C}^*$}}}
towards some  $x\in{\mathbb R}_+\subset\,${{\mbox{$\mathstrut^\backprime{}\mathbb{C}^*$}}} along a concave curve.
Then
$-1/z$, also starting from $+{\mathrm{i}}$ but further differing from $z$,
is moving around zero
in the opposite angular direction,
arriving
ultimately
at $-x^{-1}\in{\mathbb R}_-$
{\em which  does not belong to} {{\mbox{$\mathstrut^\backprime{}\mathbb{C}^*$}}}.
Thus, when dragging $z$
farther   across ${\mathbb R}_+$ 
inward 
the half-plane $\Im z<0 $, the corresponding
(A-transformed argument)
$-1/z$ {\it leaves\/} {{\mbox{$\mathstrut^\backprime{}\mathbb{C}^*$}}}
across `the upper edge' of the cut along the ray ${\mathbb R}_-$. 
Notice
that we may  not  consider
it entering {{\mbox{$\mathstrut^\backprime{}\mathbb{C}^*$}}} again through the lower cut edge
disconnected from the upper one.
This means that
in the course of the above process 
the literal applicability 
of the formulas \eqref{eq:pppA}-\eqref{eq:sssA} breaks down on the ray of positive reals.
Thus,
to ensure their meaningfulness,
one is compelled to 
obey the restriction 
$\Im z>0$. %
At the same time, it is obvious that
such a limitation is only
a consequence of certain  
simplification we had adopted for convenience.
It
would not arise
in case of consideration of {\mbox{$pqr\!s$}}-functions on their full %
 domain.
However, then  yet another
complication related to certain %
non-uniqueness of
interpretation of {{Eq.s}}~\eqref{eq:pppA}-\eqref{eq:sssA} would appear.
In total, 
 we still prefer here to restrict consideration  to
the subdomain {{\mbox{$\mathstrut^\backprime{}\mathbb{C}^*$}}}
 keeping in mind limitations induced 
by such a simplification.

\end{remark}

{{Eq.}}~\eqref{eq:050}
singles out some subset of {\mbox{$pqr\!s$}}-functions
constraining 
their values
(i.e.\ the initial data for {{Eq.s}}{~}\eqref{eq:ppp'}-\eqref{eq:sss'})
at $z={\mathrm{i}}$.
Yet 
another property 
leans on
their parameterizing
by the values
at $z=1$. It reads as follows.
\begin{theorem} \label{t:020}
Let $\lambda+\mu^2\not=0$ and
{\mbox{$pqr\!s$}}-functions obey the constraints
\begin{eqnarray}
\label{eq:110}
\mu\,{p}(1)+{q}(1)+{r}(1)&=&0,
\\
\label{eq:120}
\lambda\,{p}(1)-\mu\,{q}(1)+{s}(1)&=&0.
\end{eqnarray}
Then the following equalities 
\begin{eqnarray}
\label{eq:pppC}
{p}(1/z)
&\!\!=\!\!&
-(\lambda+\mu^2)^{-1}
z^{2 \left(1 - {\ell}\right)}
\left(
\mu{}z^2{r}(z)+{s}(z)
\right),
\\
\label{eq:qqqC}
{q}(1/z)
&\!\!=\!\!&
-(\lambda+\mu^2)^{-1}
z^{-2{\ell}}
\left(
\lambda{}z^2{r}(z)
-\mu\,{s}(z)
\right)
,
\\
\label{eq:rrrC}
{r}(1/z)
&\!\!=\!\!&
-z^{2 (1 - {\ell})}
\left(\mu z^2 {p}(z) + {q}(z)\right)
,
\\
\label{eq:sssC}
{s}(1/z)
&\!\!=\!\!&
-z^{-2{\ell}}
\left(
\lambda{}z^2{p}(z) - \mu\,{q}(z)
\right);
\end{eqnarray}
hold true. 
Conversely, {{Eq.s}}~\eqref{eq:110}, \eqref{eq:120}
follow from {{Eq.s}}~\eqref{eq:pppC}-\eqref{eq:sssC}
evaluated at  $z=1$.
\end{theorem}
\begin{remark}\label{r:030}
\hangindent=2ex
\rm
For $z=1$
the replacement
of argument of the functions
on the left in {{Eq.s}}{~}\eqref{eq:pppC}-\eqref{eq:sssC}
(we shall refer to it as \textit{ the transformation} C)
reveals no effect.
Accordingly,
there exist  open sets containing 
$+1$ which remain
invariant under the action of the transformation $\mathrm{C}$.
Then it is reasonable to consider first the equalities \eqref{eq:pppC}-\eqref{eq:sssC}
on such neighborhoods of the unity and then utilize analytic continuation
for their extending to greater domains. 
\end{remark}
\begin{remark}\label{r:040}
\hangindent=2ex
\rm
Besides $z=1$, the point $z=-1$
(excluded, by definition, from {{\mbox{$\mathstrut^\backprime{}\mathbb{C}^*$}}})
is also unaffected by the replacement
$ z \leftleftharpoons 1/z $
utilized in {{Eq.s}}{~}\eqref{eq:pppC}-\eqref{eq:sssC}, formally.
However, it can not be considered as a fixed point
of the transformation C. 
More precisely,
claiming  of $z=-1$ to be a fixed point,
one must replace the transformation $\mathrm{C}$
by ``yet another implementation'' $\tilde{\mathrm{C}}$
of the above argument replacement.
For it, 
the former fixed point
$z=1$ loses such a property.
Besides, for 
$\tilde{\mathrm{C}}$,
the associated
\hbox{(sub-)domain} of {\mbox{$pqr\!s$}}-functions, playing role of
 {\,{\mbox{$\mathstrut^\backprime{}\mathbb{C}^*$}}}, has to contain $ {\mathbb R}_- $
but not $ {\mathbb R}_+ $.
Having thus noted the presence of certain  ambiguity in the interpretation of 
  {{Eq.s}}{~}\eqref{eq:pppC}-\eqref{eq:sssC},
we limit ourselves with the above remark and
shall not consider here this issue
 in greater details.
\end{remark}

Combining
conditions of the two above theorems we obtain one more relationship
in accordance with the following. 
                                   \begin{theorem} \label{t:030}
Let {\mbox{$pqr\!s$}}-function
obey 
the conditions
of  both Theorem \ref{t:010} and Theorem \ref{t:020},
i.e.\
they
are holomorphic on a connected and simply connected open set containing
$+{\mathrm{i}}$ and $+1$ and
 meet the constraints \eqref{eq:050},
\eqref{eq:110}, and \eqref{eq:120}.
Then the equalities 
{
\hspace{-2ex}
\begin{eqnarray}
\label{eq:pppB}
{{{\mathcal{M}}}^{1/2}}[{p}]{{}}
&\!\!=\!\!&
e^{{\mathrm{i}}{ }{\ell}\pi} (\lambda+\mu^2)^{-1}
\big(\mu{}z^2{r}{{}}+{s}{{}}\big)
,
\\
\label{eq:qqqB}
{{{\mathcal{M}}}^{1/2}}[{q}]{{}}
&\!\!=\!\!&
-e^{{\mathrm{i}}{ }{\ell}\pi}
\left(
      (\mu\, z^2 {p}{{}} +{q}{{}} )
+\mu(\lambda+\mu^2)^{-1}z^2
      (\mu{}z^2{r}{{}}+{s}{{}} )
\right)
,
\\
\label{eq:rrrB}
{{{\mathcal{M}}}^{1/2}}[{r}]{{}}
&\!\!=\!\!&
-e^{{\mathrm{i}}{ }{\ell}\pi}
{r}{{}},
\\
\label{eq:sssB}
{{{\mathcal{M}}}^{1/2}}[{s}]{{}}
&\!\!=\!\!&
e^{{\mathrm{i}}{ }{\ell}\pi}
\left(
(\lambda+\mu^2) {p}{{}}+\mu{}z^2{r}{{}}
\right),
\end{eqnarray}
}%
hold true, 
where the arguments $z$ of all the functions coincide and hence are suppressed,
and
where the operator 
${{{\mathcal{M}}}^{1/2}}{}  $
carries out analytic continuation
of the function it acts to
along the circular arc started at $z$,
centered at zero,
subtending an angle $\pi$, and
 oriented counter-clockwise. %
Moreover,
the products of {\mbox{$pqr\!s$}}-functions times
the power function $ z^{-{\ell}}  $
are single-valued and holomorphic on $\mathbb{C}^*$.
\end{theorem}
\begin{remark}\label{r:050}
\hangindent=2ex
\rm
As opposed to transformations of {\mbox{$pqr\!s$}}-functions
treated
by
Theorems  \ref{t:010} and \ref{t:020},
the transformation of arguments of functions on the left in
{{Eq.s}}{~}\eqref{eq:pppB}-\eqref{eq:sssB}
(let us call it \textit{the transformation} $\mathrm{B}$)
  admits 
no fixed points 
and is not involutive.
Moreover, applying
the transformation $\mathrm{B}$
twice, the resulting effect
turns into the analytic continuation of the function to be transformed
{\it along the loop projected to
{\rm(essentially, coinciding with)}
the full circle}. 
Such kind of analytic continuation %
around a singular point (in our case, the center $z=0$)
is commonly named the 
{\it monodromy} transformation. We denote it by the
symbol ${\mathcal{M}}$.
We have therefore ${{{\mathcal{M}}}^{1/2}}\circ{{{\mathcal{M}}}^{1/2}}={\mathcal{M}}$ by definition.
The effect of the operator 
${{{\mathcal{M}}}^{1/2}}$ can thus be named semi-monodromy transformation.

In our case
${\mathcal{M}}$
is the linear operator which
sends, in particular, the values of {\mbox{$pqr\!s$}}-functions on the ``lower'' edge of the cut along
the ray
${\mathbb R}_-$
to the (generally speaking, distinct)
values they assume on its ``upper'' edge.
Since {\mbox{$pqr\!s$}}-functions obey on the both edges the same system
\eqref{eq:ppp'}-\eqref{eq:sss'}
of linear homogeneous ODEs
(since their coefficients are invariant with respect to ${\mathcal{M}}$)
such a transformation
is represented
by a {\em constant} $4\times 4$ matrix.
\end{remark}
\begin{remark}\label{r:060}
\hangindent=2ex
\rm
If $z\in\,${{\mbox{$\mathstrut^\backprime{}\mathbb{C}^*$}}} and $\Im z<0$ then
${{{\mathcal{M}}}^{1/2}}{}z \;\big(={{{\mathcal{M}}}^{1/2}}[\mathrm{Id}](z)\,\big) =e^{{\mathrm{i}}\pi}z=  -z\in\,${{\mbox{$\mathstrut^\backprime{}\mathbb{C}^*$}}}.
However, if $\Im z\ge 0$ then 
an application  
of ${{{\mathcal{M}}}^{1/2}}  $
would yield 
the argument of {\mbox{$pqr\!s$}}-functions
in {{Eq.s}}{~}\eqref{eq:pppB}-\eqref{eq:sssB}
on the left which does not belong to
{{\mbox{$\mathstrut^\backprime{}\mathbb{C}^*$}}}.
Evading such a  complication,
we shall assume $\Im z<0$
for simplicity
unless otherwise specified.
Analytic continuation has to be applied
for relaxation
of the 
limitation %
and
extending
the local form of the equalities \eqref{eq:pppB}-\eqref{eq:sssB}
in which ${{{\mathcal{M}}}^{1/2}}  $-transformation is regarded as  the
inversion of sign of the function argument
to a greater domain.
\end{remark}
\begin{remark}\label{r:070}
\hangindent=2ex
\rm
In general case, given a prescribed set of constant parameters,
simultaneous
fulfillment of   {{Eq.}}~\eqref{eq:050} and
{{Eq.s}}~\eqref{eq:110}, \eqref{eq:120}
for the same instance of {\mbox{$pqr\!s$}}-functions
should
be achievable by  means of their %
appropriate selection.
Indeed,
the set of all
{\mbox{$pqr\!s$}}-functions can be indexed
by the quad of their values at $z=1$
fixed up to multiplication by an
insignificant (associated with a decoupled degree of freedom)
non-zero common factor, i.e.\ by  points of a
projective space $\mathbb{C}\mathrm{P}^3_{\{1\}}$.
The two linear equations  \eqref{eq:110},\eqref{eq:120}
single out the projective line
embedded therein.
This projective line
is conveyed (pushedforward) by
the vector flow associated with 
the equations
\eqref{eq:ppp'}-\eqref{eq:sss'}
into another projective space $\mathbb{C}\mathrm{P}^3_{\{{\mathrm{i}}\}}$
indexing the same set of
{\mbox{$pqr\!s$}}-functions %
by their values
(also considered up to a common constant factor)
at $z={\mathrm{i}}$.
In the latter projective space, the equation
\eqref{eq:050}
singles out
certain embedded
projective plane.
The question
equivalent to  the issue
of consistency of
{{Eq.}}~\eqref{eq:050} with
{{Eq.s}}~\eqref{eq:110} and \eqref{eq:120}
reads: whether the
former (conveyed) 
projective line
intersects the latter projective plane or not?
This problem remains open
yet but
numerical computations
point in favor of the
affirmative upshot, at least,
under apparently generic conditions.
Thus,
most plausibly,
inconsistency of {{Eq.}}~\eqref{eq:050} with
{{Eq.s}}~\eqref{eq:110} and \eqref{eq:120}
and the subsequent
inanity 
of Theorem \ref{t:030}, if any,
could only occur
 under the very special conditions (currently unknown).
We may state therefore the following.
 \end{remark}
\begin{corollary}\label{c:010}
There exists a 
set of 
pairwise
linearly independent quads of
holomorphic functions 
${p},{q},{r},{s}$
parameterized by 
points of
$\mathbb{C}\mathrm{P}^2$
such that the equations \eqref{eq:pppA}-\eqref{eq:sssA}
are fulfilled.
\end{corollary}
\begin{corollary}\label{c:020}
There exists a   set of
pairwise
linearly independent quads of
of holomorphic functions 
${p},{q},{r},{s}$
parameterized by 
points of
$\mathbb{C}\mathrm{P}^1$
such that the equations \eqref{eq:pppC}-\eqref{eq:sssC}
are fulfilled.
\end{corollary}
\begin{conj}
For almost all values of the constant parameters
there exists a discrete set of quads of  functions
${p},{q},{r},{s}$ holomorphic  on the universal cover of $\mathbb{C}^*$
such that
the equations \eqref{eq:pppA}-\eqref{eq:sssA},
\eqref{eq:pppC}-\eqref{eq:sssC},
and \eqref{eq:pppB}-\eqref{eq:sssB}
are fulfilled.
  \end{conj}
The last assertion of the Theorem \ref{t:030} says
how the 
{\mbox{$pqr\!s$}}-functions referred to in the above Conjecture
are expressed through functions 
which are
single-valued and 
holomorphic on $\mathbb{C}^*$.

  \medskip

We proceed now with proofs of the three above theorems.

\smallskip

\noindent
{\sl
Proof of
Theorem \ref{t:010}. } 
Let us denote the four differences of the left- and right-hand sides
of
{{Eq.s}}{~}\eqref{eq:pppA},
\eqref{eq:qqqA},
\eqref{eq:rrrA},
\eqref{eq:sssA},
by the symbols
$
{\mbox{$\mathstrut^{A}\!\Delta_{p}$}}
,
{\mbox{$\mathstrut^{A}\!\Delta_{p}$}}
,
{\mbox{$\mathstrut^{A}\!\Delta_{r}$}}
,
{\mbox{$\mathstrut^{A}\!\Delta_{s}$}}
$
respectively, considering them, as they stand, as the functions of $z$.
For example, one  of such definitions reads 
$
{\mbox{$\mathstrut^{A}\!\Delta_{p}$}}(z)=
{p}(-1/z)
+e^{{\mathrm{i}}{\ell}\pi}
z^{2 \left(1 - {\ell}\right)}
{p}(z), \mbox{\it etc}.
$
As it is shown in   Appendix \ref{a:010},
they obey the following system of linear homogeneous
ODEs
 \begin{equation}
                      \label{eq:210}
\relax\hspace{-0.2em}
\begin{aligned}
z^2
\frac{d}{d z}
{\mbox{$\mathstrut^{A}\!\Delta_{p}$}}{{}}
=&\;
z(1 - {\ell} + \mu{}z) 
\hspace{0.05em}{\mbox{$\mathstrut^{A}\!\Delta_{p}$}}{{}}
-  z^2 \,{\mbox{$\mathstrut^{A}\!\Delta_{q}$}}{{}}
+{\mbox{$\mathstrut^{A}\!\Delta_{r}$}}{{}},
\\
z^3
\frac{d}{d z}
{\mbox{$\mathstrut^{A}\!\Delta_{q}$}}{{}}
=&\;
(({\ell} - 1)\mu + \lambda{}z)
\hspace{0.05em} {\mbox{$\mathstrut^{A}\!\Delta_{p}$}}{{}}
+ \mu{}z
  \,{\mbox{$\mathstrut^{A}\!\Delta_{q}$}}{{}}
+ z
\,{\mbox{$\mathstrut^{A}\!\Delta_{s}$}}{{}},
\\
z^2
\frac{d}{d z}
{\mbox{$\mathstrut^{A}\!\Delta_{r}$}}{{}}
=&
- (\lambda+\mu^2)z^2
 \,{\mbox{$\mathstrut^{A}\!\Delta_{p}$}}{{}}
 -\big( \mu + 2({\ell}-1)z\big)
 {\mbox{$\mathstrut^{A}\!\Delta_{r}$}}{{}}
  - z^2\,{\mbox{$\mathstrut^{A}\!\Delta_{s}$}}{{}}
 \\
z^3
\frac{d}{d z}
{\mbox{$\mathstrut^{A}\!\Delta_{s}$}}{{}}
=&
- (\lambda+\mu^2)z^3
 \,{\mbox{$\mathstrut^{A}\!\Delta_{q}$}}{{}}
+
( ({\ell}+1)\mu + \lambda{}z )
\hspace{0.05em} {\mbox{$\mathstrut^{A}\!\Delta_{r}$}}{{}}
-z^2
({\ell} - 1 + \mu{}z)
\hspace{0.1em} \rlap{$ {\mbox{$\mathstrut^{A}\!\Delta_{s}$}}{{}}, $}
 \end{aligned}
 \end{equation}
provided {{Eq.s}}{~}\eqref{eq:ppp'}-\eqref{eq:sss'}
are fulfilled.

Using the explicit definitions, let us compute
the particular values of
$  {\mbox{$\mathstrut^{A}\!\Delta_{{{\mbox{\tiny\ding{74}}}}}$}}({\mathrm{i}})=
 {\mbox{$\mathstrut^{A}\!\Delta_{{{\mbox{\tiny\ding{74}}}}}$}}(e^{{\mbox{\scriptsize${\mathrm{i}}\over2$}}\pi})
 $
for  {{\mbox{\small\ding{74}}}}$\;\in\{ {p},{q},{r},{s} \} $.
Notice that for such a choice of the argument $z$ 
one has 
 $-1/z=-e^{-{\mbox{\scriptsize${\mathrm{i}}\over2$}}\pi}
={\mathrm{i}},
z^{-2{\ell}}= e^{-{\mathrm{i}}{\ell}\pi},
z^{2(1-{\ell})} 
=-e^{-{\mathrm{i}}{\ell}\pi} $.
Then it follows from  {{Eq.}}~\eqref{eq:pppA}
that
${\mbox{$\mathstrut^{A}\!\Delta_{p}$}}({\mathrm{i}})=0$.
 The values  $ {\mbox{$\mathstrut^{A}\!\Delta_{{{\mbox{\tiny\ding{74}}}}}$}}({\mathrm{i}}) $
 of the other differences are not automatically
zero but one  easily finds that %
in accordance with definitions
$$
 {\mbox{$\mathstrut^{A}\!\Delta_{{{\mbox{\tiny\ding{74}}}}}$}}({\mathrm{i}})=\zeta_{{{\mbox{\tiny\ding{74}}}}}\cdot
\big(
{q}({\mathrm{i}})-\mu{p}({\mathrm{i}}) +{r}({\mathrm{i}})
\big) \mbox{ for {\small\ding{74}}}\in\{{q},{r},{s}\},
\mbox{ and }\zeta_{q}=\zeta_{r}=1, \zeta_{s}=\mu.
$$
Thus if
{{Eq.}}~\eqref{eq:050} is fulfilled
then $ {\mbox{$\mathstrut^{A}\!\Delta_{{{\mbox{\tiny\ding{74}}}}}$}}({\mathrm{i}})=0 $ 
for all `the indices' {{\mbox{\small\ding{74}}}}{}$\;\in \{{p},{q},{r},{s}\} $.
This implies the vanishing everywhere of
all the functions
$ {\mbox{$\mathstrut^{A}\!\Delta_{{{\mbox{\tiny\ding{74}}}}}$}}(z) $ 
in view of  uniqueness
of solutions of the Cauchy problem for {{Eq.s}}{~}\eqref{eq:210} 
with the null initial data posed at $z={\mathrm{i}}$.
\hfill$\square$

\medskip
\noindent
{\sl Proof of
Theorem \ref{t:020}.}
Building on the 
notations utilized in the preceding proof,
we
denote the differences of the left- and right-hand sides
of
{{Eq.s}}{~}\eqref{eq:pppC},
\eqref{eq:qqqC},
\eqref{eq:rrrC},
\eqref{eq:sssC}
by the 
 symbols
$
{\mbox{$\mathstrut^{C\!}\!\Delta_{p}$}}%
(z)
,
{\mbox{$\mathstrut^{C\!}\!\Delta_{q}$}}%
(z)
,
{\mbox{$\mathstrut^{C\!}\!\Delta_{r}$}}%
(z)
,
{\mbox{$\mathstrut^{C\!}\!\Delta_{s}$}}%
(z)
$,
respectively.
It is shown in Appendix \ref{a:020}
that they obey
the following system of linear homogeneous
ODEs
\\
%
%
%
%
%
%
%
%
%
\hspace{-3em}
\begin{equation} 
\label{eq:220}
\begin{aligned}
z^2
\frac{d}{d z}
{\mbox{$\mathstrut^{C\!}\!\Delta_{p}$}}{{}}
=&\; 
-z({\ell} - 1 + \mu{}z) 
\hspace{0.05em}{\mbox{$\mathstrut^{C\!}\!\Delta_{p}$}}{{}}
+z^2\,{\mbox{$\mathstrut^{C\!}\!\Delta_{q}$}}{{}}
-{\mbox{$\mathstrut^{C\!}\!\Delta_{r}$}}{{}}
,
 %
 \\
%
z^3
\frac{d}{d z}
{\mbox{$\mathstrut^{C\!}\!\Delta_{q}$}}{{}}
=&\; 
 (({\ell}+1)\mu - \lambda{}z)
\hspace{0.05em}{\mbox{$\mathstrut^{C\!}\!\Delta_{p}$}}{{}}
-\mu{} z \,{\mbox{$\mathstrut^{C\!}\!\Delta_{q}$}}{{}}
-z\,{\mbox{$\mathstrut^{C\!}\!\Delta_{s}$}}{{}},
 \\
%
z^2
\frac{d}{d z}
{\mbox{$\mathstrut^{C\!}\!\Delta_{r}$}}{{}}
=&\; 
(\lambda+\mu^2) z^2\,{\mbox{$\mathstrut^{C\!}\!\Delta_{p}$}}{{}}
+(\mu + 2(1 - {\ell})z)
\hspace{0.05em}{\mbox{$\mathstrut^{C\!}\!\Delta_{r}$}}{{}}
+z^2\,{\mbox{$\mathstrut^{C\!}\!\Delta_{s}$}}{{}},
\\ 
%
z^3
\frac{d}{d z}
{\mbox{$\mathstrut^{C\!}\!\Delta_{s}$}}{{}}
=&\; 
   (\lambda+\mu^2) z^3 \,{\mbox{$\mathstrut^{C\!}\!\Delta_{q}$}}{{}}
  +(({\ell} + 1)\mu - \lambda{} z)
\hspace{0.05em}{\mbox{$\mathstrut^{C\!}\!\Delta_{r}$}}{{}}
  +z^2 (1 - {\ell} + \mu{}z)
\hspace{0.05em}{\mbox{$\mathstrut^{C\!}\!\Delta_{s}$}}{{}},
\end{aligned}
\end{equation} 
provided {{Eq.s}}{~}\eqref{eq:ppp'}-\eqref{eq:sss'}
are fulfilled.

We compute now
the values  the differences
{\mbox{$\mathstrut^{C\!}\!\Delta_{{{\mbox{\tiny\ding{74}}}}}$}}
({{\mbox{\small\ding{74}}}}$\;\in\{{p}, {q}, {r}, {s}\}$)
acquire
after the
plugging $z\leftleftharpoons1$ in their definitions.
The result is as follows:
%
%
\begin{eqnarray} 
(\lambda+\mu^2)
\hspace{0.05em}{\mbox{$\mathstrut^{C\!}\!\Delta_{{{\mbox{\tiny\ding{74}}}}}$}}(1)
\hspace{-0.1em}&=&\hspace{-0.1em}
\zeta_{ {{{\mbox{\tiny\ding{74}}}}} }\cdot
\big({s}(1) - \mu {q}(1) + \lambda{p}(1)   \big)
\nonumber
\\&&
+
\sigma_{ {{{\mbox{\tiny\ding{74}}}}}}\cdot
\big({r}(1) +    {q}(1)  + \mu{p}(1)\big),
\nonumber
\\
&\llap{whe}\mbox{re}&
\zeta_{{p}}=1, \zeta_{{q}}=-\mu, \zeta_{{r}}=0, \zeta_{{s}}= \lambda+\mu^2,
\nonumber \\
&&
\sigma_{{p}}=\mu, \sigma_{{q}}=\lambda, \sigma_{{r}}=\lambda+\mu^2, \sigma_{{s}}=0.
\nonumber
\end{eqnarray}
Thus if the constraints
\eqref{eq:110} and \eqref{eq:120}
are fulfilled then
all the 
 differences
${\mbox{$\mathstrut^{C\!}\!\Delta_{{{\mbox{\tiny\ding{74}}}}}$}}(z)$ vanish at $z=1$.
But then
they are the identically zero functions,
 ${\mbox{$\mathstrut^{C\!}\!\Delta_{{{\mbox{\tiny\ding{74}}}}}$}}(z)\equiv0$,
as a consequence of
{{Eq.s}}{~}\eqref{eq:220}.
This means exactly that {{Eq.s}}{~}\eqref{eq:pppC}-\eqref{eq:sssC}
hold true.
\hfill$\square$

\medskip
\noindent
{\sl Proof of
Theorem \ref{t:030}.}
As above,
let us
denote the differences of the left- and right-hand sides
of
{{Eq.s}}{~}\eqref{eq:pppB},
\eqref{eq:qqqB},
\eqref{eq:rrrB},
\eqref{eq:sssB}
by the 
 symbols
$
{\mbox{$\mathstrut^{B\!}\!\Delta_{p}$}}
(z)
,$ $ {\mbox{$\mathstrut^{B\!}\!\Delta_{q}$}}
(z)
,$ $ {\mbox{$\mathstrut^{B\!}\!\Delta_{r}$}}
(z)
,$ $ {\mbox{$\mathstrut^{B\!}\!\Delta_{s}$}}
(z)$,
respectively.
It is shown in  Appendix \ref{a:030} that
in case of fulfillment of
{{Eq.s}}{~}\eqref{eq:ppp'}-\eqref{eq:sss'}
they obey the following system of linear homogeneous
ODEs
\begin{equation} 
\label{eq:230}
\begin{aligned}
z^2
\frac{d}{d z}
{\mbox{$\mathstrut^{B\!}\!\Delta_{p}$}}{{}}
=&\; 
(({\ell} - 1)z - \mu)    {\mbox{$\mathstrut^{B\!}\!\Delta_{p}$}}{{}}
+                       {\mbox{$\mathstrut^{B\!}\!\Delta_{q}$}}{{}}
-z^2                    {\mbox{$\mathstrut^{B\!}\!\Delta_{r}$}}{{}},
\\
\frac{d}{d z}
{\mbox{$\mathstrut^{B\!}\!\Delta_{q}$}}{{}}
=&\; 
-(\lambda + ({\ell}+1)\mu{} z)     {\mbox{$\mathstrut^{B\!}\!\Delta_{p}$}}{{}}
-\mu                            \,{\mbox{$\mathstrut^{B\!}\!\Delta_{q}$}}{{}}
-                                 {\mbox{$\mathstrut^{B\!}\!\Delta_{s}$}}{{}},
\\
z^2
\frac{d}{d z}
{\mbox{$\mathstrut^{B\!}\!\Delta_{r}$}}{{}}
=&\; 
(\lambda+\mu^2)                {\mbox{$\mathstrut^{B\!}\!\Delta_{p}$}}{{}}
+(2({\ell} - 1) + \mu{}z)z    \,{\mbox{$\mathstrut^{B\!}\!\Delta_{r}$}}{{}}
+                              {\mbox{$\mathstrut^{B\!}\!\Delta_{s}$}}{{}},
\\
\hspace{-2em}
z^2
\frac{d}{d z}
{\mbox{$\mathstrut^{B\!}\!\Delta_{s}$}}{{}}
=&\; 
(\lambda+\mu^2)                         {\mbox{$\mathstrut^{B\!}\!\Delta_{q}$}}{{}}
-(\lambda +  ({\ell} + 1)\mu{}z)z^2    \,{\mbox{$\mathstrut^{B\!}\!\Delta_{r}$}}{{}}
+(\mu + ({\ell} - 1)z) \hspace{0.05em}    {\mbox{$\mathstrut^{B\!}\!\Delta_{s}$}}{{}}.
\end{aligned}
\end{equation} 

The next step should assume
 computation of the particular values
 $
 {\mbox{$\mathstrut^{B\!}\!\Delta_{{{\mbox{\tiny\ding{74}}}}}$}}
(-{\mathrm{i}}),$  ${{\mbox{\small\ding{74}}}}\;\in\{{p},{q},{r},{s}\}$.
 However,
 carrying out this by means of the mere
 substitutions $z\leftleftharpoons -{\mathrm{i}}  $
 into the definitions of $ {\mbox{$\mathstrut^{B\!}\!\Delta_{{{\mbox{\tiny\ding{74}}}}}$}}$,
 some ambiguity may arise
 due to possibility of 
 overlapping 
of sheets of the branching domain {\mbox{$pqr\!s$}}-functions live on. 
 To make 
the computation univocal, we consider
first the ``deformed''
versions
{\raisebox{0.95ex}{$\raisebox{-0.2ex}[0ex][0ex]{\scriptsize${{\epsilon}}$}
\atop\raisebox{-0.2ex}[0ex][0ex]{$\mathstrut^{B\hspace{-0.2ex}}\hspace{-0.5ex}\Delta_{{{\mbox{\tiny\ding{74}}}}}$}$}}
of the differences
 $ {\mbox{$\mathstrut^{B\!}\!\Delta_{{{\mbox{\tiny\ding{74}}}}}$}}
 $, where $\epsilon$ plays role of the deformation parameter.
 Their distinction
 is that
in case of
{\raisebox{0.95ex}{$\raisebox{-0.2ex}[0ex][0ex]{\scriptsize${{\epsilon}}$}
\atop\raisebox{-0.2ex}[0ex][0ex]{$\mathstrut^{B\hspace{-0.2ex}}\hspace{-0.5ex}\Delta_{{{\mbox{\tiny\ding{74}}}}}$}$}}
the 
factor in argument of the {\mbox{$pqr\!s$}}-function
on the left  is distinct of the one
involved
 in {{Eq.s}}~\eqref{eq:pppB} - \eqref{eq:sssB}, see Remark \ref{r:060}.
The common exponential multiplier on the right is also modified.
 Namely, let
 the factor $ e^{{\mathrm{i}}{{\epsilon}}\pi} $, where
${{\epsilon}}\in [0,1] $  is the auxiliary real parameter,
 be used instead of
  $-1= e^{{\mathrm{i}}\pi} $.
For example, one has
$
{\raisebox{0.95ex}{$\raisebox{-0.2ex}[0ex][0ex]{\scriptsize${{\epsilon}}$}
\atop\raisebox{-0.2ex}[0ex][0ex]{$\mathstrut^{B\hspace{-0.2ex}}\hspace{-0.5ex}\Delta_{r}$}$}}
(z)
= {r}(e^{{\mathrm{i}}{{\epsilon}}\pi}z )
+ e^{{\mathrm{i}}{\ell}{{\epsilon}}\pi}{r}(z)  $ while
$
  {\mbox{$\mathstrut^{B\!}\!\Delta_{r}$}}
(z)= {r}(e^{{\mathrm{i}}\pi}z) + e^{{\mathrm{i}}{\ell}\pi}{r}(z)$, {\it etc}.
Explicit
definitions of all
{\raisebox{0.95ex}{$\raisebox{-0.2ex}[0ex][0ex]{\scriptsize${{\epsilon}}$}
\atop\raisebox{-0.2ex}[0ex][0ex]{$\mathstrut^{B\hspace{-0.2ex}}\hspace{-0.5ex}\Delta_{{{\mbox{\tiny\ding{74}}}}}$}$}}
are
furnished
by {{Eq.s}}~\eqref{eq:580}.

Now let us notice that
in the case ${{\epsilon}}=0$ all the arguments of {\mbox{$pqr\!s$}}-functions  utilized
for  computation of
{\raisebox{0.95ex}{$\raisebox{-0.2ex}[0ex][0ex]{\tiny${{0}}$\,}
\atop\raisebox{-0.2ex}[0ex][0ex]{$\mathstrut^{B\hspace{-0.2ex}}\hspace{-0.5ex}\Delta_{{{\mbox{\tiny\ding{74}}}}}$}$}}
coincide with $z$
and no ambiguity
in their evaluation
can thus
arise.
Then, starting from these values,
we carry out analytic continuation
varying $ {{\epsilon}} $ through the segment $[0,1]$.
We 
define
$ {\mbox{$\mathstrut^{B\!}\!\Delta_{{{\mbox{\tiny\ding{74}}}}}$}}(z)$
to be 
``the final values'' the functions
{\raisebox{0.95ex}{$\raisebox{-0.2ex}[0ex][0ex]{\scriptsize${{\epsilon}}$}
\atop\raisebox{-0.2ex}[0ex][0ex]{$\mathstrut^{B\hspace{-0.2ex}}\hspace{-0.5ex}\Delta_{{{\mbox{\tiny\ding{74}}}}}$}$}}$(z)$
arrive at
as ${{\epsilon}}\nearrow 1 $.
Such an interpretation
leaves no room for ambiguity 
in the meaning of
definitions of
$
{\mbox{$\mathstrut^{B\!}\!\Delta_{{{\mbox{\tiny\ding{74}}}}}$}}
$
and, more generally,
the relations {{Eq.s}}~\eqref{eq:pppB}-\eqref{eq:sssB}
represent.

Assuming 
the above interpretation of
$
{\mbox{$\mathstrut^{B\!}\!\Delta_{{{\mbox{\tiny\ding{74}}}}}$}}
$,
it is shown in Appendix \ref{a:040}
that the following equations are fulfilled for arbitrary
functions ${p},{q},{r},{s}$ holomorphic on the circular arc
passing through the point $-{\mathrm{i}},+1,$ and $+{\mathrm{i}}$.
\begin{equation}\label{eq:240}
\begin{aligned}
{\mbox{$\mathstrut^{B\!}\!\Delta_{p}$}}(-{\mathrm{i}})
=& \;
e^{{\mathrm{i}}{}{\ell}\pi} (\lambda+\mu^2)^{-1}
\big(
\mu \,{\mbox{$\mathstrut^{C\!}\!\Delta_{r}$}}({\mathrm{i}})
-
\,{\mbox{$\mathstrut^{C\!}\!\Delta_{s}$}}({\mathrm{i}})
 \big)
,
\\
{\mbox{$\mathstrut^{B\!}\!\Delta_{q}$}}(-{\mathrm{i}})
            =  &\;
\big(
{q}( {\mathrm{i}})
-\mu\,{p}( {\mathrm{i}})
+{r}({\mathrm{i}})
\big)
  \\&
+e^{{\mathrm{i}}{\ell}\pi}
\big(
{\mbox{$\mathstrut^{C\!}\!\Delta_{q}$}}({\mathrm{i}})
- \mu\,{\mbox{$\mathstrut^{C\!}\!\Delta_{p}$}}({\mathrm{i}})
+\mu\,(\lambda+\mu^2)^{-1}
(
\mu
\,{\mbox{$\mathstrut^{C\!}\!\Delta_{r}$}}({\mathrm{i}})
-
{\mbox{$\mathstrut^{C\!}\!\Delta_{s}$}}({\mathrm{i}}
))\big)
,
\\
{\mbox{$\mathstrut^{B\!}\!\Delta_{r}$}}(-{\mathrm{i}})
= &\;
\big(
{q}( {\mathrm{i}})
-\mu\,{p}( {\mathrm{i}})
+{r}({\mathrm{i}})
\big)
+ e^{{\mathrm{i}}{\ell}\pi}
\,{\mbox{$\mathstrut^{C\!}\!\Delta_{r}$}}({\mathrm{i}})
,
\\
{\mbox{$\mathstrut^{B\!}\!\Delta_{s}$}}(-{\mathrm{i}})
 = &\;
 \mu
\big(
{q}( {\mathrm{i}})
-\mu\,{p}( {\mathrm{i}})
+{r}({\mathrm{i}})
\big)
-
 e^{{\mathrm{i}} {\ell}\pi}(\lambda+\mu^2){\mbox{$\mathstrut^{C\!}\!\Delta_{p}$}}({\mathrm{i}})
+\mu\,
 e^{{\mathrm{i}} {\ell}\pi}
{\mbox{$\mathstrut^{C\!}\!\Delta_{r}$}}({\mathrm{i}}).
\end{aligned}
\end{equation}
The symbols
${\mbox{$\mathstrut^{C\!}\!\Delta_{{{\mbox{\tiny\ding{74}}}}}$}}$,
{{\mbox{\small\ding{74}}}}$\;\in\{ {p},{q},{r},{s} \} $,
were already
utilized
in the proof of Theorem \ref{t:020}. They
denote the differences of the left- and right-hand sides of
{{Eq.s}}~\eqref{eq:pppC}-\eqref{eq:sssC}, considered, as they stand,
 as the functions of $z$.
Every equation from the system \eqref{eq:240}
can therefore be regarded as
the coincidence, upon simplifications, of a pair of certain
linear combinations of 4+4 instances of {\mbox{$pqr\!s$}}-functions
of which some are
 evaluated at
$z={\mathrm{i}}$ and
others
 at $z=-{\mathrm{i}}$.

On the other hand,
the conditions of the theorem to be proven
imply, in particular, the fulfilment
of the assertion of Theorem \ref{t:020}
which establishes the vanishing of all the
four functions
$ {\mbox{$\mathstrut^{C\!}\!\Delta_{{{\mbox{\tiny\ding{74}}}}}$}}(z)$
irrespectively of the choice 
of their arguments. Thus all
the terms in
{{Eq.s}}~\eqref{eq:240}
involving those factors
may 
be discarded.

Now, taking into account the fulfillment of
{{Eq.}}~\eqref{eq:050}, we see that all the expressions
on the left in
\eqref{eq:240},
i.e.\ the functions $ {\mbox{$\mathstrut^{B\!}\!\Delta_{{{\mbox{\tiny\ding{74}}}}}$}}(z)$, 
{{\mbox{\small\ding{74}}}}$\;\in\{{p},{q},{r},{s}\} $,
evaluated at $z=-{\mathrm{i}})$,
actually vanish.
Since these functions
obey the system of linear homogeneous first order ODEs
(see {{Eq.s}}~\eqref{eq:230}) they
reduce to identical zero.
This means exactly that the equalities \eqref{eq:pppB}-\eqref{eq:sssB} hold true.

Let us consider the second claim
of the theorem which
establishes, under the restrictions assumed,
a  simpler domain $\mathbb{C^*}$
for the
products 
of {\mbox{$pqr\!s$}}-functions times  $z^{-l}$
as compared to  {\mbox{$pqr\!s$}}-functions themselves
which
are not single-valued on $\mathbb{C^*}$ 
and hence must be considered 
on
its universal cover.
We note first
that the four-element vector
consisting of the right-hand sides of {{Eq.s}}~\eqref{eq:pppB}-\eqref{eq:sssB}
can be obtained by means of the multiplication
of the vector $ ({p}(z),{q}(z),{r}(z),{s}(z))^\top $
by the matrix $ M^{1/2}_{{\ell}}(z)=e^{{\mathrm{i}}{\ell}\pi}M^{1/2}(z) $, where
\begin{equation*}
M^{1/2}(z)=
\begin{pmatrix}
0& 0 & \mu{}z^2(\lambda+\mu^2)^{-1}&(\lambda+\mu^2)^{-1}
\\
-\mu{}z^2 & -1 & -\mu^2{}z^4(\lambda+\mu^2)^{-1} & -\mu{}z^2(\lambda+\mu^2)^{-1}
\\
0&0&-1&0
\\
\lambda+\mu^2 & 0 & \mu^2{}z & 0
\end{pmatrix}.
\end{equation*}
Thus under the conditions assumed the action of the operator $ {{{\mathcal{M}}}^{1/2}} $
to {\mbox{$pqr\!s$}}-functions is completely described by the matrix $ M^{1/2}_{{\ell}}$.
Let us examine   the effect of this action applied twice.
In the language of matrices
it is described by the product
of the ones associated with the mentioned transformation.
However whereas the first
operator $ {{{\mathcal{M}}}^{1/2}} $
of the composition
defined as, in a sense, a rotation of the function argument,
 is associated with  $ M^{1/2}_{{\ell}}(z)$
the second one `starts' with the arguments already rotated
acting separately to the matrix $ M^{1/2}_{{\ell}}(z)$ and to the vector of
{\mbox{$pqr\!s$}}-functions. In other words,
the operator composition $ {{{\mathcal{M}}}^{1/2}} \circ {{{\mathcal{M}}}^{1/2}}$ has to be
associated with
the matrix product
$ {{{\mathcal{M}}}^{1/2}} [M^{1/2}_{{\ell}} ](z)\cdot M^{1/2}_{{\ell}}(z)=
e^{2{\mathrm{i}}{\ell}\pi} {{{\mathcal{M}}}^{1/2}}[M^{1/2}](z)\cdot M^{1/2}(z)$.
We can compute it making use of the following {\it identity}
\begin{equation*}
 M^{1/2}(e^{{\mathrm{i}}\epsilon{\ell}\pi}z)\cdot
 M^{1/2}(z)
\equiv
I
+{(1-e^{2{\mathrm{i}}\epsilon\pi})\mu{}z^2\over\lambda+\mu^2}
\begin{pmatrix}
0&0&1&0\\
\lambda+\mu^2&0&(1-e^{2{\mathrm{i}}\epsilon\pi})\mu{}z^2&1\\
0&0&0&0\\
0&0&\lambda+\mu^2&0
\end{pmatrix}\rlap{,}
\end{equation*}
where $ I $ is the unit matrix and
$\epsilon\in[0,1]$ is the real parameter.
The analytic continuation (``the rotation of the argument'') carried out by the operator
$ {{{\mathcal{M}}}^{1/2}} $ can be represented
by the evaluation of
the limit as $ {{{\epsilon}}\nearrow 1} $. Then the factor in parenthesis
in front of the second summand on the right goes to zero
and we obtain ${{{\mathcal{M}}}^{1/2}}[M^{1/2}](z)\cdot M^{1/2}(z)=I$.
We see therefore that
{\em
under the conditions of the theorem
the monodromy
transformation ${\mathcal{M}}={{{\mathcal{M}}}^{1/2}}\circ{{{\mathcal{M}}}^{1/2}}$
of {\mbox{$pqr\!s$}}-functions reduces to their multiplication by the constant\/}
$ e^{2{\mathrm{i}}{\ell}\pi} $.
Accordingly, the products of {\mbox{$pqr\!s$}}-functions times 
$z^{-{\ell}}$ reveal the {\em trivial\/} (identical) monodromy transformations.
Thus they can be continuously extended in both directions across the cut
along $\mathbb{R}_{\le0}  $ 
distinguishing ${\mbox{$\mathstrut^\backprime{}\mathbb{C}^*$}}$ from $\mathbb{C}^*$.
Since they also obey a system of
the {\em first order\/}
ODEs (which can easily be derived from {{Eq.s}}~\eqref{eq:ppp'}-\eqref{eq:sss'})
no branching appear showing that %
they are actually single-valued holomorphic in $\mathbb{C}^*$ itself.
The theorem is proved.
\hfill$\square$

\section{The first integral}\label{s:020}

It proves sometimes to be 
useful to take into account the following noteworthy
property of all solutions to {{Eq.s}}{~}\eqref{eq:ppp'}-\eqref{eq:sss'}.
    \begin{theorem}\label{t:040}
Let $\lambda+\mu^2\not=0$. Then the following statements
hold true.
\begin{enumerate}
               \item[0.]
If holomorphic
 functions ${p},{q},{r},{s}  $ obey
{{Eq.s}}{~}\eqref{eq:ppp'}-\eqref{eq:sss'}
then
the value of expression
\begin{equation}
\label{delTa}
{\hspace{0.08ex}{\mathfrak D}}=z^{2(1-{\ell})}\big({p}(z){s}(z)-{q}(z){r}(z)\big)
\end{equation}
does not depends on $z$;
                 \item
if holomorphic functions ${p},{q},{r},{s}  $ obey
{{Eq.s}}~\eqref{eq:pppA}-\eqref{eq:sssA}
then
\begin{equation}
\label{AdelTa}
\raisebox{-0.9ex}{$\Bigl\lfloor$}
\raisebox{-1.3ex}[0ex][2ex]{$\mathstrut$}_{  z\leftleftharpoons -1/z }
\hspace{-2.85em}
{\hspace{0.08ex}{\mathfrak D}}
\hspace{0.2ex}
\rfloor
=
\lfloor\hspace{0.2ex}
{\hspace{0.08ex}{\mathfrak D}}\hspace{0.2ex}
\rfloor,
\end{equation}
where and in what follows $ \lfloor {\hspace{0.08ex}{\mathfrak D}}\rfloor $ denotes {
the right-hand side of {{Eq.}}~\eqref{delTa}} considered as a function of $z$;
                    \item
if holomorphic functions ${p},{q},{r},{s}  $ obey
{{Eq.s}}~\eqref{eq:pppC}-\eqref{eq:sssC}
then
\begin{equation}
\label{CdelTa}
\raisebox{-0.9ex}{$\Bigl\lfloor$}
\raisebox{-1.3ex}[0ex][2ex]{$\mathstrut$}_{  z\leftleftharpoons 1/z }
\hspace{-2.25em}
{\hspace{0.08ex}{\mathfrak D}}
\hspace{0.2ex}
\rfloor
=
\lfloor\hspace{0.2ex}
{\hspace{0.08ex}{\mathfrak D}}
\hspace{0.2ex}
\rfloor;
\end{equation}
\item
if holomorphic functions ${p},{q},{r},{s}  $ obey
{{Eq.s}}~\eqref{eq:pppB}-\eqref{eq:sssB}
then
\begin{equation}
\label{BdelTa}
\raisebox{-0.9ex}{$\Bigl\lfloor$}
\raisebox{-1.3ex}[0ex][2ex]{$\mathstrut$}_{  z\leftleftharpoons {{{\mathcal{M}}}^{1/2}} z }
\hspace{-3.4em}
{\hspace{0.08ex}{\mathfrak D}}
\hspace{0.2ex}
\rfloor
=
\lfloor
\hspace{0.2ex}
{\hspace{0.08ex}{\mathfrak D}}
\hspace{0.2ex}
\rfloor.
\end{equation}
\end{enumerate}
      \end{theorem}
\noindent
It has to be added that the precise meaning of the argument replacements
$ z\leftleftharpoons -1/z   $,
$ z\leftleftharpoons 1/z $, and
$ z\leftleftharpoons {{{\mathcal{M}}}^{1/2}} z \,(\simeq -z)$
involved in the above formulas is the same as
in the corresponding systems of the equations
claimed to be fulfilled.
%
{\sl Proof}.
We shall consider
the
 above
assertions one by one.
~\phantom{.}
\\ 
Assertion {\it 0}. 
Let us expand the expression of  
the derivative of the right-hand side of 
{{Eq.}}~\eqref{delTa}
in case of arbitrary
holomorphic functions ${p},{q},{r},{s}$.
A straightforward computation establishes 
the following 
identity
\begin{equation}
              \label{eq:290}
    \begin{aligned}
z^{2{\ell}}
    {d\over d\,z}  \lfloor {\hspace{0.08ex}{\mathfrak D}} \rfloor  
    \equiv\;
     {s}\,\,
\!\Delta_{{p}}
    -z^2{r}\,\,
\!\Delta_{{q}}
    -{q}\,\,
\!\Delta_{{r}}
    +{p}\,\,
\!\Delta_{{s}}.
    \end{aligned}
\end{equation}
Here 
the symbols
$
\;
\!\Delta{{{\mbox{\tiny\ding{74}}}}}
$, where {{\mbox{\small\ding{74}}}}$\;\in\{{p},{q},{r},{s}\}$,
denote 
the differences of the left- and right-hand sides of {{Eq.s}}{~}%
\eqref{eq:ppp'},
\eqref{eq:qqq'},
\eqref{eq:rrr'},
\eqref{eq:sss'}, respectively, as they stand.
Hence if the latter equations are fulfilled then the 
derivative \eqref{eq:290}
vanishes and 
the value of
$ \lfloor{\hspace{0.08ex}{\mathfrak D}}\rfloor $
does not depend on $z$.

\medskip
\noindent
Assertion {\it 1}.
Its validity
follows from the equality
%
\begin{equation}
\label{eq:300}
    \begin{aligned}
\raisebox{-0.9ex}{$\Bigl\lfloor$}
\raisebox{-1.3ex}[0ex][2ex]{$\mathstrut$}_{  z\leftleftharpoons -1/z }
\hspace{-2.85em}
{\hspace{0.08ex}{\mathfrak D}}
\hspace{0.2ex}
\rfloor
 =   
& \;
\lfloor
\hspace{0.2ex}
{\hspace{0.08ex}{\mathfrak D}}
\hspace{0.2ex}
\rfloor
+e^{-{\mathrm{i}} {\ell}\pi}
\big\lgroup
e^{-{\mathrm{i}} {\ell}\pi}z^{2({\ell}-1)}
\bigl({s}(-1 /z)\,{\mbox{$\mathstrut^{A}\!\Delta_{p}$}}(z)
-{r}(-1 /z)\,{\mbox{$\mathstrut^{A}\!\Delta_{q}$}}(z)\bigr)
\\[-1.0ex]
      & \hphantom { \;\lfloor\,{\hspace{0.08ex}{\mathfrak D}}\;\rfloor   +e^{-{\mathrm{i}} {\ell}\pi}   \big\lgroup      }
-\bigl({r}(z)+\mu\,z^{-2}{p}(z)\bigr)
{\mbox{$\mathstrut^{A}\!\Delta_{r}$}}(z)
-{p}(z)\,{\mbox{$\mathstrut^{A}\!\Delta_{s}$}}(z)
\big\rgroup,
    \end{aligned}
\end{equation}
holding true for arbitrary functions $  {p},{q},{r},{s} $,
holomorphic at (and in the vicinity of) $z={\mathrm{i}}$. 
Here 
the symbols 
${\mbox{$\mathstrut^{A}\!\Delta_{{{\mbox{\tiny\ding{74}}}}}$}}, 
{{\mbox{\small\ding{74}}}}\in\{{p},{q},{r},{s}\} $,
denote the 
differences of the left- and right-hand sides
of
{{Eq.s}}{~}\eqref{eq:pppA},
\eqref{eq:qqqA},
\eqref{eq:rrrA},
\eqref{eq:sssA},
respectively, considered, as they stand, as the functions of $z$.

{{Eq.}}~\eqref{eq:300} 
follows, in turn, from the identity \eqref{eq:Adelta} given
in Appendix \ref{a:050}.

Thus if 
{{Eq.s}}{~}\eqref{eq:pppA}-\eqref{eq:sssA}
are fulfilled then
the equality 
\eqref{AdelTa} holds true. 

\medskip
\noindent
Assertion {\it  2}.
Let us consider the following identity
\begin{equation}
\label{eq:310}
    \begin{aligned}
\raisebox{-0.9ex}{$\Bigl\lfloor$}
\raisebox{-1.3ex}[0ex][2ex]{$\mathstrut$}_{  z\leftleftharpoons 1/z }
\hspace{-2.25em}
{\hspace{0.08ex}{\mathfrak D}}
\hspace{0.2ex}
\rfloor
\equiv&\;
\lfloor\hspace{0.2ex}
{\hspace{0.08ex}{\mathfrak D}}
\hspace{0.2ex}
\rfloor
+z^{2({\ell}-1)}
\bigl({s}(1/z){\mbox{$\mathstrut^{C\!}\!\Delta_{p}$}}(z)
-
{r}(1/z){\mbox{$\mathstrut^{C\!}\!\Delta_{q}$}}(z)\bigr)
\\[-1ex]  &     \hphantom{   \lfloor\,{\hspace{0.08ex}{\mathfrak D}}\;\rfloor }
 +(\lambda+\mu^2)^{-1}
\bigl(
(\lambda\,{r}(z) - \mu\,z^{-2}{s}(z)){\mbox{$\mathstrut^{C\!}\!\Delta_{r}$}}(z)
\\&  \hphantom{ \;  \lfloor\,{\hspace{0.08ex}{\mathfrak D}}\;\rfloor  +  (\lambda+\mu^2)^{-1} \,  }
-
(\mu\,z^2{r}(z) + {s}(z)){\mbox{$\mathstrut^{C\!}\!\Delta_{s}$}}(z)
\bigr)
\end{aligned}
\end{equation}
which is verifiable by straightforward computation.
Here 
${\mbox{$\mathstrut^{C\!}\!\Delta_{{{\mbox{\tiny\ding{74}}}}}$}}(z)$,
{{\mbox{\small\ding{74}}}}$\;\in\{{p},{q},{r},{s}\}$,
denote 
the differences of the left- and right-hand sides
of {{Eq.s}}{~}\eqref{eq:pppC},\eqref{eq:qqqC},\eqref{eq:rrrC},\eqref{eq:sssC},
respectively, considered, as they stand, as the functions of $z$.
The equality \eqref{eq:310}
holds true for arbitrary functions  ${p},{q},{r},{s}$
holomorphic at (and in the vicinity of) $z=1$.
It is extended to any other $z\not=0$ by means of analytic 
continuation.

In view of \eqref{eq:310},
it  is  
obvious that if 
{{Eq.s}}{~}\eqref{eq:pppC}--\eqref{eq:sssC}
are fulfilled then 
{{Eq.}}~\eqref{CdelTa} holds true.

Turning to the assertion {\it 3}, let us consider the equation %
%
%
%
%
\begin{equation}
  \label{eq:320}
    \begin{aligned}
\raisebox{-0.9ex}{$\Bigl\lfloor$}
\raisebox{-1.3ex}[0ex][2ex]{$\mathstrut$}_{  z\leftleftharpoons {{{\mathcal{M}}}^{1/2}} z }
\hspace{-3.4em}
{\hspace{0.08ex}{\mathfrak D}}
\hspace{0.2ex}%
\rfloor
 = 
& \;
\lfloor
\hspace{0.2ex}%
{\hspace{0.08ex}{\mathfrak D}}
\hspace{0.2ex}%
\rfloor
+
e^{-{\mathrm{i}}{\ell}\pi}
z^{2(1-{\ell})}
\big\lgroup
e^{-{\mathrm{i}}{\ell}\pi}
(\, 
{\hspace{-0.2ex}\raisebox{0.8ex}{\!$\raisebox{0.5ex}[0.0ex][0.0ex]%
{$ \rotatebox{185}{\scalebox{0.7}[0.9]{$\curvearrowleft$}} $}
\atop\mbox{$ {s} $}
$}\hspace{-0.2ex}}
\,{\mbox{$\mathstrut^{B\!}\!\Delta_{p}$}} %
-
{\hspace{-0.2ex}\raisebox{0.8ex}{\!$\raisebox{0.5ex}[0.0ex][0.0ex]%
{$ \rotatebox{185}{\scalebox{0.7}[0.9]{$\curvearrowleft$}} $}
\atop\mbox{$ {r} $}
$}\hspace{-0.2ex}}
\,{\mbox{$\mathstrut^{B\!}\!\Delta_{q}$}} %
)
+(\mu\,z^2{p} + {q})
\, {\mbox{$\mathstrut^{B\!}\!\Delta_{r}$}} %
\\[-1.2ex]
&  \hphantom{ \;\lfloor\,{\hspace{0.08ex}{\mathfrak D}}\;\rfloor + e^{-{\mathrm{i}}{\ell}\pi}z^{2(1-{\ell})}\big\lgroup  \! }
+
(\lambda+\mu^2)^{-1}
(\mu\,z^2{r} + {s})
(\mu\,z^2 
\, {\mbox{$\mathstrut^{B\!}\!\Delta_{r}$}} %
+
 {\mbox{$\mathstrut^{B\!}\!\Delta_{s}$}} %
 )
 \big\rgroup.
    \end{aligned}
\end{equation}
Here the symbols 
$ {\mbox{$\mathstrut^{B\!}\!\Delta_{{{\mbox{\tiny\ding{74}}}}}$}}$,
where  {{\mbox{\small\ding{74}}}}$\;\in\{{p},{q},{r},{s}\}  $,
denote the differences
of the left- and right-hand sizes of 
the equations \eqref{eq:pppB}, \eqref{eq:qqqB}, \eqref{eq:rrrB}, \eqref{eq:sssB}, 
respectively. `The diacritic mark'
\raisebox{2ex}{\rotatebox{185}{\scalebox{1}[1.1]{$\curvearrowleft$}}}
denotes the transformation
of the function argument 
defined as follows:
$
{\hspace{-0.2ex}\raisebox{0.8ex}{\!$\raisebox{0.7ex}[0.0ex][0.0ex]%
{$ \rotatebox{185}{\scalebox{0.7}[0.9]{$\curvearrowleft$}} $}
\atop\mbox{$ {{\mbox{\small\ding{74}}}} $}
$}\hspace{-0.2ex}}
= 
{\hspace{-0.2ex}\raisebox{0.8ex}{\!$\raisebox{0.7ex}[0.0ex][0.0ex]%
{$ \rotatebox{185}{\scalebox{0.7}[0.9]{$\curvearrowleft$}} $}
\atop\mbox{$ {{\mbox{\small\ding{74}}}} $}
$}\hspace{-0.2ex}}
(z)=
\lim_{{{\epsilon}}\nearrow 1}{{{\mbox{\small\ding{74}}}}}(e^{{\mathrm{i}}{{\epsilon}}\pi}z)  $.
Here $ \lim_{{{\epsilon}}\nearrow 1} $ should be understood  as
the analytic continuation 
along the image of the segment $[0,1]\ni{{\epsilon}}$ to the  end point
corresponding to ${{\epsilon}}=1$.
In Theorem \ref{t:030} such a transformation is associated with the operator ${{{\mathcal{M}}}^{1/2}}$.

{{Eq.}}{~}\eqref{eq:320}
follows from the identity \eqref{eq:Bdelta} given in Appendix \ref{a:050}.
In turn,
under conditions of the theorem, 
{{Eq.}}~\eqref{BdelTa} is the obvious consequence of {{Eq.}}{~}\eqref{eq:320}.
\hfill$\square$

\medskip

The constant ${\hspace{0.08ex}{\mathfrak D}}$ (in fact, the first integral for the
system \eqref{eq:ppp'}-\eqref{eq:sss'})
may vanish.
Indeed, if it is null at some point
(this is a quadratic constraint to values of {\mbox{$pqr\!s$}}-functions thereat)
then it is zero everywhere.
Such a case bears
many
signs of a degeneracy
--- being nevertheless 
in no way meaningless. 
Following 
here the requirement of genericity,
we assume throughout that  ${\hspace{0.08ex}{\mathfrak D}}\not=0$
without separate mentioning.
It is also worth noting that there is another
case for which many of the relations
discussed here
 degenerate.
Namely, this takes place if $\lambda+\mu^2=0$.
We evade 
here
clarification of its specialties
as well.

\section{On applications of {\mbox{$pqr\!s$}}-functions}\label{s:030}

The properties of {\mbox{$pqr\!s$}}-functions established above
make them an object of notable interest in itself.
However, they arose originally in
the context of 
another important problem, namely,
the study of
symmetries of the space of solutions to
the following ordinary second-order
linear homogeneous differential equation
\begin{equation}
\label{sDCHE}
z^2 {{E''}}
+\big(({\ell}+1)
z
+ \mu\, (1-z^2) \big) {{E'}}
+\big
( -({\ell}+1)\mu\, z
+
\lambda
\big
)
{{E}}
=0.
\end{equation}
Here ${{E}}={{E}}(z)$ is
the unknown holomorphic function,
 ${\ell}, \lambda, \mu $
are the
constant parameters
which
may be identified
with the ones 
involved in {{Eq.s}}{~}\eqref{eq:ppp'}-\eqref{eq:sss'}.

{{Eq.}}~\eqref{sDCHE}
belongs to the family of
double confluent
Heun equations ({{\sf DCHE}}{}).
They 
are
discussed 
in Ref.s~\cite{SW,SL},
 Ref.s~\cite{HeunProject,H}
contain some more
recent
bibliography.
Since a
 generic {{\sf DCHE}}{}
is 
characterized
by the
\textit{four\/}
constant parameters,
whereas {{Eq.}}~\eqref{sDCHE} involves
only \textit{three} ones,
{{Eq.}}~\eqref{sDCHE}
was named 
a
\textit{special\/}
double confluent Heun equation (s{{\sf DCHE}}{}). 
This naming 
is 
adopted 
in the present work as well.

It should be noted that
{{Eq.}}{~}\eqref{sDCHE} was
segregated
within the
{{\sf DCHE}}{} family because of
its intimate
relation (in fact, equivalence)
to
the following
nonlinear first-order ODE%
                               \footnote{
The basic points of such a relationship
are discussed in Ref.~\cite{T2}.
                               } 
\begin{equation*}
\dot\varphi +\sin\varphi=B+A\cos\omega t,
\end{equation*}
in which $\varphi=\varphi(t) $ is the unknown function, 
the symbols
$A,B,\omega$ denote real constants, $t$ is a free real
variable, and
the %
overdot %
denotes %
derivating
with respect to $t$.
The latter
equation 
and its generalizations
are, in turn,
well known
due to their 
emerging 
in a number of problems in physics
(most notably in the modeling of
Josephson junctions) \cite{St,McC,Ba},
mechanics \cite{Foo,FLT}, dynamical systems theory \cite{GI},
and geometry \cite{BLPT}.

In earlier
investigations the functions, obeying equations equivalent
to  {{Eq.s}}{~}\eqref{eq:ppp'}-\eqref{eq:sss'},
were utilized for the constructing of linear
operators sending 
the
space of solutions to {{Eq.}}{~}\eqref{sDCHE}
into itself \cite{BT1}.
It was found that the 
transformations
they determine
generate
a group
which can be regarded as a discrete symmetry of the noted space of solutions.
(More precisely, in case of real parameters,
one of the three groups arises depending on their values).

The principal limitation of those considerations was, however, the
restriction
of the parameter ${\ell}$ 
(sometimes called {\it the order\/} of {{Eq.}}{~}\eqref{sDCHE})
to integers only.
The simplification following from this assumption
(the starting point of derivation of the mentioned symmetry
transformations, in fact) is the reducing
of the functions equivalent to our
{\mbox{$pqr\!s$}}-functions to polynomials%
              \footnote{
The positive integer order
${\ell}\in \mathbb{Z}_{+}$
determines their  degrees %
linearly increasing with its incrementing.
As to the case of negative ${\ell}\in \mathbb{Z}_{-}$,
there is a trick enabling one 
to convert it to the case of the positive order
equal to
 $|{\ell}|$.
                          }
 in $z$ as well as in the parameters
$\lambda, \mu  $. Moreover, there exists the 
recurrent scheme enabling one to compute these polynomials
for any given positive 
 integer ${\ell}$.%

The definition of {\mbox{$pqr\!s$}}-functions considered in the present work
needs 
no such a restriction 
that enables us 
to make
a crucial
step
in 
revealing
of discrete symmetries of the noted space of solutions in
case on non-integer ${\ell}$.
We apply the approach closely following
the one utilized in the case of integer order
although
some specific
subtleties still have to be taken
into account.

To that end, let us define the  two families
of linear operators,
$ 
{\raisebox{1.0ex}{\strut}\raisebox{0.81ex}[0ex][0ex]{${\mbox{\scriptsize
$ {\epsilon} $%
\;\;\;\;}\atop\mbox{$\mathrm{L}_{ {A} }$}}$}}
$  
and
$ 
{\raisebox{1.0ex}{\strut}\raisebox{0.81ex}[0ex][0ex]{${\mbox{\scriptsize
$ {\epsilon} $%
\;\;\;\;}\atop\mbox{$\mathrm{L}_{ {B} }$}}$}}%
$,
 depending on the real parameter ${{\epsilon}}\in[-1,1]$.
They act to arbitrary functions (denoted $E$)
holomorphic in
{{\mbox{$\mathstrut^\backprime{}\mathbb{C}^*$}}
in accordance with the following formulas.
\begin{eqnarray}
\label{eq:eopA}
&&\begin{aligned}
 {\raisebox{1.0ex}{\strut}\raisebox{0.81ex}[0ex][0ex]{${\mbox{\scriptsize
$ {\epsilon} $%
\;\;\;\;}\atop\mbox{$\mathrm{L}_{ {A} }$}}$}}%
 [{{E}}\hspace{0.3ex}](z)
=&\;
e^{\mu(z+1/z)}
\raisebox{-0.9ex}{$\Bigl\lfloor$}
\raisebox{-1.3ex}[0ex][2ex]{$\mathstrut$}_{  z\leftleftharpoons e^{{\mathrm{i}}{{\epsilon}}\pi}/z }
\hspace{-3.3em}
{
z^2{p}(z) {{E'}}(z) + {q}(z){{E}}(z)
}
\,\bigr\rfloor,
\end{aligned}
\\
\label{eq:eopB}
&&\begin{aligned}
 {\raisebox{1.0ex}{\strut}\raisebox{0.81ex}[0ex][0ex]{${\mbox{\scriptsize
$ {\epsilon} $%
\;\;\;\;}\atop\mbox{$\mathrm{L}_{ {B} }$}}$}}%
 [{{E}}\hspace{0.3ex}](z)
=&\;
z^{1 - {\ell}}e^{\mu(z+1/z)}
\raisebox{-0.9ex}{$\Bigl\lfloor$}
\raisebox{-1.3ex}[0ex][2ex]{$\mathstrut$}_{  z\leftleftharpoons e^{{\mathrm{i}}{{\epsilon}}\pi} z }
\hspace{-3.0em}
{
 z^2{r}(z) {{E'}}(z)+ {s}(z){{E}}(z)
}
\hspace{0.02em}\bigr\rfloor.
\end{aligned}
\end{eqnarray}
The functions ${p}, {q}, {r}, {s}  $
are assumed to be holomorphic in the same domain.

If 
${{\epsilon}}=0$
then
the common argument of the functions ${p},{q},{{E}}$ and
${r},{s},{{E}}$ in right-hand sides of \eqref{eq:eopA} and
\eqref{eq:eopB}
coincide with $1/z$ and ${z}$,  respectively.
Then
the values of 
$ 
{\raisebox{1.0ex}{\strut}\raisebox{0.81ex}[0ex][0ex]{${\mbox{\scriptsize
$ {\mbox{\tiny 0}} $%
\;\;\;\;}\atop\mbox{$\mathrm{L}_{ {A} }$}}$}}%
[E\,]$ and 
$
{\raisebox{1.0ex}{\strut}\raisebox{0.81ex}[0ex][0ex]{${\mbox{\scriptsize
$ {\mbox{\tiny 0}} $%
\;\;\;\;}\atop\mbox{$\mathrm{L}_{ {B} }$}}$}}%
[E\,]$
are correctly defined everywhere in {{\mbox{$\mathstrut^\backprime{}\mathbb{C}^*$}}}. 

We introduce now the particular instances  
$
\mathrm{L}_{{A}}
$ and 
$
\mathrm{L}_{{B}}
$
of the operators
 $
 {\raisebox{1.0ex}{\strut}\raisebox{0.81ex}[0ex][0ex]{${\mbox{\scriptsize
$ {\epsilon} $%
\;\;\;\;}\atop\mbox{$\mathrm{L}_{ {A} }$}}$}}%
 $ and  
 $
 {\raisebox{1.0ex}{\strut}\raisebox{0.81ex}[0ex][0ex]{${\mbox{\scriptsize
$ {\epsilon} $%
\;\;\;\;}\atop\mbox{$\mathrm{L}_{ {B} }$}}$}}%
 $
 as the results of the analytic continuations,
 starting from  
 $
 {\raisebox{1.0ex}{\strut}\raisebox{0.81ex}[0ex][0ex]{${\mbox{\scriptsize
$ {\mbox{\tiny 0}} $%
\;\;\;\;}\atop\mbox{$\mathrm{L}_{ {A} }$}}$}}%
 $ and 
 $
 {\raisebox{1.0ex}{\strut}\raisebox{0.81ex}[0ex][0ex]{${\mbox{\scriptsize
$ {\mbox{\tiny 0}} $%
\;\;\;\;}\atop\mbox{$\mathrm{L}_{ {B} }$}}$}}\!%
 $,
along the images of the segment $[0,1]\ni{{\epsilon}}$
through the corresponding maps
$ z \mapsto e^{{\mathrm{i}}{{\epsilon}}}/ z  $
and
$ z \mapsto e^{{\mathrm{i}}{{\epsilon}}} z $.
We may write down 
these relationships
as follows.
\begin{equation}
\label{eq:036}
\mathrm{L}_{{A}}
= \lim_{{{\epsilon}}\nearrow 1} 
{\raisebox{1.0ex}{\strut}\raisebox{0.81ex}[0ex][0ex]{${\mbox{\scriptsize
$ {\epsilon} $%
\;\;\;\;}\atop\mbox{$\mathrm{L}_{ {A} }$}}$}}%
,
\mathrm{L}_{{B}}
= \lim_{{{\epsilon}}\nearrow 1} 
{\raisebox{1.0ex}{\strut}\raisebox{0.81ex}[0ex][0ex]{${\mbox{\scriptsize
$ {\epsilon} $%
\;\;\;\;}\atop\mbox{$\mathrm{L}_{ {B} }$}}$}}%
.
\end{equation}
              \begin{theorem} \label{t:050}
Let 
{\mbox{$pqr\!s$}}-functions
verify {{Eq.s}}{~}\eqref{eq:ppp'}-\eqref{eq:sss'}
and ${{E}}$ verify {{Eq.}}{~}\eqref{sDCHE}.
Then the functions
$
\mathrm{L}_{{A}}
[{{E}}
\hspace{0.3ex}]$ and 
$
\mathrm{L}_{{B}}
[{{E}}
\hspace{0.3ex}]  $
also
verify {{Eq.}}{~}\eqref{sDCHE}.
\end{theorem}
Proof.
Let us  introduce
the operator $ {{\mathcal H}} $ associated with {{Eq.}}{~}\eqref{sDCHE},
i.e.~let
\begin{equation} 
\label{sDCHop}
{{\mathcal H}}[{{E}}\hspace{0.3ex}](z) = z^2{{E''}}(z) + \big(({\ell}+1) z+\mu(1-z^2)\big){{E'}}(z)
+\big(\lambda-\mu({\ell}+1)z\big){{E}}(z).
\end{equation} 
Composing 
it with
the operator
 ${\mathrm L}_A$, the following
expansion of the slightly modified 
result of their combined action to an arbitrary holomorphic
function ${{E}} $ 
can be obtained
    \begin{eqnarray}
    \label{eq:380} 
\raisebox{-0.7ex}{ \scalebox{1.}[0.8]{ $\biggl\lfloor$ }  }\hspace{-1,7ex}
\raisebox{.0ex}[0ex][1.7ex]{$\mathstrut$}_{ z\leftleftharpoons -1/z }
\hspace{-2.65em}
e^{-\mu(z+1/z)}
             ({{\mathcal H}} \circ \!
             \mathrm{L}_{{A}}
             \!)
             [{{E}}\hspace{0.3ex}]{{}}
\big\rfloor
                  &=&
       z^2{p}{{}}\,{{\mathcal H}}'[ {{E}}]{{}}
\nonumber
 \\[-2ex]
 &&\relax    
+
\big(2\mu\,{p}{{}} - {q}{{}} + 2z^2{r}{{}}\big)
   {{\mathcal H}}[ {{E}}]{{}}
\nonumber
\\
&&\relax             
+ z^2{{E'}}{{}}
 \!\Delta_{{p}}
\!\!\mathstrut'{{}}
+ z^2{{E}}{{}}\,
 \!\Delta_{{q}}
\!\!\mathstrut'{{}}
\nonumber
\\
&&
\relax \hspace{-11.2ex}
+\big( (\lambda - ({\ell}+1)\mu{}z){{E}}{{}} + 2 z^2{{E''}}{{}}
+(\mu(2 - z^2) + 2z){{E'}}{{}}
\big)
 \Delta_{{p}}{{}}
\nonumber\\
  &&\relax  
+
\big(
(\mu(1 - z^2)
-z({\ell} - 1 -\mu{}z)){{E}}{{}}
+
z^2{{E'}}{{}}
\big)
\Delta_{{q}}{{}}
\nonumber\\
&&\relax       
+z^2{{E'}}{{}}
  \Delta_{{r}}{{}}
+{{E}}{{}}\,
  \Delta_{{s}}{{}},
\end{eqnarray}
provided the functions ${{E}}, {p}, {q}, {r}, {s}  $
of the variable $z$ are holomorphic at (and in the vicinity of)
$z={\mathrm{i}}$.
Here
the symbols
$
\;\Delta{{{\mbox{\tiny\ding{74}}}}}
{{}}$, where {{\mbox{\small\ding{74}}}}$\;\in\{{p},{q},{r},{s}\},
$
denote 
the differences of the left- and right-hand sides of 
{{Eq.s}}{~}\eqref{eq:ppp'}-\eqref{eq:sss'}
considered as the functions of $z$
(they were already used in the proof of Theorem \ref{t:040}),
${{\mathcal H}}'=d/d z \circ {{\mathcal H}}  $.

In case of 
the operator 
${\mathrm L}_B$
similar expansion looks as follows.
\begin{eqnarray}
\label{eq:390}
\hspace{-2em}
z^{{\ell}-1}e^{-\mu(z+1/z)}
({{\mathcal H}} \circ
\!
\mathrm{L}_{{B}}
\!)
[{{E}}\hspace{0.3ex}]{{}}
  &=&
z^2
\,
{\hspace{-0.2ex}\raisebox{0.8ex}{\!$\raisebox{0.5ex}[0.0ex][0.0ex]%
{$ \rotatebox{185}{\scalebox{0.7}[0.9]{$\curvearrowleft$}} $}
\atop\mbox{$ {r} $}
$}\hspace{-0.2ex}}
\,
{{\hspace{-0.1ex}\raisebox{1.07ex}{\!$\raisebox{0.4ex}[0.0ex][0.0ex]%
 {$\rotatebox{185}{\scalebox{0.7}[0.9]{$\curvearrowleft$}}$}\atop
 \mbox{${\mathcal H}'[\!{{E}}]$}
 \,$}\hspace{-0.63ex}}}
\nonumber\\&&  \hspace{-2ex}
+\left(
2({\ell} - 1)z
\,
{\hspace{-0.2ex}\raisebox{0.8ex}{\!$\raisebox{0.5ex}[0.0ex][0.0ex]%
{$ \rotatebox{185}{\scalebox{0.7}[0.9]{$\curvearrowleft$}} $}
\atop\mbox{$ {r} $}
$}\hspace{-0.2ex}}
\, %
+
\,
{\hspace{-0.2ex}\raisebox{0.8ex}{\!$\raisebox{0.5ex}[0.0ex][0.0ex]%
{$ \rotatebox{185}{\scalebox{0.7}[0.9]{$\curvearrowleft$}} $}
\atop\mbox{$ {s} $}
$}\hspace{-0.2ex}}
\, %
+2z^2
\,
{\hspace{-0.2ex}\raisebox{0.9ex}{\!$ \raisebox{0.8ex}[0.0ex][0.0ex]%
 {$\rotatebox{185}{\scalebox{0.7}[0.9]{$\curvearrowleft$}} $}\atop
 \mbox{${{r}}'$}
 $}\hspace{-0.2ex}}
\,%
\right) 
{{\hspace{-0.1ex}\raisebox{1.07ex}{\!$ \raisebox{0.4ex}[0.0ex][0.0ex]%
 {$\rotatebox{185}{\scalebox{0.7}[0.9]{$\curvearrowleft$}} $} \atop
 \mbox{${\mathcal H}[\!{{E}}]$}\,$}\hspace{-0.63ex}}}
\nonumber
\\
   &&             \hspace{-2ex}
+z^2
{{\hspace{-0.1ex}\raisebox{1.05ex}{\!${\raisebox{0.4ex}[0.0ex][0.0ex]
{$\rotatebox{185}{\scalebox{0.7}[0.9]{$\curvearrowleft$}}$} \atop
 {\mbox{${{E}}$}}\raisebox{0ex}[0ex][0ex]{$\vphantom{{{E}}}'$}
 \,}$}\hspace{-0.09ex}}}
\!\Delta_{{r}}
\!\!\mathstrut' %
+
{{\hspace{-0.1ex}\raisebox{1.07ex}{\!$ \raisebox{0.4ex}[0.0ex][0.0ex]%
 {$ \rotatebox{185}{\scalebox{0.7}[0.9]{$\curvearrowleft$}} $} \atop
 \mbox{${{E}}$}
 \,$}\hspace{-0.08ex}}}
\!\Delta_{{s}}
\!\!\mathstrut' %
-(\lambda+\mu^2)
\bigl(
{{\hspace{-0.1ex}\raisebox{1.05ex}{\!${\raisebox{0.4ex}[0.0ex][0.0ex]
{$\rotatebox{185}{\scalebox{0.7}[0.9]{$\curvearrowleft$}}$} \atop
 {\mbox{${{E}}$}}\raisebox{0ex}[0ex][0ex]{$\vphantom{{{E}}}'$}
 \,}$}\hspace{-0.09ex}}}
\!\Delta_{{p}} %
+
{{\hspace{-0.1ex}\raisebox{1.07ex}{\!$ \raisebox{0.4ex}[0.0ex][0.0ex]%
 {$ \rotatebox{185}{\scalebox{0.7}[0.9]{$\curvearrowleft$}} $} \atop
 \mbox{${{E}}$}
 \,$}\hspace{-0.08ex}}}
\!\Delta_{{q}} %
\bigr)
\nonumber\\
         &&        \hspace{-2ex}
-\big(
(\mu - ({\ell} - 1)z)
{{\hspace{-0.1ex}\raisebox{1.05ex}{\!${\raisebox{0.4ex}[0.0ex][0.0ex]
{$\rotatebox{185}{\scalebox{0.7}[0.9]{$\curvearrowleft$}}$} \atop
 {\mbox{${{E}}$}}\raisebox{0ex}[0ex][0ex]{$\vphantom{{{E}}}'$}
 \,}$}\hspace{-0.09ex}}}
+
(\lambda+({\ell}+1)\mu{}z)
{{\hspace{-0.1ex}\raisebox{1.07ex}{\!$ \raisebox{0.4ex}[0.0ex][0.0ex]%
 {$ \rotatebox{185}{\scalebox{0.7}[0.9]{$\curvearrowleft$}} $} \atop
 \mbox{${{E}}$}
 \,$}\hspace{-0.08ex}}}
\!
\big)
\,\,\!\Delta_{{r}} %
\nonumber\\ 
 &&                \hspace{-2ex}
+ (
{{\hspace{-0.1ex}\raisebox{1.05ex}{\!${\raisebox{0.4ex}[0.0ex][0.0ex]
{$\rotatebox{185}{\scalebox{0.7}[0.9]{$\curvearrowleft$}}$} \atop
 {\mbox{${{E}}$}}\raisebox{0ex}[0ex][0ex]{$\vphantom{{{E}}}'$}
 \,}$}\hspace{-0.09ex}}}
-\mu
{{\hspace{-0.1ex}\raisebox{1.07ex}{\!$ \raisebox{0.4ex}[0.0ex][0.0ex]%
 {$ \rotatebox{185}{\scalebox{0.7}[0.9]{$\curvearrowleft$}} $} \atop
 \mbox{${{E}}$}
 \,$}\hspace{-0.08ex}}}
  \!  )  %
\,\, \!\Delta_{{q}}
. %
\end{eqnarray}
The symbols $
\,\Delta{{{\mbox{\tiny\ding{74}}}}}$
have the same meaning as in {{Eq.}}{~}\eqref{eq:380}.
`The diacritic mark'
\raisebox{2ex}{\rotatebox{185}{\scalebox{1}[1.1]{$\curvearrowleft$}}}
denoting the semi-monodromy transformation
was also used in the proof of Theorem \ref{t:040}.
 It is worth  reminding  that
$
{\hspace{-0.2ex}\raisebox{0.8ex}{\!$\raisebox{0.7ex}[0.0ex][0.0ex]%
{$ \rotatebox{185}{\scalebox{0.7}[0.9]{$\curvearrowleft$}} $}
\atop\mbox{$ {{\mbox{\small\ding{74}}}} $}
$}\hspace{-0.2ex}}
=
{\hspace{-0.2ex}\raisebox{0.8ex}{\!$\raisebox{0.7ex}[0.0ex][0.0ex]%
{$ \rotatebox{185}{\scalebox{0.7}[0.9]{$\curvearrowleft$}} $}
\atop\mbox{$ {{\mbox{\small\ding{74}}}} $}
$}\hspace{-0.2ex}}
(z)
=
\lim_{{{\epsilon}}\nearrow 1}{{{\mbox{\small\ding{74}}}}}(e^{{\mathrm{i}}{{\epsilon}}\pi}z)$ for any holomorphic
function ${{\mbox{\small\ding{74}}}}$.
If the functions ${{E}}, {p}, {q}, {r}, {s}$
are holomorphic on {{\mbox{$\mathstrut^\backprime{}\mathbb{C}^*$}}}
and $\Im z<0$ then the argument of evaluation of
`semi-monodromy-transformed'
functions also belongs to {{\mbox{$\mathstrut^\backprime{}\mathbb{C}^*$}}}
and all the constituents of {{Eq.}}{~}\eqref{eq:390} are well defined.
Other values of $z$ are to be handled by means of the analytic continuation.

The equalities \eqref{eq:380}, \eqref{eq:390}
follow from the identities \eqref{eq:HoL_A} and
                           \eqref{eq:HoL_B}, respectively,
given in Appendix \ref{a:060}.
In turn, the theorem's assertion follows 
from {{Eq.s}}{~}\eqref{eq:380} and \eqref{eq:390}
since  the fulfillment of
{{Eq.s}}{~}\eqref{eq:ppp'}-\eqref{eq:sss'}
implies 
$\;
\Delta{{{\mbox{\tiny\ding{74}}}}}
=0$
and the identical vanishing of ${{\mathcal H}}[ {{E}}]$
is equivalent to fulfillment of {{Eq.}}{~}\eqref{sDCHE}
that had also been assumed.
\hfill$\square$

\begin{remark}
\hangindent=2ex
\rm
The transformations
realized by the operators
$\mathrm{L}_{{A}}$
and
$\mathrm{L}_{{B}}$
carry out the
(lifted) replacements
$  z\leftleftharpoons -1/z  $ and $  z\leftleftharpoons -z $
of arguments of the functions involved.
There exists the third operator which we denote
$\mathrm{L}_{{C}}$
also sending any solution
to {{Eq.}}{~}\eqref{sDCHE} to some its solution
and utilizing the missed replacement
$  z\leftleftharpoons 1/z$
of arguments
expressing the composition of the preceding ones
and constituting in conjunction with them
the Klein group of maps naturally acting on $\mathbb{C}^*$.
$\mathrm{L}_{{C}} $
is not linked to 
{\mbox{$pqr\!s$}}-functions and is well defined for any choice of constant parameters.
It
can be represented 
by the following formula. 
\begin{equation}
 \label{eq:400}
\mathrm{L}_{{C}}
[{{E}}\hspace{0.3ex}](z)=
z^{-{\ell}-1}
\raisebox{-0.9ex}{$\Bigl\lfloor$}
\raisebox{-1.2ex}[0ex][2ex]{$\mathstrut$}_{  z\leftleftharpoons 1/z }
\hspace{-2.5em}
 {{E'}}(z) - \mu{{E}}(z)
\hspace{0.02em}\rfloor.
\end{equation}
In view of the nontrivial structure of the domain of solutions to
{{Eq.}}{~}\eqref{sDCHE}
``the implementation'' of the rule \eqref{eq:400} is not unique.
In particular, for one of them  (the lift of) $+1$
is the fixed point
of the map indicated
by the argument replacement $  z \leftleftharpoons 1/z $
whereas  for the other one it is
(the lift of) $-1$ which plays a similar role.
\end{remark}

The transformations 
of the space of solutions to  {{Eq.}}{~}\eqref{sDCHE}
associated with {\mbox{$pqr\!s$}}-functions
possess the properties of quasi-involutions
similar to ones found earlier in the case of integer order 
${\ell}$, {{\it cf\/}} Ref. \cite{T3}, {{Eq.s}}{~}(34), (35).
                       \begin{theorem} \label{t:060}
Let $\lambda+\mu^2\not=0$,
the function
${{E}}$ obey {{Eq.}}{~}\eqref{sDCHE},
and 
{\mbox{$pqr\!s$}}-functions obey {{Eq.s}}{~}\eqref{eq:ppp'}-\eqref{eq:sss'}.
If, additionally, 
\begin{enumerate}
\item {{Eq.s}}{~}\eqref{eq:pppA}-\eqref{eq:sssA}
hold true then
\begin{equation}
\label{eq:410}
({\mathrm L}_A\circ{\mathrm L}_A)[{{E}}\hspace{0.3ex}]
= - e^{{\mathrm{i}}{\ell}\pi}{\hspace{0.08ex}{\mathfrak D}}\cdot{{E}};
\end{equation}
\item
{{Eq.s}}{~}\eqref{eq:pppB}-\eqref{eq:sssB}
hold true then
\begin{equation}
\label{eq:420}
({\mathrm L}_B\circ{\mathrm L}_B)[{{E}}\hspace{0.3ex}]
= -(\lambda+\mu^2) e^{2{\mathrm{i}}{\ell}\pi}{\hspace{0.08ex}{\mathfrak D}}\cdot{\mathcal{M}}[{{E}}\hspace{0.3ex}].
\end{equation}
\end{enumerate}
                         \end{theorem}
Proof.
The above claims 
follow from the equalities 
%
%
%
%
\begin{eqnarray}
\label{eq:430}\hspace{-4ex}
(
 {{\mathrm L}_A } 
\circ
 {{\mathrm L}_A } 
)[  {{E}} ](z)
+e^{{\mathrm{i}}{\ell}\pi}
\lfloor{\hspace{0.08ex}{\mathfrak D}}\rfloor
{{E}} (z)
\hspace{-1ex}
&=&
\hspace{-1ex}
\big(({s}(z)+\mu{\,}z^{-2}{q}(z)) {{E}}(z)
\nonumber
\\ && \hphantom{ +\big(    }
+(z^2{r}(z)+\mu\,{p}(z)) {{E'}}(z)
   \big)%
          {\mbox{$\mathstrut^{A}\!\Delta_{p}$}}(z)
\nonumber \\ &&
+\big({q}(z){{E}}(z)+z^2{p}(z){{E'}}(z)\big)
           {\mbox{$\mathstrut^{A}\!\Delta_{q}$}}(z)
\nonumber \\ &&
+
{\hspace{-0.2ex}\raisebox{0.8ex}{\!$\raisebox{0.5ex}[0.0ex][0.0ex]%
{$ \rotatebox{185}{\scalebox{0.7}[0.9]{$\curvearrowleft$}} $}
\atop\mbox{$ {p} $}
$}\hspace{-0.2ex}}
(z^{-1})
\bigl({p}(z)
{{\mathcal H}}[{{E}}\hspace{0.3ex}](z)
\nonumber \\ && \hphantom{    \hspace{2.2ex}  
                          (z^{-1})\bigl( }
+  %
{{E'}}(z)
\Delta_{{p}}(z)
+ 
{{E}}(z)
 \Delta_{{q}}(z)
\bigr),
\end{eqnarray}
\vspace{-6ex}
\begin{eqnarray}
\label{eq:440} \hspace{-5ex}
e^{{\mathrm{i}}{\ell}\pi}
 z^{2({\ell}-1)}
\big\lgroup
(
{\mathrm L}_B
\circ
{\mathrm L}_B
)[  {{E}}\, ]
+(\lambda+\mu^2) e^{2{\mathrm{i}}{\ell}\pi}\,
{\mbox{$\lfloor
{\hspace{0.05em}\raisebox{1.35ex}{\!$\raisebox{0.4ex}[0.0ex][0.0ex]%
 {$ \rotatebox{185}{\scalebox{0.7}[0.9]{$\curvearrowleft$}} $} \atop
 \raisebox{-0.5ex}[0.0ex][0.0ex]{$  {\hspace{0.08ex}{\mathfrak D}}  $}
$}\hspace{-0.15em}}
\rfloor$}}
{ \raisebox{0.9ex}{$\mbox{\tiny\;$\circlearrowleft$}\atop\mbox{${{E}}$}$} }
\big\rgroup \!
                              &=&
\nonumber
\\[-0.0ex]
&&  \hspace{-14ex}
-(\lambda+\mu^2)\,
{\hspace{-0.2ex}\raisebox{0.8ex}{\!$\raisebox{0.5ex}[0.0ex][0.0ex]%
{$ \rotatebox{185}{\scalebox{0.7}[0.9]{$\curvearrowleft$}} $}
\atop\mbox{$ {r} $}
$}\hspace{-0.2ex}}
\cdot
\bigl(
z^2\! 
{\mbox{$\raisebox{0.9ex}{$\mbox{\tiny\;$\circlearrowleft$}\atop\mbox{${{E}}$}$}\mathstrut'$}}
\,\,
{\hspace{-0.2em}\raisebox{1.3ex}{\!$\raisebox{0.4ex}[0.0ex][0.0ex]%
 {$\rotatebox{185}{\scalebox{0.7}[0.9]{$\curvearrowleft$}} $} \atop
 \raisebox{-0.5ex}[0.0ex][0.0ex]{$
 {\mbox{$\mathstrut^{{B}\!}\!\Delta_{ {p} }$}}
 $}$}\hspace{-0.2em}}
+
 { \raisebox{0.9ex}{$\mbox{\tiny\;$\circlearrowleft$}\atop\mbox{${{E}}$}$} }
 \,\,
 {\hspace{-0.2em}\raisebox{1.3ex}{\!$\raisebox{0.4ex}[0.0ex][0.0ex]%
 {$\rotatebox{185}{\scalebox{0.7}[0.9]{$\curvearrowleft$}} $} \atop
 \raisebox{-0.5ex}[0.0ex][0.0ex]{$
 {\mbox{$\mathstrut^{{B}\!}\!\Delta_{ {q} }$}}
 $}$}\hspace{-0.2em}}
\bigr)
\nonumber\\ 
&&  \hspace{-14ex}
-(\mu\,
z^2\, %
{\hspace{-0.2ex}\raisebox{0.8ex}{\!$\raisebox{0.5ex}[0.0ex][0.0ex]%
{$ \rotatebox{185}{\scalebox{0.7}[0.9]{$\curvearrowleft$}} $}
\atop\mbox{$ {r} $}
$}\hspace{-0.2ex}}
\!+\!
{\hspace{-0.2ex}\raisebox{0.8ex}{\!$\raisebox{0.5ex}[0.0ex][0.0ex]%
{$ \rotatebox{185}{\scalebox{0.7}[0.9]{$\curvearrowleft$}} $}
\atop\mbox{$ {s} $}
$}\hspace{-0.2ex}}
 )
\bigl(
z^2\! 
{\mbox{$\raisebox{0.9ex}{$\mbox{\tiny\;$\circlearrowleft$}\atop\mbox{${{E}}$}$}\mathstrut'$}}
\,\,
{\hspace{-0.2em}\raisebox{1.3ex}{\!$\raisebox{0.4ex}[0.0ex][0.0ex]%
 {$\rotatebox{185}{\scalebox{0.7}[0.9]{$\curvearrowleft$}} $} \atop
 \raisebox{-0.5ex}[0.0ex][0.0ex]{$
 {\mbox{$\mathstrut^{{B}\!}\!\Delta_{ {r} }$}}
 $}$}\hspace{-0.2em}}
+
 { \raisebox{0.9ex}{$\mbox{\tiny\;$\circlearrowleft$}\atop\mbox{${{E}}$}$} }
 \,\,
 {\hspace{-0.2em}\raisebox{1.3ex}{\!$\raisebox{0.4ex}[0.0ex][0.0ex]%
 {$\rotatebox{185}{\scalebox{0.7}[0.9]{$\curvearrowleft$}} $} \atop
 \raisebox{-0.5ex}[0.0ex][0.0ex]{$
 {\mbox{$\mathstrut^{{B}\!}\!\Delta_{ {s} }$}}
 $}$}\hspace{-0.2em}}
\bigr)
\nonumber\\[0.2ex]
&&\hspace{-14ex}
+\,
{\hspace{-0.2ex}\raisebox{0.8ex}{\!$\raisebox{0.5ex}[0.0ex][0.0ex]%
{$ \rotatebox{185}{\scalebox{0.7}[0.9]{$\curvearrowleft$}} $}
\atop\mbox{$ {r} $}
$}\hspace{-0.2ex}}
\cdot\big(\!
  z^2
{ \raisebox{0.805ex}{$\mbox{\tiny\;$\circlearrowleft$}\atop\mbox{$ {r} $}$} }
{ \raisebox{0.95ex}{$\mbox{\tiny\;$\circlearrowleft$}\atop\mbox{${{\mathcal H}}$}$} }
+
  z^2
{\mbox{$\raisebox{0.9ex}{$\mbox{\tiny\;$\circlearrowleft$}\atop\mbox{${{E}}$}$}\mathstrut'$}}
{\mbox{\hspace{0.3em}%
\raisebox{0.8ex}{\!$\raisebox{0.2ex}[0.0ex][0.0ex]%
{\scalebox{0.7}[0.7]{$\hspace{-1.5ex}\circlearrowleft$} }\atop
 \raisebox{0.0ex}[0.0ex][0.0ex]{$
\!\! \Delta_{{r}}
 $}$}\hspace{-0.1em}}}
+
{ \raisebox{0.9ex}{$\mbox{\tiny\;$\circlearrowleft$}\atop\mbox{${{E}}$}$} }
{\mbox{\hspace{0.3em}%
\raisebox{0.8ex}{\!$\raisebox{0.2ex}[0.0ex][0.0ex]%
{\scalebox{0.7}[0.7]{$\hspace{-1.5ex}\circlearrowleft$} }\atop
 \raisebox{0.0ex}[0.0ex][0.0ex]{$
 \!\! \Delta_{{s}}
 $}$} \hspace{-0.9em}
 } }\big),
\end{eqnarray}
involving arbitrary functions ${{E}}, {p}, {q}, {r}, {s}  $
and their derivatives.
These 
are, in turn, the consequences of the identities
\eqref{eq:L_A o L_A} and \eqref{eq:L_B o L_B},
given in Appendix \ref{a:070}.
Concerning the notations utilized therein,
let us remind that $ \lfloor{\hspace{0.08ex}{\mathfrak D}}\rfloor $
denote the right-hand side of {{Eq.}}{~}\eqref{delTa}
considered as a function of $z$.
The symbols
$
\;\Delta{{{\mbox{\tiny\ding{74}}}}}
$,
$
  {\mbox{$\mathstrut^{A}\!\Delta_{{{\mbox{\tiny\ding{74}}}}}$}}
$%
, and
$ {\mbox{$\mathstrut^{B\!}\!\Delta_{{{\mbox{\tiny\ding{74}}}}}$}}
$, %
where  {{\mbox{\small\ding{74}}}}$\;\in\{{p},{q},{r},{s}\}$,
denote the differences of the left- and right-hand sides
for
{{Eq.s}}{~}\eqref{eq:ppp'}-\eqref{eq:sss'},
for {{Eq.s}}{~}\eqref{eq:pppA}-\eqref{eq:sssA},
and for
{{Eq.s}}{~}\eqref{eq:pppB}-\eqref{eq:sssB},
respectively. They are also considered as the functions of $z$.

There are also the two kinds of `diacritic marks' in use.
Of them, `the accent'
$ \raisebox{1.6ex} {\rotatebox{185}{\scalebox{0.99}[0.99]{$\curvearrowleft$}}} $
indicates the transformation of the function argument
carrying out its continuous anti-clockwise
rotation in the complex plane
at an angle $\pi$. It was earlier named the semi-monodromy map.
In Theorem \ref{t:030}
such a transformation is associated with the operator
${{{\mathcal{M}}}^{1/2}}$.
Evidently, 
if $\Im z<0$ then
         $
         {\hspace{-0.2ex}\raisebox{0.8ex}{\!$\raisebox{0.7ex}[0.0ex][0.0ex]%
{$ \rotatebox{185}{\scalebox{0.7}[0.9]{$\curvearrowleft$}} $}
\atop\mbox{$ {{\mbox{\small\ding{74}}}} $}
$}\hspace{-0.2ex}}
         (z)$
is simply ${{\mbox{\small\ding{74}}}}(-z)$.
However, if $\Im z \ge 0$ then
the
semi-monodromy transformation
sends such argument out the subdomain {{\mbox{$\mathstrut^\backprime{}\mathbb{C}^*$}}}
and this can not be expressed by the inversion of the sign.
It worth noting here that 
 $
 {\hspace{-0.2ex}\raisebox{0.8ex}{\!$\raisebox{0.5ex}[0.0ex][0.0ex]%
{$ \rotatebox{185}{\scalebox{0.7}[0.9]{$\curvearrowleft$}} $}
\atop\mbox{$ {p} $}
$}\hspace{-0.2ex}}
 (z^{-1})$
(see the last but one line in
{{Eq.}}~\eqref{eq:430}) is well defined, provided $\Im z>0$. Indeed,
then $\Im z^{-1}<0  $ and the argument of evaluation of the function
${p}$ when computing
$
{\hspace{-0.2ex}\raisebox{0.8ex}{\!$\raisebox{0.5ex}[0.0ex][0.0ex]%
{$ \rotatebox{185}{\scalebox{0.7}[0.9]{$\curvearrowleft$}} $}
\atop\mbox{$ {p} $}
$}\hspace{-0.2ex}}
(z^{-1})= \lim_{{{\epsilon}}\nearrow 1}{p}(e^{{\mathrm{i}}{{\epsilon}}\pi} z^{-1})$
belongs to {{\mbox{$\mathstrut^\backprime{}\mathbb{C}^*$}}}.

The second `accent'
{\raisebox{0.7ex}{\scalebox{0.7}[0.7]{$\circlearrowleft$}}}
has a similar meaning but ``the rotation angle'' of a function argument is here
twice as much  amounting to
$2\pi$. Such a transformation 
looks like
a full revolution in $\mathbb{C}^*  $ around zero 
but it does not lead to 
the identical map 
in view of non-trivial structure
(distinction of %
complex plane or any subset of the complex
plane) of the domains of the functions we consider.
Rather it corresponds to the {\it monodromy transformation}.

For some reasons
we had agreed above to consider
 {\mbox{$pqr\!s$}}-functions on their
subdomain 
$\mathstrut^\backprime\mathbb{C}^* = \mathbb{C}^*{\,\fgebackslash\,}{\mathbb R}_- $%
. 
Here, however, this is not enough and we are forced to introduce for a time
a somehow extended one.
Indeed, if $z\in {{\mbox{$\mathstrut^\backprime{}\mathbb{C}^*$}}} $
then the point of evaluation of
a monodromy-transformed function does not
belong
to
{{\mbox{$\mathstrut^\backprime{}\mathbb{C}^*$}}}
due to the cut along the ray of negative reals which
the circular path
of analytical continuation inevitably meets. 
``The minimally extended subdomain'',
where the monodromy map can still be consistently defined,
is constructed, for instance, 
by means of addition of another copy of
{{\mbox{$\mathstrut^\backprime{}\mathbb{C}^*$}}} 
and
the gluing of it to the original one 
along the opposite edges of their cuts (the two complementary ones remain free).
Then if $z$ belongs to the ``lower'' (original) sheet
of this ``double\,-{{\mbox{$\mathstrut^\backprime{}\mathbb{C}^*$}}}''
then
the point of evaluation of analytic continuation of the
function to be
monodromy 
 transformed
belongs to the upper one and
in this way
all the constituents of {{Eq.}}{~}\eqref{eq:440}
can be consistently 
computed (and it is finally fulfilled).

The important circumstance is, however, that under conditions of the theorem
the evaluation of many functions and the handling
of the associated subtleties it implies
is superfluous.
Indeed,
the
fulfillment of certain equations
required by the theorem conditions
means 
the vanishing of the expressions
$
{{\mathcal H}}[{{E}}\hspace{0.3ex}],\,
\Delta{{{\mbox{\tiny\ding{74}}}}},
 {\mbox{$\mathstrut^{A}\!\Delta_{{{\mbox{\tiny\ding{74}}}}}$}},$ and\, $
 {\mbox{$\mathstrut^{B\!}\!\Delta_{{{\mbox{\tiny\ding{74}}}}}$}}
$
which yield  zero independently of the point of their evaluation.
They constitute a full collection of the factors
in the right-hand sides of {{Eq.s}}~\eqref{eq:430} and \eqref{eq:440}
such that
if all they are zero then
all the terms 
on the right
vanish.
The only function with transformed argument
which
`survives'
is
the monodromy transformed function ${{E}}$ involved in the left-hand side of
{{Eq.}}~\eqref{eq:440} in the form 
$ 
{ \raisebox{0.9ex}{$\mbox{\tiny\;$\circlearrowleft$}\atop\mbox{${{E}}$}$} }
(={\mathcal{M}}[{{E}}\hspace{0.3ex}]) $.
Besides, we know that the fulfillment of
{{Eq.s}}{~}\eqref{eq:ppp'}-\eqref{eq:sss'}
leads to the independence of $ \lfloor{\hspace{0.08ex}{\mathfrak D}}\rfloor $
on the point of evaluation, see Theorem \ref{t:040}.
Thus we may replace it, as well as
\raisebox{0ex}[2.5ex][0ex]{$
{\mbox{$\lfloor
{\hspace{0.05em}\raisebox{1.35ex}{\!$\raisebox{0.4ex}[0.0ex][0.0ex]%
 {$ \rotatebox{185}{\scalebox{0.7}[0.9]{$\curvearrowleft$}} $} \atop
 \raisebox{-0.5ex}[0.0ex][0.0ex]{$  {\hspace{0.08ex}{\mathfrak D}}  $}
$}\hspace{-0.15em}}
\rfloor$}}
$,}
with the
corresponding
constant ${\hspace{0.08ex}{\mathfrak D}}$ (whose value is actually determined
in a complicated way
by the parameters ${\ell}, \lambda, \mu$).

Now, as the right-hand sides of {{Eq.s}}~\eqref{eq:430}, \eqref{eq:440}
are zero,
the vanishing of their left-hand sides leads just to
{{Eq.s}}~\eqref{eq:410} and \eqref{eq:420}.
\hfill$\square$
                         \begin{corollary}\label{c:030}
If the conditions of Theorem \ref{t:060} are fulfilled then 
the operators 
$
\mathrm{L}_{{A}},\mathrm{L}_{{B}}
$
it concerns
determine automorphisms of the space of solutions to {{Eq.}}{~}\eqref{sDCHE}.
                          \end{corollary}
Indeed, due to \eqref{eq:410} and \eqref{eq:420}
they have no zero eigenvalues 
and thus their  kernels are trivial.
                          
\section{Summary}\label{s:040}

We define a family of quads of holomorphic functions
(referred to, for brevity, as {\mbox{$pqr\!s$}}-functions) as the non-trivial solutions to
the system of linear homogeneous first order ODEs
\eqref{eq:ppp'}-\eqref{eq:sss'}.
Each instance of such functions can be constructed
as a solution of the Cauchy problem for the initial data
specified at any given point $z_0$ except zero.
It is shown that, fixing the initial data
at $z_0={\mathrm{i}}$ and claiming fulfillment of the linear homogeneous constraint
\eqref{eq:050}, one obtains {\mbox{$pqr\!s$}}-functions
which obey the equalities
\eqref{eq:pppA}-\eqref{eq:sssA} (Theorem \ref{t:010}).
Similarly,
if the initial data
are specified at $z_0=1$
and obey thereat the two linear homogeneous constraints
\eqref{eq:110}, \eqref{eq:120}
then
{\mbox{$pqr\!s$}}-functions
obey the equalities
\eqref{eq:pppC}-\eqref{eq:sssC} (Theorem \ref{t:020}).
Lastly, if all the three mentioned linear constraints 
(imposed at two distinct locations) 
are met then the equalities
\eqref{eq:pppB}-\eqref{eq:sssB}
involving semi-monodromy map
take place as well (Theorem \ref{t:030}). 
This case is most important since for it
the monodromy transformation can also be easely computed.
It turns out coinciding with multiplication to a known numerical
factor
showing that {\mbox{$pqr\!s$}}-functions
are the products of certain power function
and functions holomorphic on $\mathbb{C}^*$ (instead of the uinversal
cover of $\mathbb{C}^*$).

{\mbox{$P\!qr\!s$}}-functions had found application (if fact, arose)
in frameworks of investigation of
properties of solutions to special double confluent Heun equation
\eqref{sDCHE}. 
Under conditions here assumed  
the operators 
$ \mathrm{L}_{{A}},\mathrm{L}_{{B}} $
defined by the formulas
\eqref{eq:036},
\eqref{eq:eopA},
\eqref{eq:eopB}
turn out to define the  
maps of 
the space of
its solutions 
into itself
(Theorem \ref{t:050}). 
Moreover, they
possess quite remarkable composition properties.
It particular, the operator 
$\mathrm{L}_{{A}}$
is ``almost involutive''
(see Theorem \ref{t:060}, {{Eq.}}~\eqref{eq:410})
while 
$\mathrm{L}_{{B}}$%
,
being applied twice, reduces, up to a known constant factor,
to the monodromy transformation (see {{Eq.}}~\eqref{eq:420}).
Besides, they define automorphisms of the space of solutions to 
{{Eq.}}~\eqref{sDCHE} (Corollary \ref{c:030}).

In the special case
of integer values of the constant parameter ${\ell}$
the functions almost identical to our {\mbox{$pqr\!s$}}-functions
were originally introduced in Ref.~\cite{BT1}.
The distinction of functions with the same notations
considered therein against the present ones reduces to   different
normalizations of the functions ${p}$ and ${q}$.
It is worth mentioning that the variant of
{\mbox{$pqr\!s$}}-functions considered in \cite{BT1}
deals exclusively with polynomials.
Moreover, they are
polynomial not only in $z$
but also in the parameters $\lambda$ and $\mu$ (while ${\ell}$ determines
the polynomial degrees).
Thus we may claim that in the case of a (positive) integer ${\ell}$
{{Eq.s}}~\eqref{eq:ppp'}-\eqref{eq:sss'}
admit a polynomial solution.

\vfill

\appendix

\section{Identities leading to {{Eq.s}}{~}\eqref{eq:210}}\label{a:010}

The equations \eqref{eq:210}
the proof of Theorem \ref{t:010}
leans on 
are the straightforward
consequences of the
four \textit{identities\/}
displayed below
which are, in principle, verifiable
by explicit
computations%
                    \footnote{\label{ft10}
As a matter of fact,
they were handled 
with  help of the computer algebra.
Similar remarks concern
the majority of formulas in the present
paper or, at least, all more or less lengthy ones.
}.
Namely, the following equalities hold true
%
%
%
%
%
%
\begin{eqnarray}
\label{eq:idAp}
z^2
\frac{d}{d z}
{\raisebox{0.95ex}{$\raisebox{-0.2ex}[0ex][0ex]{\scriptsize$\epsilon$}
\atop\raisebox{-0.2ex}[0ex][0ex]{$\mathstrut^{A\hspace{-0.1ex}}\hspace{-0.5ex}\Delta_{p}$}$}}
(z)
\!&\equiv&\!
-z({\ell} - 1 +e^{-{\mathrm{i}}{{\epsilon}}\pi}\mu{}z) 
\,
{\raisebox{0.95ex}{$\raisebox{-0.2ex}[0ex][0ex]{\scriptsize$\epsilon$}
\atop\raisebox{-0.2ex}[0ex][0ex]{$\mathstrut^{A\hspace{-0.1ex}}\hspace{-0.5ex}\Delta_{p}$}$}}
(z)
\nonumber
\\[-1.5ex]&&
+e^{-{\mathrm{i}}{{\epsilon}}\pi} z^2 \,
{\raisebox{0.95ex}{$\raisebox{-0.2ex}[0ex][0ex]{\scriptsize$\epsilon$}
\atop\raisebox{-0.2ex}[0ex][0ex]{$\mathstrut^{A\hspace{-0.1ex}}\hspace{-0.5ex}\Delta_{q}$}$}}
(z)
-e^{{\mathrm{i}}{{\epsilon}}\pi}
{\raisebox{0.95ex}{$\raisebox{-0.2ex}[0ex][0ex]{\scriptsize$\epsilon$}
\atop\raisebox{-0.2ex}[0ex][0ex]{$\mathstrut^{A\hspace{-0.1ex}}\hspace{-0.5ex}\Delta_{r}$}$}}
 (z)
\nonumber\\&&\relax   
- e^{-{\mathrm{i}}{{\epsilon}}\pi} z^2
 \,\!\Delta_{{p}}
( e^{{\mathrm{i}}{{\epsilon}}\pi}\!/z )
 \nonumber
+
e^{{\mathrm{i}}{{\epsilon}}{\ell}\pi} z^{2(1-{\ell})}\,
\!\Delta_{{p}}( z )
\nonumber
\nonumber
  \\
&&
\relax\hspace{-9ex}
\llap{$+$}
(1+e^{-{\mathrm{i}}{{\epsilon}}\pi})\,
e^{{\mathrm{i}}{{\epsilon}}{\ell}\pi}
z^{2(1-{\ell})}\times
\nonumber\\&&
\;\;
\left(
\mu(1+(1-e^{{\mathrm{i}}{{\epsilon}}\pi})z^2){p}(z)
-e^{{\mathrm{i}}{{\epsilon}}\pi}{q}(z)
+z^2{r}(z)
\right),
\end{eqnarray}
%
%
%
%
%
\begin{eqnarray}
\label{eq:idAq}
z^3
\frac{d}{d z}
{\raisebox{0.95ex}{$\raisebox{-0.2ex}[0ex][0ex]{\scriptsize$\epsilon$}
\atop\raisebox{-0.2ex}[0ex][0ex]{$\mathstrut^{A\hspace{-0.1ex}}\hspace{-0.5ex}\Delta_{q}$}$}}
(z)
\!&\equiv&\!
e^{{\mathrm{i}}{{\epsilon}}\pi}
(e^{{\mathrm{i}}{{\epsilon}}\pi}({\ell} + 1)\mu - \lambda{}z)
 \,
 {\raisebox{0.95ex}{$\raisebox{-0.2ex}[0ex][0ex]{\scriptsize$\epsilon$}
\atop\raisebox{-0.2ex}[0ex][0ex]{$\mathstrut^{A\hspace{-0.1ex}}\hspace{-0.5ex}\Delta_{p}$}$}}
 (z)
\nonumber
\\[-1.5ex]&&
-e^{{\mathrm{i}}{{\epsilon}}\pi}\mu{}z
  \,
{\raisebox{0.95ex}{$\raisebox{-0.2ex}[0ex][0ex]{\scriptsize$\epsilon$}
\atop\raisebox{-0.2ex}[0ex][0ex]{$\mathstrut^{A\hspace{-0.1ex}}\hspace{-0.5ex}\Delta_{q}$}$}}
  (z)
-e^{{\mathrm{i}}{{\epsilon}}\pi}z
\,
{\raisebox{0.95ex}{$\raisebox{-0.2ex}[0ex][0ex]{\scriptsize$\epsilon$}
\atop\raisebox{-0.2ex}[0ex][0ex]{$\mathstrut^{A\hspace{-0.1ex}}\hspace{-0.5ex}\Delta_{s}$}$}}
(z)
\nonumber\\&&\relax 
-
e^{{\mathrm{i}}{{\epsilon}}\pi}z
\,\,
\!\Delta_{{q}}   ( e^{{\mathrm{i}}{{\epsilon}}\pi}\!/z )
-e^{{\mathrm{i}}{{\epsilon}}{\ell}\pi}
z^{1-2{\ell}}
\big(
\mu
\,\,\!\Delta_{{p}}(z)
+
z^2
\,\!\Delta_{{r}}(z)
\big)
\nonumber
  \\
  &&\relax  \hspace{-2.5ex}
\llap{$+$}
          (1+e^{{\mathrm{i}}{{\epsilon}}\pi})\,
e^{{\mathrm{i}}{{\epsilon}}{\ell}\pi}
z^{1-2{\ell}}  \times
\nonumber\\&& \relax\hspace{-7ex} \phantom{ (1+e^{{\mathrm{i}}{{\epsilon}}\pi}) }
\left(
       \bigl((\lambda+\mu^2)z^2+(1-e^{{\mathrm{i}}{{\epsilon}}\pi})({\ell}+1)\mu{}z-\mu^2\bigr)
{p}(z)
\right.\nonumber\\&&\relax\hspace{2.9ex} \phantom{ (1+e^{{\mathrm{i}}{{\epsilon}}\pi}) }
\left.
+\mu{q}(z)
+z^2\left(\mu(z^2-1){r}(z)+{s}(z)\right)
\vphantom{      ((\lambda+\mu^2)z^2+(1-e^{{\mathrm{i}}{{\epsilon}}\pi})({\ell}+1)\mu{}z-\mu^2)      }\!
\right),
\end{eqnarray}
%
%
%
%
%
\begin{eqnarray}
\label{eq:idAr}
z^2
\frac{d}{d z}
{\raisebox{0.95ex}{$\raisebox{-0.2ex}[0ex][0ex]{\scriptsize$\epsilon$}
\atop\raisebox{-0.2ex}[0ex][0ex]{$\mathstrut^{A\hspace{-0.1ex}}\hspace{-0.5ex}\Delta_{r}$}$}}
(z)
\!\!&\equiv&
e^{-{\mathrm{i}}{{\epsilon}}\pi}(\lambda+\mu^2)z^2
 \,
{\raisebox{0.95ex}{$\raisebox{-0.2ex}[0ex][0ex]{\scriptsize$\epsilon$}
\atop\raisebox{-0.2ex}[0ex][0ex]{$\mathstrut^{A\hspace{-0.1ex}}\hspace{-0.5ex}\Delta_{p}$}$}}
(z)
  +e^{-{\mathrm{i}}{{\epsilon}}\pi}z^2\,
{\raisebox{0.95ex}{$\raisebox{-0.2ex}[0ex][0ex]{\scriptsize$\epsilon$}
\atop\raisebox{-0.2ex}[0ex][0ex]{$\mathstrut^{A\hspace{-0.1ex}}\hspace{-0.5ex}\Delta_{s}$}$}}
(z)
 \nonumber
\\[-1.5ex]
 &&
 +\big(e^{{\mathrm{i}}{{\epsilon}}\pi}\mu - 2({\ell}-1)z\big)
{\raisebox{0.95ex}{$\raisebox{-0.2ex}[0ex][0ex]{\scriptsize$\epsilon$}
\atop\raisebox{-0.2ex}[0ex][0ex]{$\mathstrut^{A\hspace{-0.1ex}}\hspace{-0.5ex}\Delta_{r}$}$}}
 (z)
    \nonumber
    \\
     &&  
 -e^{-{\mathrm{i}}{{\epsilon}}\pi}z^2
 \,\!\Delta_{{r}}
   ( e^{{\mathrm{i}}{{\epsilon}}\pi}\!/z )
 -e^{{\mathrm{i}}{{\epsilon}}{\ell}\pi} z^{2(2 -{\ell})}\big(\mu \,\,
      \!\Delta_{{p}}( z )
 + \,\!\Delta_{{q}}( z )
 \big)
  \nonumber
  \\
  &  &  \relax \hspace{-9ex}
  \llap{$-$}
  (1+e^{-{\mathrm{i}}{{\epsilon}}\pi})\,
  e^{{\mathrm{i}}{{\epsilon}}{\ell} \pi} z^{2(1-{\ell})}\times
  \nonumber
   \\
 &&
 \relax\;\; 
 \left(
 (\lambda+(2-e^{{\mathrm{i}}{{\epsilon}}\pi})\mu^2)z^2{p}(z) +(1-e^{{\mathrm{i}}{{\epsilon}}\pi})\mu{q}(z)
 \right.
   \nonumber\\
 &&
 \relax\hspace{29.3ex}
 \left.
 +z^2(\mu{}z^2{r}(z)+{s}(z) )
 \vphantom{     (\lambda+(2-e^{{\mathrm{i}}{{\epsilon}}\pi})\mu^2)z^2{p}(z)      }
 \right),
 \end{eqnarray}
%
%
%
%
%
 \begin{eqnarray}
\label{eq:idAs}
z^3
\frac{d}{d z}
{\raisebox{0.95ex}{$\raisebox{-0.2ex}[0ex][0ex]{\scriptsize$\epsilon$}
\atop\raisebox{-0.2ex}[0ex][0ex]{$\mathstrut^{A\hspace{-0.1ex}}\hspace{-0.5ex}\Delta_{s}$}$}}
(z)
\!\!&\equiv&
e^{-{\mathrm{i}}{{\epsilon}}\pi}(\lambda+\mu^2)z^3
{\raisebox{0.95ex}{$\raisebox{-0.2ex}[0ex][0ex]{\scriptsize$\epsilon$}
\atop\raisebox{-0.2ex}[0ex][0ex]{$\mathstrut^{A\hspace{-0.1ex}}\hspace{-0.5ex}\Delta_{q}$}$}}
(z)
 \nonumber
\\[-1.5ex] &&\!\!
-e^{{\mathrm{i}}{{\epsilon}}\pi}
(\lambda{}z - e^{{\mathrm{i}}{{\epsilon}}\pi}({\ell}+1)\mu)
{\raisebox{0.95ex}{$\raisebox{-0.2ex}[0ex][0ex]{\scriptsize$\epsilon$}
\atop\raisebox{-0.2ex}[0ex][0ex]{$\mathstrut^{A\hspace{-0.1ex}}\hspace{-0.5ex}\Delta_{r}$}$}}
(z)
 \nonumber\\
 &&\!\!
-z^2
({\ell} - 1 - e^{-{\mathrm{i}}{{\epsilon}}\pi}\mu{}z)
 \,
{\raisebox{0.95ex}{$\raisebox{-0.2ex}[0ex][0ex]{\scriptsize$\epsilon$}
\atop\raisebox{-0.2ex}[0ex][0ex]{$\mathstrut^{A\hspace{-0.1ex}}\hspace{-0.5ex}\Delta_{s}$}$}}
(z)
   \nonumber
\\
 &&
- 
 e^{-{\mathrm{i}}{{\epsilon}}\pi}z^3
  \,\,\!\Delta_{{s}}
                 ( e^{{\mathrm{i}}{{\epsilon}}\pi}\!/z )
  \nonumber\\
 &&
+ 
e^{{\mathrm{i}}{{\epsilon}}{\ell}\pi}
    z^{3-2{\ell}}
 \big(
\mu^2\,
\!\Delta_{{p}}(z)
 + \mu \,\,
\!\Delta_{{q}}(z)
 +\mu{}z^2\,\!\Delta_{{r}}(z)
 +\!\Delta_{{s}}(z)
 \big)
 \nonumber
\\
&&
   \relax\hspace{-8.8ex}
  \llap{$+$}  (1 + e^{-{\mathrm{i}}{{\epsilon}}\pi})\,
  e^{{\mathrm{i}}{{\epsilon}}{\ell}\pi} z^{2(1-{\ell})}\times
 \nonumber\\
&&
   \relax\hspace{-8ex} \hphantom{ (1 + e^{-{\mathrm{i}}{{\epsilon}}\pi}) }
\left(
\mu{}z\big(\lambda+\mu^2 - (e^{{\mathrm{i}}{{\epsilon}}\pi}\lambda+\mu^2)z^2
 \right.
 \nonumber
 \\
&&
\relax\hspace{-3.2ex}    \hphantom{ (1 + e^{-{\mathrm{i}}{{\epsilon}}\pi}) }
+
\mu
e^{{\mathrm{i}}{{\epsilon}}\pi}
\big(e^{{\mathrm{i}}{{\epsilon}}\pi} -1)({\ell} +1)z\big){p}(z)
 \nonumber
 \\
 &&\relax\hspace{-8ex}   \hphantom{ (1 + e^{-{\mathrm{i}}{{\epsilon}}\pi}) \;}
 +(e^{2{\mathrm{i}}{{\epsilon}}\pi}
 ({\ell}+1)\mu - \mu^2 z
  \nonumber
\end{eqnarray}\begin{eqnarray}
 &&
\relax\hspace{-3.9ex}  \phantom{  (1 + e^{-{\mathrm{i}}{{\epsilon}}\pi}) + ( }
 -e^{{\mathrm{i}}{{\epsilon}}\pi}
 (({\ell}+1)\mu+\lambda{}z)\big){q}(z)
 \nonumber
 \\
 && \left.
   \relax\hspace{-8ex}   \hphantom{ (1 + e^{-{\mathrm{i}}{{\epsilon}}\pi}) \;}
 +z^3
 \big((\lambda+\mu^2(1 - z^2)){r}(z)-\mu{}{s}(z)\big)
 \vphantom{  \big(\lambda+\mu^2 - (e^{{\mathrm{i}}{{\epsilon}}\pi}\lambda+\mu^2)z^2  }\!
\right).
 \end{eqnarray}
The symbols
\,$\Delta{{{\mbox{\tiny\ding{74}}}}} (z)$,
where {{\mbox{\small\ding{74}}}}$\;\in \{{p},{q},{r},{s}\},$
stand for the differences of the left- and right-hand sides of the equations
\eqref{eq:ppp'},
\eqref{eq:qqq'},
\eqref{eq:rrr'},
\eqref{eq:sss'}, respectively.
The symbol
$ {{\epsilon}} $ denotes the real parameter, $ {{\epsilon}} \in [-1,1]$.
The definitions of the functions
$
{\raisebox{0.95ex}{$\raisebox{-0.2ex}[0ex][0ex]{\scriptsize$\epsilon$}
\atop\raisebox{-0.2ex}[0ex][0ex]{$\mathstrut^{A\hspace{-0.1ex}}\hspace{-0.5ex}
\Delta_{{{\mbox{\tiny\ding{74}}}}}$}$}}
(z)
$
read
\begin{equation}
\label{eq:490}
\begin{aligned}
{\raisebox{0.95ex}{$\raisebox{-0.2ex}[0ex][0ex]{\scriptsize$\epsilon$}
\atop\raisebox{-0.2ex}[0ex][0ex]{$\mathstrut^{A\hspace{-0.1ex}}\hspace{-0.5ex}\Delta_{p}$}$}}
(z) =&\;
{p}(e^{{\mathrm{i}}{{\epsilon}}\pi}/z)
+
e^{{\mathrm{i}}{{\epsilon}}{\ell}\pi}
z^{2 \left(1 - {\ell}\right)}
{p}(z)
,
\\
{\raisebox{0.95ex}{$\raisebox{-0.2ex}[0ex][0ex]{\scriptsize$\epsilon$}
\atop\raisebox{-0.2ex}[0ex][0ex]{$\mathstrut^{A\hspace{-0.1ex}}\hspace{-0.5ex}\Delta_{q}$}$}}
(z) =&\;
{q}(e^{{\mathrm{i}}{{\epsilon}}\pi}/z)
-
e^{{\mathrm{i}}{{\epsilon}}{\ell}\pi}
z^{-2{\ell}}
\left(
\mu\, {p}(z)
+z^2 {r}(z)\right)
,
\\
{\raisebox{0.95ex}{$\raisebox{-0.2ex}[0ex][0ex]{\scriptsize$\epsilon$}
\atop\raisebox{-0.2ex}[0ex][0ex]{$\mathstrut^{A\hspace{-0.1ex}}\hspace{-0.5ex}\Delta_{r}$}$}}
(z) =&\;
{r}(e^{{\mathrm{i}}{{\epsilon}}\pi}/z)
-
e^{{\mathrm{i}}{{\epsilon}}{\ell}\pi}
z^{2 (1 - {\ell})} \left(\mu z^2 {p}(z) + {q}(z)\right)
,
\\
{\raisebox{0.95ex}{$\raisebox{-0.2ex}[0ex][0ex]{\scriptsize$\epsilon$}
\atop\raisebox{-0.2ex}[0ex][0ex]{$\mathstrut^{A\hspace{-0.1ex}}\hspace{-0.5ex}\Delta_{s}$}$}}
(z) =&\;
{s}(e^{{\mathrm{i}}{{\epsilon}}\pi}/z)
+
e^{{\mathrm{i}}{{\epsilon}}{\ell}\pi}
z^{-2{\ell}} \left(
\mu\left(\mu z^2 {p}(z) + {q}(z)\right)
+
z^2 \left(\mu z^2 {r}(z) + {s}(z)\right)
\right).
\end{aligned}
\end{equation}
Hence, as a matter of fact,
the equalities
\eqref{eq:idAp}-\eqref{eq:idAs}
signify the four
pairwise
coincidences, upon simplification,
of certain expressions 
constructed in two different ways
from
{\em arbitrary\/}
holomorphic functions
 $ {p},{q},{r},{s} $
and their first order derivatives.

Let us notice now that 
$\lim_{{{\epsilon}}\nearrow 1  }\!\!
{\raisebox{0.95ex}{$\raisebox{-0.2ex}[0ex][0ex]{\scriptsize$\epsilon$}
\atop\raisebox{-0.2ex}[0ex][0ex]{$\mathstrut^{A\hspace{-0.1ex}}\hspace{-0.5ex}
\Delta_{{{\mbox{\tiny\ding{74}}}}}$}$}}
(z)
$
coincide with the functions
$
{\mbox{$\mathstrut^{A}\!\Delta_{{{\mbox{\tiny\ding{74}}}}}$}}(z)
$
introduced 
 in the beginning of the proof of Theorem \ref{t:010} %
and involved in {{Eq.s}}{~}\eqref{eq:210}.
It is worth reminding
that they 
were defined as the differences of the
left- and right-hand sides of the equations
\eqref{eq:pppA}-\eqref{eq:sssA}.
They are correctly defined if
the functions  $ {p},{q},{r},{s} $
are holomorphic in the vicinity of $z={\mathrm{i}}$.
Besides, for 
 $ {{\epsilon}}=1 $, the last summands in the right-hand sides of
the equalities
\eqref{eq:idAp}-\eqref{eq:idAs}, which are
proportional to either
$ (1+e^{{\mathrm{i}}{{\epsilon}}\pi})  $ or $ (1+e^{-{\mathrm{i}}{{\epsilon}}\pi})  $,
vanish. 
Finally, it remains to note that if the functions
$ {p},{q},{r},{s} $ obey
 {{Eq.s}}{~}\eqref{eq:ppp'}-%
\eqref{eq:sss'}
then
the differences
$\,\Delta{{{\mbox{\tiny\ding{74}}}}}(z)$
become identically zero
for all {{\mbox{\small\ding{74}}}}$\;\in\{{p},{q},{r},{s}\}$
 and
all the summands which contain them can also be dropped out.
After such simplifications, 
comparing the resulting form of
{{Eq.s}}{~}\eqref{eq:idAp}-\eqref{eq:idAp} with {{Eq.s}}{~}\eqref{eq:210},
one easily finds that they coincide.
Thus the equalities \eqref{eq:210} hold true.

\section{Identities leading to {{Eq.s}}{~}\eqref{eq:220}}\label{a:020}

{{Eq.s}}{~}\eqref{eq:220} follow  from the
identities given below which can be, in principle, verified
by straightforward computations. %
Namely,
for any functions
$ {p},{q},{r},{s} $
holomorphic, at least, in the vicinity of $z=1$
the following identities 
hold true
%
%
%
%
 \begin{eqnarray}
\label{eq:idCp}
z^2
\frac{d}{d z}
{\mbox{$\mathstrut^{C\!}\!\Delta_{p}$}}(z)
&\equiv&
-z({\ell} - 1 + \mu{}z) 
\,{\mbox{$\mathstrut^{C\!}\!\Delta_{p}$}}(z)
+z^2\,{\mbox{$\mathstrut^{C\!}\!\Delta_{q}$}}(z)
-{\mbox{$\mathstrut^{C\!}\!\Delta_{r}$}}(z)
\nonumber\\&&
+z^2\,\!\Delta_{{p}}(1/z)
-(\lambda+\mu^2)^{-1}
z^{2(1-{\ell})}
(\mu{}z^2\,\!\Delta_{{r}}(z)
+\,\!\Delta_{{s}}(z)),
 \\
%
%
%
\label{eq:idCq}
z^3
\frac{d}{d z}
{\mbox{$\mathstrut^{C\!}\!\Delta_{q}$}}(z)
&\equiv&
 (({\ell}+1)\mu - \lambda{}z)\,{\mbox{$\mathstrut^{C\!}\!\Delta_{p}$}}(z)
-\mu{} z \,{\mbox{$\mathstrut^{C\!}\!\Delta_{q}$}}(z)
-z\,{\mbox{$\mathstrut^{C\!}\!\Delta_{s}$}}(z)
\nonumber\\&&
+z \,\,\!\Delta_{{q}}(1/z)
-(\lambda+\mu^2)^{-1}
z^{1-2{\ell}}
(\lambda{}z^2\,\!\Delta_{{r}}(z)
-\mu\,\,\!\Delta_{{s}}(z)),
 \end{eqnarray}
%
%
%
%
 \begin{eqnarray}
\label{eq:idCr}
z^2
\frac{d}{d z}
{\mbox{$\mathstrut^{C\!}\!\Delta_{r}$}}(z)
&\equiv&
(\lambda+\mu^2) z^2\,{\mbox{$\mathstrut^{C\!}\!\Delta_{p}$}}(z)
+(\mu + 2(1 - {\ell})z)\,{\mbox{$\mathstrut^{C\!}\!\Delta_{r}$}}(z)
+z^2\,{\mbox{$\mathstrut^{C\!}\!\Delta_{s}$}}(z)
\nonumber\\&&
+z^2 \,\!\Delta_{{r}}(1/z)
-z^{2(2 - {\ell})}
(\mu\,\,\!\Delta_{{p}}(z)
+\,\,\!\Delta_{{q}}(z)),
\\ 
%
%
%
%
\label{eq:idCs}
z^3
\frac{d}{d z}
{\mbox{$\mathstrut^{C\!}\!\Delta_{s}$}}(z)
&\equiv&
   (\lambda+\mu^2) z^3 \,{\mbox{$\mathstrut^{C\!}\!\Delta_{q}$}}(z)
  +(({\ell} + 1)\mu - \lambda{} z)\,{\mbox{$\mathstrut^{C\!}\!\Delta_{r}$}}(z)
\nonumber\\&&
  +z^2 (1 - {\ell} + \mu{}z)\,{\mbox{$\mathstrut^{C\!}\!\Delta_{s}$}}(z)
\nonumber\\&&
+z^3 \,\!\Delta_{{s}}(1/z)
-z^{3-2{\ell}}
(\lambda\,\,\!\Delta_{{p}}(z)
 -\mu\,\,\!\Delta_{{q}}(z)).
\end{eqnarray}
Here the symbols
${\mbox{$\mathstrut^{C\!}\!\Delta_{{{\mbox{\tiny\ding{74}}}}}$}}(z)$,
where {{\mbox{\small\ding{74}}}}$\;\in\{{p},{q},{r},{s}\}$,
used already
in 
the proof of Theorem \ref{t:020},
stand for
the differences
of the left- and right-hand sides of
{{Eq.s}}{~}\eqref{eq:pppC}-\eqref{eq:sssC}.
The symbols
\,$\Delta{{{\mbox{\tiny\ding{74}}}}}(z),$
where {{\mbox{\small\ding{74}}}}$\;\in\{{p},{q},{r},{s}\}$,
denote
 the differences of the left- and right-hand sides of the equations
\eqref{eq:ppp'}-\eqref{eq:sss'}.
Thus the equalities
\eqref{eq:idCp}-\eqref{eq:idCs}
signify the
pairwise
coincidences, upon simplification,
of certain expressions 
constructed in two different ways
from
arbitrary holomorphic functions
 $ {p},{q},{r},{s} $
and their first order derivatives.
Obviously, these expressions 
are correctly defined if the
above four functions
${p},{q},{r},{s}$
are holomorphic in the vicinity of $z=1$. 

Finally,
if
the functions $ {p},{q},{r},{s} $
are not arbitrary but
verify
 {{Eq.s}}{~}\eqref{eq:ppp'}-%
\eqref{eq:sss'}
then
the differences
$\,\Delta{{{\mbox{\tiny\ding{74}}}}}(z)$
vanish and the identities
\eqref{eq:idCp}-\eqref{eq:idCs}
convert to
{{Eq.s}}{~}\eqref{eq:220}
which are therefore the
direct consequence of  {{Eq.s}}{~}\eqref{eq:ppp'}-%
\eqref{eq:sss'}.

\section{Identities leading to {{Eq.s}}{~}\eqref{eq:230}}\label{a:030}

{{Eq.s}}{~}\eqref{eq:230},
utilized in the  proof of Theorem \ref{t:030},
can be obtained from the four
identities displayed below
which are verifiable by straightforward
computations. %
Namely, it can be shown that
 {
%
%
%
%
\begin{eqnarray}
\label{eq:idBp}
z^2
\frac{d}{d z}
{\raisebox{0.95ex}{$\raisebox{-0.2ex}[0ex][0ex]{\scriptsize$\epsilon$}
\atop\raisebox{-0.2ex}[0ex][0ex]{$\mathstrut^{B\hspace{-0.2ex}}\hspace{-0.5ex}\Delta_{p}$}$}}
(z)
&\equiv&
(-\mu+({\ell} - 1)z)
\,
{\raisebox{0.95ex}{$\raisebox{-0.2ex}[0ex][0ex]{\scriptsize$\epsilon$}
\atop\raisebox{-0.2ex}[0ex][0ex]{$\mathstrut^{B\hspace{-0.2ex}}\hspace{-0.5ex}\Delta_{p}$}$}}
(z)
-e^{-{\mathrm{i}}{{\epsilon}}\pi}
{\raisebox{0.95ex}{$\raisebox{-0.2ex}[0ex][0ex]{\scriptsize$\epsilon$}
\atop\raisebox{-0.2ex}[0ex][0ex]{$\mathstrut^{B\hspace{-0.2ex}}\hspace{-0.5ex}\Delta_{q}$}$}}
(z)
-z^2
{\raisebox{0.95ex}{$\raisebox{-0.2ex}[0ex][0ex]{\scriptsize$\epsilon$}
\atop\raisebox{-0.2ex}[0ex][0ex]{$\mathstrut^{B\hspace{-0.2ex}}\hspace{-0.5ex}\Delta_{r}$}$}}
(z)
 \nonumber\\&&
+e^{-{\mathrm{i}}{{\epsilon}}\pi}
\,\!\Delta_{{p}}
          ({e^{{\mathrm{i}}{{\epsilon}}\pi}}{}z)
-e^{{\mathrm{i}}{\ell}{{\epsilon}}\pi}  (\lambda+\mu^2)^{-1}
(\mu{}z^2\,\!\Delta_{{r}}(z)
+\,\,\!\Delta_{{s}}(z))
\nonumber %
 \\
&&  \hspace{-3ex}
+(1 + e^{-{\mathrm{i}}{{\epsilon}}\pi})
\big(\mu\,{p}({e^{{\mathrm{i}}{{\epsilon}}\pi}} z)+e^{{\mathrm{i}}{{\epsilon}}\pi}z^2{r}({e^{{\mathrm{i}}{{\epsilon}}\pi}}{}z)
\nonumber\\&&  \hspace{-5ex} \hphantom{  -(1 + e^{-{\mathrm{i}}{{\epsilon}}\pi})  }
+e^{{\mathrm{i}}{\ell}{{\epsilon}}\pi}({q}(z)+\mu{}z^2{p}(z)
\nonumber\\&&  \hspace{-5ex} \hphantom{  -(1 + e^{-{\mathrm{i}}{{\epsilon}}\pi})( +e^{{\mathrm{i}}{\ell}{{\epsilon}}\pi}(   }
+\mu(\lambda+\mu^2)^{-1}z^2({s}(z)+\mu{}z^2{r}(z))) \big),
\\ 
%
%
%
\label{eq:idBq}
\phantom{z^2}
\frac{d}{d z}
{\raisebox{0.95ex}{$\raisebox{-0.2ex}[0ex][0ex]{\scriptsize$\epsilon$}
\atop\raisebox{-0.2ex}[0ex][0ex]{$\mathstrut^{B\hspace{-0.2ex}}\hspace{-0.5ex}\Delta_{q}$}$}}
(z)
&\equiv&
-(\lambda + ({\ell}+1)\mu{} z) \,
{\raisebox{0.95ex}{$\raisebox{-0.2ex}[0ex][0ex]{\scriptsize$\epsilon$}
\atop\raisebox{-0.2ex}[0ex][0ex]{$\mathstrut^{B\hspace{-0.2ex}}\hspace{-0.5ex}\Delta_{p}$}$}}
(z)
-\mu                           \,
{\raisebox{0.95ex}{$\raisebox{-0.2ex}[0ex][0ex]{\scriptsize$\epsilon$}
\atop\raisebox{-0.2ex}[0ex][0ex]{$\mathstrut^{B\hspace{-0.2ex}}\hspace{-0.5ex}\Delta_{q}$}$}}
(z)
-                              \,
{\raisebox{0.95ex}{$\raisebox{-0.2ex}[0ex][0ex]{\scriptsize$\epsilon$}
\atop\raisebox{-0.2ex}[0ex][0ex]{$\mathstrut^{B\hspace{-0.2ex}}\hspace{-0.5ex}\Delta_{s}$}$}}
(z)
 \nonumber\\&&
+e^{{\mathrm{i}}{{\epsilon}}\pi}
                             \,\!\Delta_{{q}}
             ({e^{{\mathrm{i}}{{\epsilon}}\pi}}{}z)
+e^{{\mathrm{i}}{\ell}{{\epsilon}}\pi}
\big(
\mu \,\,\!\Delta_{{p}}(z)
 + \,\,\!\Delta_{{q}}(z)
\big)
 \nonumber\\&&
+e^{{\mathrm{i}}{\ell}{{\epsilon}}\pi} \mu\,  (\lambda+\mu^2)^{-1}
\big(
\mu{}z^2 \,\!\Delta_{{r}}(z)
 + \,\,\!\Delta_{{s}}(z)
\big)
 \nonumber
 \\
 &&  \hspace{-5ex}
+ (1 + e^{{\mathrm{i}}{{\epsilon}}\pi})
\big((\lambda +(1 - e^{{\mathrm{i}}{{\epsilon}}\pi})({\ell} + 1)\mu{}z){{p}}({e^{{\mathrm{i}}{{\epsilon}}\pi}}{}z)
 \nonumber\\&&   \hspace{-5ex} \hphantom{ -(1 + e^{{\mathrm{i}}{{\epsilon}}\pi}) \big( }
+\mu\,{{q}}({e^{{\mathrm{i}}{{\epsilon}}\pi}}{}z) + {{s}}({e^{{\mathrm{i}}{{\epsilon}}\pi}}{}z)
\big),
 \\  
%
%
%
\label{eq:idBr}
z^2
\frac{d}{d z}
{\raisebox{0.95ex}{$\raisebox{-0.2ex}[0ex][0ex]{\scriptsize$\epsilon$}
\atop\raisebox{-0.2ex}[0ex][0ex]{$\mathstrut^{B\hspace{-0.2ex}}\hspace{-0.5ex}\Delta_{r}$}$}}
(z)
&\equiv&
(\lambda+\mu^2)                \,
{\raisebox{0.95ex}{$\raisebox{-0.2ex}[0ex][0ex]{\scriptsize$\epsilon$}
\atop\raisebox{-0.2ex}[0ex][0ex]{$\mathstrut^{B\hspace{-0.2ex}}\hspace{-0.5ex}\Delta_{p}$}$}}
(z)
+(2({\ell} - 1) + \mu{}z)z     \,
{\raisebox{0.95ex}{$\raisebox{-0.2ex}[0ex][0ex]{\scriptsize$\epsilon$}
\atop\raisebox{-0.2ex}[0ex][0ex]{$\mathstrut^{B\hspace{-0.2ex}}\hspace{-0.5ex}\Delta_{r}$}$}}
(z)
+                                
{\raisebox{0.95ex}{$\raisebox{-0.2ex}[0ex][0ex]{\scriptsize$\epsilon$}
\atop\raisebox{-0.2ex}[0ex][0ex]{$\mathstrut^{B\hspace{-0.2ex}}\hspace{-0.5ex}\Delta_{s}$}$}}
(z)
 \nonumber\\&&
+e^{-{\mathrm{i}}{{\epsilon}}\pi}
             \,\!\Delta_{{r}}({e^{{\mathrm{i}}{{\epsilon}}\pi}}{}z)
+e^{{\mathrm{i}}{\ell}{{\epsilon}}\pi}
       \,\!\Delta_{{r}}(z)
 \nonumber
 \end{eqnarray}\begin{eqnarray}
 &&    \hspace{-5ex}
-(1 + e^{-{\mathrm{i}}{{\epsilon}}\pi})
\big(
(\lambda+\mu^2){p}({e^{{\mathrm{i}}{{\epsilon}}\pi}}{}z)
 \nonumber\\&&  \hspace{-5ex} \hphantom{  +(1 + e^{-{\mathrm{i}}{{\epsilon}}\pi})  \big( }
+e^{{\mathrm{i}}{{\epsilon}}\pi} \mu{}z^2 {r}({e^{{\mathrm{i}}{{\epsilon}}\pi}}{}z)
+ {s}({e^{{\mathrm{i}}{{\epsilon}}\pi}}{}z)
\big),
\\ 
%
%
%
\label{eq:idBs}
z^2
\frac{d}{d z}
{\raisebox{0.95ex}{$\raisebox{-0.2ex}[0ex][0ex]{\scriptsize$\epsilon$}
\atop\raisebox{-0.2ex}[0ex][0ex]{$\mathstrut^{B\hspace{-0.2ex}}\hspace{-0.5ex}\Delta_{s}$}$}}
(z)
&\equiv&
(\lambda+\mu^2)                         \,
{\raisebox{0.95ex}{$\raisebox{-0.2ex}[0ex][0ex]{\scriptsize$\epsilon$}
\atop\raisebox{-0.2ex}[0ex][0ex]{$\mathstrut^{B\hspace{-0.2ex}}\hspace{-0.5ex}\Delta_{q}$}$}}
(z)
-(\lambda +  ({\ell} + 1)\mu{}z)z^2       
{\raisebox{0.95ex}{$\raisebox{-0.2ex}[0ex][0ex]{\scriptsize$\epsilon$}
\atop\raisebox{-0.2ex}[0ex][0ex]{$\mathstrut^{B\hspace{-0.2ex}}\hspace{-0.5ex}\Delta_{r}$}$}}
(z)
\nonumber\\[-1ex]&&
+(\mu + ({\ell} - 1)z)                   \,
{\raisebox{0.95ex}{$\raisebox{-0.2ex}[0ex][0ex]{\scriptsize$\epsilon$}
\atop\raisebox{-0.2ex}[0ex][0ex]{$\mathstrut^{B\hspace{-0.2ex}}\hspace{-0.5ex}\Delta_{s}$}$}}
(z)
 \nonumber\\&&
+ e^{-{\mathrm{i}}{{\epsilon}}\pi}
                  \,\!\Delta_{{s}}({e^{{\mathrm{i}}{{\epsilon}}\pi}}{}z)
- e^{{\mathrm{i}}{\ell}{{\epsilon}}\pi}
((\lambda+\mu^2) \,\,\!\Delta_{{p}}(z)
 + \mu{}z^2 \,\!\Delta_{{r}}(z))
 \nonumber\\&&    \hspace{-5ex}
-(1 + e^{-{\mathrm{i}}{{\epsilon}}\pi})
\big(
(\lambda+\mu^2){q}({e^{{\mathrm{i}}{{\epsilon}}\pi}}{}z) + \mu\,{s}({e^{{\mathrm{i}}{{\epsilon}}\pi}}{}z)
 \nonumber\\&& \hspace{-5ex} \hphantom{  -(1 + e^{-{\mathrm{i}}{{\epsilon}}\pi})  }\;
-e^{{\mathrm{i}}{{\epsilon}}\pi} z^2
(\lambda +(1 -  e^{{\mathrm{i}}{{\epsilon}}\pi} )({\ell} + 1)\mu{}z )
{r}({e^{{\mathrm{i}}{{\epsilon}}\pi}}{}z)
\big).
\end{eqnarray}
} 
Here $ {{\epsilon}}\in [0,1] $ is the auxiliary real parameter,
the symbols $ {p},{q},{r},{s} $ stay for
arbitrary
functions holomorphic in the vicinity of
an arc of the circle connecting $ e^{-{\mathrm{i}}{{\epsilon}}\pi/2} $ with $ e^{{\mathrm{i}}{{\epsilon}}\pi/2} $
and passing
inbetween them
through $+1$  counter-clockwise.
The functions
$
{\raisebox{0.95ex}{$\raisebox{-0.2ex}[0ex][0ex]{\scriptsize$\epsilon$}
\atop\raisebox{-0.2ex}[0ex][0ex]{$\mathstrut^{B\hspace{-0.2ex}}\hspace{-0.5ex}
\Delta_{{{\mbox{\tiny\ding{74}}}}}$}$}}
\!(z)$,
where {{\mbox{\small\ding{74}}}}$\;\in\{{p},{q},{r},{s}\},
$
are
defined as follows.
 {
%
\begin{equation}
\label{eq:580}
\begin{aligned}
{\raisebox{0.95ex}{$\raisebox{-0.2ex}[0ex][0ex]{\scriptsize$\epsilon$}
\atop\raisebox{-0.2ex}[0ex][0ex]{$\mathstrut^{B\hspace{-0.2ex}}\hspace{-0.5ex}\Delta_{p}$}$}}
(z)
=&\; 
{p}({e^{{\mathrm{i}}{{\epsilon}}\pi}}{}z)
 -
e^{{\mathrm{i}}{{\epsilon}}{\ell}\pi} (\lambda+\mu^2)^{-1}
\big(\mu{}z^2{r}(z)+{s}(z)\big)
,
%
\\
{\raisebox{0.95ex}{$\raisebox{-0.2ex}[0ex][0ex]{\scriptsize$\epsilon$}
\atop\raisebox{-0.2ex}[0ex][0ex]{$\mathstrut^{B\hspace{-0.2ex}}\hspace{-0.5ex}\Delta_{q}$}$}}
(z)
=&\; 
{q}({e^{{\mathrm{i}}{{\epsilon}}\pi}}{}z)
+
e^{{\mathrm{i}}{{\epsilon}}{\ell}\pi}
\left(
     (\mu\, z^2 {p}(z) +{q}(z))
\right.
\\
&\; 
   \hphantom{ {q}({e^{{\mathrm{i}}{{\epsilon}}\pi}}{}z) +  e^{{\mathrm{i}}{{\epsilon}}{\ell}\pi} }
\left.
+\mu(\lambda+\mu^2)^{-1}z^2
     (\mu{}z^2{r}(z)+{s}(z))
\right)
,
%
\\
{\raisebox{0.95ex}{$\raisebox{-0.2ex}[0ex][0ex]{\scriptsize$\epsilon$}
\atop\raisebox{-0.2ex}[0ex][0ex]{$\mathstrut^{B\hspace{-0.2ex}}\hspace{-0.5ex}\Delta_{r}$}$}}
(z)
=&\; 
{r}({e^{{\mathrm{i}}{{\epsilon}}\pi}}{}z)
+
e^{{\mathrm{i}}{{\epsilon}}{\ell}\pi}
{r}(z),
%
\\
{\raisebox{0.95ex}{$\raisebox{-0.2ex}[0ex][0ex]{\scriptsize$\epsilon$}
\atop\raisebox{-0.2ex}[0ex][0ex]{$\mathstrut^{B\hspace{-0.2ex}}\hspace{-0.5ex}\Delta_{s}$}$}}
(z)
=&\; 
{s}({e^{{\mathrm{i}}{{\epsilon}}\pi}}{}z)
 -
e^{{\mathrm{i}}{{\epsilon}}{\ell}\pi}
\left(
(\lambda+\mu^2) {p}(z)+\mu{}z^2{r}(z)\right).
\end{aligned}
\end{equation}%
}\\[-2ex] 
Lastly,
the symbols
\,$
\Delta{{{\mbox{\tiny\ding{74}}}}}(z)
,$ where
{{\mbox{\small\ding{74}}}}$\;\in\{{p},{q},{r},{s}\}$,
stand for the differences of the left- and right-hand sides of the equations
\eqref{eq:ppp'},
\eqref{eq:qqq'},
\eqref{eq:rrr'},
\eqref{eq:sss'}, respectively.
They were also used in 
{{Eq.s}}{~}\eqref{eq:idAp}-\eqref{eq:idAs}, \eqref{eq:idCp}-\eqref{eq:idCs}.

Thus the equalities
\eqref{eq:idBp}-\eqref{eq:idBs}
signify the
pairwise
coincidences, upon simplification,
of certain expressions 
constructed in two different ways
from
arbitrary holomorphic functions
 $ {p},{q},{r},{s} $
and their first order derivatives.

Let us now consider the case $ {{\epsilon}}=1 $.
The definitions \eqref{eq:580} are pertinent
if the domain of functions $ {p},{q},{r},{s} $
covers 
$\pm{\mathrm{i}}$ and $ 1  $, i.e.,  in particular, if they are holomorphic
on the circular arc passing 
through $-{\mathrm{i}}, 1$ and $+{\mathrm{i}}$.
The solutions
of the Cauchy problem
for
{{Eq.s}}{~}\eqref{eq:ppp'}-
\eqref{eq:sss'}
with initial data specified at $z=1$
possess such a property. 
One has for them, by definition,
$\;
\Delta{{{\mbox{\tiny\ding{74}}}}}(z)
=0.
$
Besides, due to the above choice of $ {{\epsilon}} $,
the summands
in right-hand sides 
of 
{{Eq.s}}~\eqref{eq:idBp}-\eqref{eq:idBs}
 involving 
 either
 the factor
$ (1+e^{{\mathrm{i}}{{\epsilon}}\pi})  $ or the factor $ (1+e^{-{\mathrm{i}}{{\epsilon}}\pi})  $ 
have to be discarded 
as well.

Taking  the above simplifications into account, {{Eq.s}}{~}\eqref{eq:230}
follow %
since the argument 
for which
the transformed
{\mbox{$pqr\!s$}}-functions on the left in {{Eq.s}}{~}\eqref{eq:pppB}-\eqref{eq:sssB}
have to be evaluated
is exactly the limit of $ e^{{\mathrm{i}}{{\epsilon}}\pi}z $
reached
as  ${{\epsilon}}\nearrow 1 $ (provided $z$ belongs to 
the vicinity of $-{\mathrm{i}}$, at least).

\section{Identities leading to {{Eq.s}}{~}\eqref{eq:240}}\label{a:040}

Let us consider ``the deformed differences''
$
{\raisebox{0.95ex}{$\raisebox{-0.2ex}[0ex][0ex]{\scriptsize$\epsilon$}
\atop\raisebox{-0.2ex}[0ex][0ex]{$\mathstrut^{B\hspace{-0.2ex}}\hspace{-0.5ex}
\Delta_{{{\mbox{\tiny\ding{74}}}}}$}$}}
(z)$,
 where {{\mbox{\small\ding{74}}}}$\;\in \{{p},{q},{r},{s} \},$
 and ${{\epsilon}}\in[-1,1]$ is the real parameter,
defined by the formulas \eqref{eq:580}.
If ${{{\epsilon}}\nearrow 1}$
then they %
reduce to 
the differences
$
{\mbox{$\mathstrut^{B\!}\!\Delta_{{{\mbox{\tiny\ding{74}}}}}$}}(z)
$
of the left-
and right-hand sides
of
{{Eq.s}}{~}\eqref{eq:pppB}-\eqref{eq:sssB}.
We consider their values in the case 
$z=-{\mathrm{i}}=e^{-{\mathrm{i}}\pi/2}$.
The following four double equalities hold true
for arbitrary functions ${p},{q},{r},{s}$ holomorphic on the circular arc
passing through 
$-{\mathrm{i}},+1,$ and $+{\mathrm{i}}$.
%
%
\begin{equation} 
\label{eq:difBp}
\begin{aligned}
{\raisebox{0.95ex}{$\raisebox{-0.2ex}[0ex][0ex]{\scriptsize$\epsilon$}
\atop\raisebox{-0.2ex}[0ex][0ex]{$\mathstrut^{B\hspace{-0.2ex}}\hspace{-0.5ex}\Delta_{p}$}$}}
(e^{-{\mathrm{i}}{{\epsilon}}\pi/2})
= &\;
{p}( e^{{\mathrm{i}}{{\epsilon}}\pi/2})
-
e^{{\mathrm{i}}{}{{\epsilon}}{\ell}\pi}(\lambda+\mu^2)^{-1}
\bigl(
{{\lfloor\![\mbox{-}1]\!\rfloor}}
\mu\,
{r}( e^{-{\mathrm{i}}{{\epsilon}}\pi/2})
+{s}( e^{-{\mathrm{i}}{{\epsilon}}\pi/2})
\bigr)
  \\ \equiv
&\;
{{\lfloor0\rfloor}}
\,{p}( e^{{\mathrm{i}}{{\epsilon}}\pi/2})
\\&
-
e^{{\mathrm{i}}{}{{\epsilon}}{\ell}\pi} (\lambda+\mu^2)^{-1}
\big(
{{\lfloor\![\mbox{-}1]\!\rfloor}}
             \mu \,{\mbox{$\mathstrut^{C\!}\!\Delta_{r}$}}(e^{{\mathrm{i}}{{\epsilon}}\pi/2})
+
\,{\mbox{$\mathstrut^{C\!}\!\Delta_{s}$}}(e^{{\mathrm{i}}{{\epsilon}}\pi/2})
 \big),
\end{aligned}
\end{equation}
%
%
%
%
%
%
%
\begin{equation} 
\label{eq:difBq}
\begin{aligned}
{\raisebox{0.95ex}{$\raisebox{-0.2ex}[0ex][0ex]{\scriptsize$\epsilon$}
\atop\raisebox{-0.2ex}[0ex][0ex]{$\mathstrut^{B\hspace{-0.2ex}}\hspace{-0.5ex}\Delta_{q}$}$}}
(e^{-{\mathrm{i}}{{\epsilon}}\pi/2})
= &\;
{q}( e^{{\mathrm{i}}{{\epsilon}}\pi/2})
+
e^{{\mathrm{i}}{}{{\epsilon}}{\ell}\pi}
\bigl(
{{\lfloor\![\mbox{-}1]\!\rfloor}}
              \mu\,{p}( e^{-{\mathrm{i}}{{\epsilon}}\pi/2})
+{q}( e^{-{\mathrm{i}}{{\epsilon}}\pi/2})
\\& \hspace{10.ex} 
+{{\lfloor\![\mbox{-}1]\!\rfloor}}
              (\lambda+\mu^2)^{-1}\mu\,
({{\lfloor\![\mbox{-}1]\!\rfloor}}
               \,\mu\,
{r}( e^{-{\mathrm{i}}{{\epsilon}}\pi/2})
+{s}( e^{-{\mathrm{i}}{{\epsilon}}\pi/2})
)
\bigr)
           \\\equiv&\,
{{\lfloor\raisebox{0.1ex}{-}1\rfloor}}
\big(
\mu\,{p}( e^{{\mathrm{i}}{{\epsilon}}\pi/2})
 -{q}( e^{{\mathrm{i}}{{\epsilon}}\pi/2})
 -{r}(e^{{\mathrm{i}}{{\epsilon}}\pi/2})
\big)
 \\
&
+{{\lfloor0\rfloor}}
\big(
{q}( e^{{\mathrm{i}}{{\epsilon}}\pi/2})
-\mu\,{p}( e^{{\mathrm{i}}{{\epsilon}}\pi/2})
\big)
\end{aligned}\end{equation}\begin{equation*}\begin{aligned}
&
+ e^{{\mathrm{i}}{{\epsilon}}{\ell}\pi}
\big(
{\mbox{$\mathstrut^{C\!}\!\Delta_{q}$}}(e^{{\mathrm{i}}{{\epsilon}}\pi/2})
+{{\lfloor\![\mbox{-}1]\!\rfloor}}
              \, \mu\,{\mbox{$\mathstrut^{C\!}\!\Delta_{p}$}}(e^{{\mathrm{i}}{{\epsilon}}\pi/2})
        \\& \hphantom{  e^{{\mathrm{i}}{{\epsilon}}{\ell}\pi/2} \big(  }
+{{\lfloor\![\mbox{-}1]\!\rfloor}}
               (\lambda+\mu^2)^{-1} \mu\,
(
{{\lfloor\![\mbox{-}1]\!\rfloor}}
               \,\mu
\,{\mbox{$\mathstrut^{C\!}\!\Delta_{r}$}}(e^{{\mathrm{i}}{{\epsilon}}\pi/2})
+
{\mbox{$\mathstrut^{C\!}\!\Delta_{s}$}}(e^{{\mathrm{i}}{{\epsilon}}\pi/2})
)
\big),
\end{aligned}
\end{equation*}
%
%
%
%
\begin{equation} 
\label{eq:difBr}
\begin{aligned}
{\raisebox{0.95ex}{$\raisebox{-0.2ex}[0ex][0ex]{\scriptsize$\epsilon$}
\atop\raisebox{-0.2ex}[0ex][0ex]{$\mathstrut^{B\hspace{-0.2ex}}\hspace{-0.5ex}\Delta_{r}$}$}}
(e^{-{\mathrm{i}}{{\epsilon}}\pi/2})
= &\;
{r}( e^{{\mathrm{i}}{{\epsilon}}\pi/2})
+e^{{\mathrm{i}}{{\epsilon}}{\ell}\pi}
{r}( e^{-{\mathrm{i}}{{\epsilon}}\pi/2})
             \\\equiv&\;
-\big(
\mu\,{p}( e^{{\mathrm{i}}{{\epsilon}}\pi/2})
-{q}( e^{{\mathrm{i}}{{\epsilon}}\pi/2})
-{r}(e^{{\mathrm{i}}{{\epsilon}}\pi/2})
\big)
\\&\;
+{{\lfloor0\rfloor}}
\big(
{{\lfloor2\rfloor}}
 \mu\,{p}( e^{{\mathrm{i}}{{\epsilon}}\pi/2})
-{q}( e^{{\mathrm{i}}{{\epsilon}}\pi/2})
\big)
+ e^{{\mathrm{i}}{{\epsilon}}{\ell}\pi}
\,{\mbox{$\mathstrut^{C\!}\!\Delta_{r}$}}(e^{{\mathrm{i}}{{\epsilon}}\pi/2}),
\end{aligned}
\end{equation}
%
%
%
%
\begin{equation} 
\label{eq:difBs}
\begin{aligned}
{\raisebox{0.95ex}{$\raisebox{-0.2ex}[0ex][0ex]{\scriptsize$\epsilon$}
\atop\raisebox{-0.2ex}[0ex][0ex]{$\mathstrut^{B\hspace{-0.2ex}}\hspace{-0.5ex}\Delta_{s}$}$}}
(e^{-{\mathrm{i}}{{\epsilon}}\pi/2})
= &\;
{s}( e^{{\mathrm{i}}{{\epsilon}}\pi/2})
-
 e^{{\mathrm{i}}{{\epsilon}}{\ell}\pi}
\big(
{{\lfloor\![\mbox{-}1]\!\rfloor}}
              \mu\,{r}(e^{-{\mathrm{i}}{{\epsilon}}\pi/2})
+(\lambda+\mu^2){p}(e^{-{\mathrm{i}}{{\epsilon}}\pi/2})
\big)
    \\\equiv&\;
-{{\lfloor\![\mbox{-}1]\!\rfloor}}%
               ^2
\mu
\big(
\mu\,{p}( e^{{\mathrm{i}}{{\epsilon}}\pi/2})
-{q}( e^{{\mathrm{i}}{{\epsilon}}\pi/2})
-{r}(e^{{\mathrm{i}}{{\epsilon}}\pi/2})
\big)
    \\&
+{{\lfloor0\rfloor}}
\big(
{{\lfloor2\rfloor}}
\mu\,{q}(e^{-{\mathrm{i}}{{\epsilon}}\pi/2})
+{s}(e^{-{\mathrm{i}}{{\epsilon}}\pi/2})
+{{\lfloor\raisebox{0.1ex}{-}1\rfloor}}
 \mu^2{p}(e^{-{\mathrm{i}}{{\epsilon}}\pi/2})
\big)
    \\&
-
 e^{{\mathrm{i}}{{\epsilon}}{\ell}\pi}(\lambda+\mu^2){\mbox{$\mathstrut^{C\!}\!\Delta_{p}$}}(e^{{\mathrm{i}}{{\epsilon}}\pi/2})
-
 e^{{\mathrm{i}}{{\epsilon}}{\ell}\pi}  {{\lfloor\![\mbox{-}1]\!\rfloor}}
\,\mu\,{\mbox{$\mathstrut^{C\!}\!\Delta_{r}$}}(e^{{\mathrm{i}}{{\epsilon}}\pi/2}).
\end{aligned}
\end{equation}
Here we use the following auxiliary abbreviations:
\begin{equation}
\label{eq:akas}
 {{\lfloor\raisebox{0.1ex}{-}1\rfloor}}
   =e^{{\mathrm{i}}{{\epsilon}}\pi},\;
{{\lfloor\![\mbox{-}1]\!\rfloor}}
   =e^{-{\mathrm{i}}{{\epsilon}}\pi},\;
 {{\lfloor0\rfloor}}
   =1+e^{{\mathrm{i}}{{\epsilon}}\pi},\;
 {{\lfloor2\rfloor}}
 =1-e^{{\mathrm{i}}{{\epsilon}}\pi}.
\end{equation}
The expressions
${\mbox{$\mathstrut^{C\!}\!\Delta_{{{\mbox{\tiny\ding{74}}}}}$}}$,
where {{\mbox{\small\ding{74}}}}$\;\in\{{p},{q},{r},{s}\}$,
were
introduced in the proof of Theorem \ref{t:020}.
They denote the differences of the left- and right-hand sides of
{{Eq.s}}~\eqref{eq:pppC}-\eqref{eq:sssC}.

In each of the above four pairs
of equalities
the first ones 
are merely the expansions of the
corresponding definitions
\eqref{eq:580} with regard to the particular value of $z$ picked out above.
On the contrary,
the second equalities are ``the genuine identities'' in which the right-hand sides
represent
some
rearrangements of the left-hand ones
whose
several constituents are aggregated 
to the expressions
${\mbox{$\mathstrut^{C\!}\!\Delta_{{{\mbox{\tiny\ding{74}}}}}$}}$.
Thus {{Eq.s}}~\eqref{eq:difBp}-\eqref{eq:difBs}
express the coincidences, upon simplification,  of some linear combinations
of arbitrary fixed functions
${p},{q},{r},{s}$ evaluated at $z= e^{{\mathrm{i}}{{\epsilon}}\pi/2} $
and at $z= e^{-{\mathrm{i}}{{\epsilon}}\pi/2}$.

If ${{\epsilon}}=0$ then the argument of all the {\mbox{$pqr\!s$}}-functions and the
expressions
$
{\raisebox{0.95ex}{$\raisebox{-0.2ex}[0ex][0ex]{\scriptsize$\epsilon$}
\atop\raisebox{-0.2ex}[0ex][0ex]{$\mathstrut^{B\hspace{-0.2ex}}\hspace{-0.5ex}
\Delta_{{{\mbox{\tiny\ding{74}}}}}$}$}}
$
considered as the functions of $z$
 is $+1$.
Let $ {{\epsilon}} $ be further varied through the segment $[0,1]$.
Then the arguments of the functions involved in {{Eq.s}}~\eqref{eq:difBp}-\eqref{eq:difBs}
move along the
circular arcs, either clockwise of counter-clockwise.
The values the functions assume
thereat can be regarded as the result of their analytic continuation
from 
the vicinity of $+1$.
At end points of the noted arcs corresponding to ${{\epsilon}}=1$ the arguments of
the functions become either $  e^{{\mathrm{i}}\pi/2}={\mathrm{i}} $ or  $  e^{-{\mathrm{i}}\pi/2}=-{\mathrm{i}} $
while the expressions
$
{\raisebox{0.95ex}{$\raisebox{-0.2ex}[0ex][0ex]{\scriptsize$\epsilon$}
\atop\raisebox{-0.2ex}[0ex][0ex]{$\mathstrut^{B\hspace{-0.2ex}}\hspace{-0.5ex}
\Delta_{{{\mbox{\tiny\ding{74}}}}}$}$}}
$
on the left
turn into
$
{\mbox{$\mathstrut^{B\!}\!\Delta_{{{\mbox{\tiny\ding{74}}}}}$}}
=\lim_{{{{\epsilon}}\nearrow 1}}
{\raisebox{0.95ex}{$\raisebox{-0.2ex}[0ex][0ex]{\scriptsize$\epsilon$}
\atop\raisebox{-0.2ex}[0ex][0ex]{$\mathstrut^{B\hspace{-0.2ex}}\hspace{-0.5ex}
\Delta_{{{\mbox{\tiny\ding{74}}}}}$}$}}
$
evaluated at $-{\mathrm{i}}$.
Besides, it holds 
$
 {{\lfloor\raisebox{0.1ex}{-}1\rfloor}}
=
{{\lfloor\![\mbox{-}1]\!\rfloor}}
=
-1,\;
 {{\lfloor0\rfloor}}
=0,\;
 {{\lfloor2\rfloor}}
=2$ thereat.

Taking all these simplifications into account, one finds
that in the particular case under consideration
the equalities of the first and the last expressions
in each of the formulas
{{Eq.s}}~\eqref{eq:difBp}-\eqref{eq:difBs}
combine 
to {{Eq.s}}~\eqref{eq:240}.

\section{Identities utilized in the proof of Theorem \ref{t:040} }\label{a:050}

The following two
identities
verifiable by straightforward computation
hold true for arbitrary holomorphic
functions ${p}, {q}, {r}, {s} $.
%
%
%
%
\begin{equation} \label{eq:Adelta}
\hspace{-2.5em}
    \begin{aligned}
e^{-2{\mathrm{i}}{{\epsilon}}\pi}
\raisebox{-1.0ex}{$\Bigl\lfloor $}
\raisebox{+0.0ex}[0ex][2.2ex]{$\mathstrut$}_{  z\leftleftharpoons e^{{\mathrm{i}}{{\epsilon}}\pi}/z }
\hspace{-3.2em}
{\hspace{0.08ex}{\mathfrak D}}\hspace{0.5ex}
\rfloor
\equiv & \;
  \lfloor\,{\hspace{0.08ex}{\mathfrak D}}\,\rfloor
+
e^{-{\mathrm{i}}{{\epsilon}}{\ell}\pi}z^{2({\ell}-1)}
\bigl(
\!
{\raisebox{0.8ex}{\!$ \raisebox{1ex}[0.0ex][0.0ex]%
 {$\mbox{\raisebox{-2.1ex}[1ex][0.0ex]{\scalebox{0.6}{$\hspace{1.3ex}{{\epsilon}}$} }}
 \atop\rotatebox{180}{\scalebox{0.7}[0.7]{$\curvearrowleft$}}$}\atop
 \mbox{$ {{s}} $}$}\!}
\!(1/z) %
{\raisebox{0.95ex}{$\raisebox{-0.2ex}[0ex][0ex]{\scriptsize$\epsilon$}
\atop\raisebox{-0.2ex}[0ex][0ex]{$\mathstrut^{A\hspace{-0.1ex}}\hspace{-0.5ex}\Delta_{p}$}$}}
\!(z)
-
\!
{\raisebox{0.8ex}{\!$ \raisebox{1ex}[0.0ex][0.0ex]%
 {$\mbox{\raisebox{-2.1ex}[1ex][0.0ex]{\scalebox{0.6}{$\hspace{1.3ex}{{\epsilon}}$} }}
 \atop\rotatebox{180}{\scalebox{0.7}[0.7]{$\curvearrowleft$}}$}\atop
 \mbox{$ {{r}} $}$}\!}
\!(1/z) %
{\raisebox{0.95ex}{$\raisebox{-0.2ex}[0ex][0ex]{\scriptsize$\epsilon$}
\atop\raisebox{-0.2ex}[0ex][0ex]{$\mathstrut^{A\hspace{-0.1ex}}\hspace{-0.5ex}\Delta_{q}$}$}}
\!(z)
\bigr)
\\[-2.0ex]
      & \hphantom {   \;\lfloor\,{\hspace{0.08ex}{\mathfrak D}}\,\rfloor             }
-\bigl({r}(z)+\mu\,z^{-2}{p}(z)\bigr)
{\raisebox{0.95ex}{$\raisebox{-0.2ex}[0ex][0ex]{\scriptsize$\epsilon$}
\atop\raisebox{-0.2ex}[0ex][0ex]{$\mathstrut^{A\hspace{-0.1ex}}\hspace{-0.5ex}\Delta_{r}$}$}}
\!(z)%
-
{p}(z)
{\raisebox{0.95ex}{$\raisebox{-0.2ex}[0ex][0ex]{\scriptsize$\epsilon$}
\atop\raisebox{-0.2ex}[0ex][0ex]{$\mathstrut^{A\hspace{-0.1ex}}\hspace{-0.5ex}\Delta_{s}$}$}}
\!(z)
,
    \end{aligned}
\end{equation}
\begin{equation}
                \label{eq:Bdelta}
    \begin{aligned}
e^{-2{\mathrm{i}}{{\epsilon}}\pi}
\raisebox{-1.1ex}{$\Bigl\lfloor$}
\raisebox{-1.3ex}[0ex][2.4ex]{$\mathstrut$}_{  z\leftleftharpoons e^{{\mathrm{i}}{{\epsilon}}\pi} z }
\hspace{-2.9em}
{\hspace{0.08ex}{\mathfrak D}}\;
\rfloor
\equiv &
 \;\lfloor\,{\hspace{0.08ex}{\mathfrak D}}\,\rfloor
+ e^{-{\mathrm{i}}{{\epsilon}}{\ell}\pi}
z^{2(1-{\ell})}
\big\lgroup
e^{-{\mathrm{i}}{{\epsilon}}{\ell}\pi}
(\!
{\raisebox{0.8ex}{\!$ \raisebox{1ex}[0.0ex][0.0ex]%
 {$\mbox{\raisebox{-2.1ex}[1ex][0.0ex]{\scalebox{0.6}{$\hspace{1.3ex}{{\epsilon}}$} }}
 \atop\rotatebox{180}{\scalebox{0.7}[0.7]{$\curvearrowleft$}}$}\atop
 \mbox{$ {{s}} $}$}\!}
\!\!
\,
{\raisebox{0.95ex}{$\raisebox{-0.2ex}[0ex][0ex]{\scriptsize$\epsilon$}
\atop\raisebox{-0.2ex}[0ex][0ex]{$\mathstrut^{B\hspace{-0.2ex}}\hspace{-0.5ex}\Delta_{p}$}$}}
-
\!
{\raisebox{0.8ex}{\!$ \raisebox{1ex}[0.0ex][0.0ex]%
 {$\mbox{\raisebox{-2.1ex}[1ex][0.0ex]{\scalebox{0.6}{$\hspace{1.3ex}{{\epsilon}}$} }}
 \atop\rotatebox{180}{\scalebox{0.7}[0.7]{$\curvearrowleft$}}$}\atop
 \mbox{$ {{r}} $}$}\!}
\!\!
\,
{\raisebox{0.95ex}{$\raisebox{-0.2ex}[0ex][0ex]{\scriptsize$\epsilon$}
\atop\raisebox{-0.2ex}[0ex][0ex]{$\mathstrut^{B\hspace{-0.2ex}}\hspace{-0.5ex}\Delta_{q}$}$}}
)
+(\mu\,z^2{p} + {q})
{\raisebox{0.95ex}{$\raisebox{-0.2ex}[0ex][0ex]{\scriptsize$\epsilon$}
\atop\raisebox{-0.2ex}[0ex][0ex]{$\mathstrut^{B\hspace{-0.2ex}}\hspace{-0.5ex}\Delta_{r}$}$}}
\\[-1.5ex]
&     \hphantom{ \;\lfloor\,{\hspace{0.08ex}{\mathfrak D}}\,\rfloor  + e^{-{\mathrm{i}}{{\epsilon}}\pi} z^{2(1-{\ell})}\big\lgroup \!\!\!\!\! }
+
(\lambda+\mu^2)^{-1}
(\mu\,z^2{r} + {s})
(\mu\,z^2
{\raisebox{0.95ex}{$\raisebox{-0.2ex}[0ex][0ex]{\scriptsize$\epsilon$}
\atop\raisebox{-0.2ex}[0ex][0ex]{$\mathstrut^{B\hspace{-0.2ex}}\hspace{-0.5ex}\Delta_{r}$}$}}
+
{\raisebox{0.95ex}{$\raisebox{-0.2ex}[0ex][0ex]{\scriptsize$\epsilon$}
\atop\raisebox{-0.2ex}[0ex][0ex]{$\mathstrut^{B\hspace{-0.2ex}}\hspace{-0.5ex}\Delta_{s}$}$}}
 )
 \big\rgroup.
    \end{aligned}
\end{equation}
Here $ {{\epsilon}} \in [-1,1]$ is the auxiliary real parameter.
$ \lfloor{\hspace{0.08ex}{\mathfrak D}}\rfloor $
denotes the right-hand side of {{Eq.}}{~}\eqref{delTa}. 
{{Eq.s}}~\eqref{eq:490} play
role of definitions of the symbols
$
{\raisebox{0.95ex}{$\raisebox{-0.2ex}[0ex][0ex]{\scriptsize$\epsilon$}
\atop\raisebox{-0.2ex}[0ex][0ex]{$\mathstrut^{A\hspace{-0.1ex}}\hspace{-0.5ex}
\Delta_{{{\mbox{\tiny\ding{74}}}}}$}$}}
,$ where
{{\mbox{\small\ding{74}}}} $\in\{ {p},{q},{r},{s} \}.$
Similarly,
{{Eq.s}}~\eqref{eq:580} explain 
the meaning of the symbols
$
{\raisebox{0.95ex}{$\raisebox{-0.2ex}[0ex][0ex]{\scriptsize$\epsilon$}
\atop\raisebox{-0.2ex}[0ex][0ex]{$\mathstrut^{B\hspace{-0.2ex}}\hspace{-0.5ex}
\Delta_{{{\mbox{\tiny\ding{74}}}}}$}$}}
$. 
The 
renderings 
of all these abbreviations, as they stand,
are considered as the functions of $z$.
Here
we also employ in recording some tricks 
allowing 
somewhat more compact
presentation
of formulas than in the preceding Appendices. 
In particular,  `the accent' 
\raisebox{2.3ex}[0ex][0ex]{$
 \mbox{\raisebox{-2.1ex}[1ex][0.0ex]{\scalebox{0.7}{$\hspace{0.9ex}{{\epsilon}}$} }}
 \atop
 \rotatebox{180}{\scalebox{0.9}[0.9]{$\curvearrowleft$}}
 $}
denotes the transformation of rotation of the function argument
at an angle ${{\epsilon}}\pi$, i.e.\
$ 
{\raisebox{0.8ex}{\!$ \raisebox{1.3ex}[0.0ex][0.0ex]%
 {$\mbox{\raisebox{-2.1ex}[1ex][0.0ex]{\scalebox{0.6}{$\hspace{1.3ex}{{\epsilon}}$} }}
 \atop\rotatebox{180}{\scalebox{0.7}[0.7]{$\curvearrowleft$}}$}\atop
 \mbox{$  {{\mbox{\small\ding{74}}}}  $}$}\!}
\!(z)^{\vphantom{I^I}}\!
= {{\mbox{\small\ding{74}}}}(e^{{\mathrm{i}}{{\epsilon}}\pi} z)$.
Note that the arguments of functions are  displayed in {{Eq.}}{~}\eqref{eq:Adelta}
(except for $ \lfloor{\hspace{0.08ex}{\mathfrak D}}\rfloor $) but they are suppressed 
in {{Eq.}}~\eqref{eq:Bdelta}
because in the latter case
the arguments of all the functions
coincide and are
equal to $z$.

To guarantee the meaningfulness of the formulas \eqref{eq:Adelta} and \eqref{eq:Bdelta}
one has to ensure the belonging of the values of arguments, for which
the functions involved in them are evaluated, to the appropriate domain.
These values depend on ${{\epsilon}}$.
In particular, if $ {{\epsilon}} = 0 $ then all the functions are evaluated at
either  $z$ or  $1/z$. In such a case
one may get any $z\in{\mbox{$\mathstrut^\backprime{}\mathbb{C}^*$}}$
for which 
both formulas
\eqref{eq:Adelta}, \eqref{eq:Bdelta} prove to be correctly defined ---
and the equalities they represent hold true. 
Further, starting
with
${{\epsilon}} = 0$,
we
carry out analytic continuations of all the
constituents of {{Eq.s}}{~}\eqref{eq:Adelta} and \eqref{eq:Bdelta}
varying $ {{\epsilon}}\in[0,1] $ from 0 to 1. In the limit $ {{{\epsilon}}\nearrow 1} $
(i.e.\ at the end point of the arc of analytic continuation)
the expressions denoted
$
{\raisebox{0.95ex}{$\raisebox{-0.2ex}[0ex][0ex]{\scriptsize$\epsilon$}
\atop\raisebox{-0.2ex}[0ex][0ex]{$\mathstrut^{A\hspace{-0.1ex}}\hspace{-0.5ex}
\Delta_{{{\mbox{\tiny\ding{74}}}}}$}$}}
$
and
$
{\raisebox{0.95ex}{$\raisebox{-0.2ex}[0ex][0ex]{\scriptsize$\epsilon$}
\atop\raisebox{-0.2ex}[0ex][0ex]{$\mathstrut^{B\hspace{-0.2ex}}\hspace{-0.5ex}
\Delta_{{{\mbox{\tiny\ding{74}}}}}$}$}}
$
become identical
to the expressions
${\mbox{$\mathstrut^{A}\!\Delta_{{{\mbox{\tiny\ding{74}}}}}$}}$
and
$ {\mbox{$\mathstrut^{B\!}\!\Delta_{{{\mbox{\tiny\ding{74}}}}}$}}$,
respectively (see the discussion following {{Eq.}}~\eqref{eq:300} and {{Eq.}}~\eqref{eq:320}),
while
the 
transformation
indicated by `the accent'
\raisebox{2.1ex}[0ex][0ex]{$
 \mbox{\raisebox{-2.1ex}[1ex][0.0ex]{\scalebox{0.9}{$\hspace{0.9ex}{{\epsilon}}$} }}
 \atop
 \rotatebox{180}{\scalebox{0.9}[0.9]{$\curvearrowleft$}}$}%
converts to
the semi-monodromy transformation denoted earlier by `the accent'
\raisebox{2ex}{\rotatebox{185}{\scalebox{0.9}[1.0]{$\curvearrowleft$}}}
(indicating 
application 
of the operator ${{{\mathcal{M}}}^{1/2}}$).

Inspecting the result of the outlined analytic continuation along the
image of the segment $[0,1]\in{{\epsilon}}$, one finds
that this is nothing else but
the equations \eqref{eq:300} and \eqref{eq:320},
provided that
$z$ belongs to the vicinity of $+{\mathrm{i}}$ 
in the former case
and $\Im z<0$ in the latter one.

\section{Identities utilized in the proof of Theorem \ref{t:050} }\label{a:060}

The following
identity,
which is
verifiable by straightforward computation,
holds true for arbitrary holomorphic
functions ${{E}}, {p}, {q}, {r}, {s} $.
  \begin{equation}
               \label{eq:HoL_A}
  \begin{aligned}
\raisebox{-0.9ex}{$\Bigl\lfloor$}
\raisebox{-1.3ex}[0ex][1.7ex]{$\mathstrut$}_{  z\leftleftharpoons e^{{\mathrm{i}}{{\epsilon}}\pi}/z }
\hspace{-2.9em}
e^{-\mu(z+1/z)}
({{\mathcal H}}
\circ\!
{\raisebox{1.0ex}{\strut}\raisebox{0.81ex}[0ex][0ex]{${\mbox{\scriptsize
$ {\epsilon} $%
\;\;\;\;}\atop\mbox{$\mathrm{L}_{ {A} }$}}$}}
\!)
[{{E}}\hspace{0.3ex}]{{{}}}
\raisebox{-0.2ex}{$\bigr\rfloor$}
                  &\equiv\;
z^2{p}{{{}}}\,{{\mathcal H}}'[ {{E}}]{{{}}}
  \\[-1.2ex]
 &\relax      \hspace{-12ex}
+{{\lfloor\![\mbox{-}1]\!\rfloor}}
\big( 
\mu(
{{\lfloor\raisebox{0.1ex}{-}1\rfloor}}
{{\lfloor2\rfloor}}
+
{{\lfloor0\rfloor}}
z^2){p}{{{}}}
-
{{\lfloor\raisebox{0.1ex}{-}1\rfloor}}
{q}{{{}}}
+2{{\lfloor\raisebox{0.1ex}{-}1\rfloor}}
z^2{r}{{{}}}
\big)
{{\mathcal H}}[ {{E}}]{{{}}}
  \\&\relax\hspace{-12ex}
+ z^2{{E'}}{{{}}}
 \Delta_{{p}}\!\!'{{{}}}\,\,
+ z^2{{E}}{{{}}}
 \Delta_{{q}}\!\!'{{{}}}\,\,
  \\
  &\relax\hspace{-12ex}
+{{\lfloor\![\mbox{-}1]\!\rfloor}}
\big(
{{\lfloor\raisebox{0.1ex}{-}1\rfloor}}
(\lambda-({\ell}+1)\mu{}z){{E}}{{{}}}
+2
{{\lfloor\raisebox{0.1ex}{-}1\rfloor}}
  z^2{{E''}}{{{}}}
  \\
&\relax\hspace{-10.2ex} \hphantom{-{{\lfloor\![\mbox{-}1]\!\rfloor}}\big(  }
+(\mu(
     {{\lfloor0\rfloor}}
     {{\lfloor2\rfloor}}
                      -1+z^2)
+
{{\lfloor\raisebox{0.1ex}{-}1\rfloor}}
(\mu + 2z)){{E'}}{{{}}}
\big)
 \Delta_{{p}}{{{}}}
  \\
& \relax\hspace{-24ex}
+{{\lfloor\![\mbox{-}1]\!\rfloor}}
\big(
\big(\mu(
{{\lfloor0\rfloor}}
{{\lfloor2\rfloor}}
 -1 + z^2)
-
{{\lfloor\raisebox{0.1ex}{-}1\rfloor}}
\mu
({\ell} - 1 -\mu{}z)
\big){{E}}{{{}}}
+
{{\lfloor\raisebox{0.1ex}{-}1\rfloor}}
z^2{{E'}}{{{}}}
\big)
\Delta_{{q}}{{{}}}
  \\&
+z^2{{E'}}{{{}}}
\Delta_{{r}}{{{}}}
+{{E}}{{{}}}
\Delta_{{s}}{{{}}}
  \\&
-{{\lfloor0\rfloor}}
 \mu
\big(
(\mu+
{{\lfloor\![\mbox{-}1]\!\rfloor}}
  {}z({\ell}-1-\mu{}z))
(z^2{p}{{{}}}{{E'}}{{{}}}+{q}{{{}}}{{E}}{{{}}})
   \\&\relax\hspace{-11ex}
\hphantom{  +{{\lfloor0\rfloor}}
              \mu  \big(   }\hspace{4.8em}
+(1- {{\lfloor\![\mbox{-}1]\!\rfloor}}
   {}z^2)
(z^2{r}{{{}}}{{E'}}{{{}}}+{s}{{{}}}{{E}}{{{}}})
\big).
\end{aligned}
\end{equation}
Here ${{\epsilon}}\in[-1,1]$ is the auxiliary real parameter.
The abbreviations
$ 
 {{\lfloor\raisebox{0.1ex}{-}1\rfloor}}
  , 
 {{\lfloor\![\mbox{-}1]\!\rfloor}}
, 
 {{\lfloor0\rfloor}}
, 
 {{\lfloor2\rfloor}}
$ 
are to be expanded in accordance with
formulas \eqref{eq:akas}.
The operator ${{\mathcal H}}$ is defined by {{Eq.}}~\eqref{sDCHop} (see also {{Eq.}}{~}\eqref{sDCHE}),
${{\mathcal H}}'=d/d z\, \circ\, {{\mathcal H}}  $,
the operator
{\raisebox{1.0ex}{\strut}\raisebox{0.81ex}[0ex][0ex]{${\mbox{\scriptsize
$ {\epsilon} $%
\;\;\;\;}\atop\mbox{$\mathrm{L}_{ {A} }$}}$}}
is defined by
{{Eq.}}~\eqref{eq:eopA}.
The symbols
$
\;\Delta{{{\mbox{\tiny\ding{74}}}}}(z)$, where {{\mbox{\small\ding{74}}}}$\;\in\{{p},{q},{r},{s}\}
$,
denote
the differences of the left- and right-hand sides of {{Eq.s}}{~}%
\eqref{eq:ppp'},
\eqref{eq:qqq'},
\eqref{eq:rrr'},
\eqref{eq:sss'}, respectively, considered as the functions of $z$.

Similar identity describing this time
the composition 
of ${{\mathcal H}}$ with the operator $\!
{\raisebox{1.0ex}{\strut}\raisebox{0.81ex}[0ex][0ex]{${\mbox{\scriptsize
$ {\epsilon} $%
\;\;\;\;}\atop\mbox{$\mathrm{L}_{ {B} }$}}$}}
$
(see {{Eq.}}~\eqref{eq:eopB})
reads
\begin{equation}
               \label{eq:HoL_B}
\begin{aligned}
z^{{\ell}-1}e^{-\mu(z+1/z)}
({{\mathcal H}}
\circ\!
{\raisebox{1.0ex}{\strut}\raisebox{0.81ex}[0ex][0ex]{${\mbox{\scriptsize
$ {\epsilon} $%
\;\;\;\;}\atop\mbox{$\mathrm{L}_{ {B} }$}}$}}
\!)
[  {{E}}\,]{{}}
  \equiv&\;
z^2
{\raisebox{0.8ex}{\!$ \raisebox{0.9ex}[0.0ex][0.0ex]%
 {$\mbox{\raisebox{-2.1ex}[1ex][0.0ex]{\scalebox{0.6}{$\hspace{1.3ex}{{\epsilon}}$} }}
 \atop\rotatebox{180}{\scalebox{0.7}[0.7]{$\curvearrowleft$}}$}\atop
 \mbox{$ {{r}} $}$}\!}
\!\!  \raisebox{0.2ex}{
{\raisebox{0.8ex}{\!$ \raisebox{1.2ex}[0.0ex][0.0ex]%
 {$\mbox{\raisebox{-2.1ex}[1ex][0.0ex]{\scalebox{0.6}{$\hspace{1.3ex}{{\epsilon}}$} }}
 \atop\rotatebox{180}{\scalebox{0.7}[0.7]{$\curvearrowleft$}}$}\atop
 \mbox{$ {{\mathcal H}}' $}$}\!}
}\hspace{-0.7ex}[{{E}}\hspace{0.2ex}]
+\left(
2({\ell} - 1)z
\!
{\raisebox{0.8ex}{\!$ \raisebox{0.9ex}[0.0ex][0.0ex]%
 {$\mbox{\raisebox{-2.1ex}[1ex][0.0ex]{\scalebox{0.6}{$\hspace{1.3ex}{{\epsilon}}$} }}
 \atop\rotatebox{180}{\scalebox{0.7}[0.7]{$\curvearrowleft$}}$}\atop
 \mbox{$ {{r}} $}$}\!}
\!
+
\!
{\raisebox{0.8ex}{\!$ \raisebox{0.9ex}[0.0ex][0.0ex]%
 {$\mbox{\raisebox{-2.1ex}[1ex][0.0ex]{\scalebox{0.6}{$\hspace{1.3ex}{{\epsilon}}$} }}
 \atop\rotatebox{180}{\scalebox{0.7}[0.7]{$\curvearrowleft$}}$}\atop
 \mbox{$ {{s}} $}$}\!}
\!
+2z^2
 \!\! \raisebox{0.13ex}{  
 {\raisebox{0.8ex}{\!$ \raisebox{1.1ex}[0.0ex][0.0ex]%
 {$\mbox{\raisebox{-2.1ex}[1ex][0.0ex]{\scalebox{0.6}{$\hspace{1.3ex}{{\epsilon}}$} }}
 \atop\rotatebox{180}{\scalebox{0.7}[0.7]{$\curvearrowleft$}}$}\atop
 \mbox{$ {r}\rlap{$\mathstrut'$} $}$}\!}
 }\!\!
\right) 
 \!\!\! \raisebox{0.0ex}{  
 {\raisebox{0.8ex}{\!$ \raisebox{1.4ex}[0.0ex][0.0ex]%
 {$\mbox{\raisebox{-2.1ex}[1ex][0.0ex]{\scalebox{0.6}{$\hspace{1.3ex}{{\epsilon}}$} }}
 \atop\rotatebox{180}{\scalebox{0.7}[0.7]{$\curvearrowleft$}}$}\atop
 \mbox{$ {{\mathcal H}} $}$}\!}
 }\!\![{{E}}\hspace{0.2ex}] %
              \\[-2.5ex]  
 &
+z^2
\!\! \raisebox{0.06ex}{ 
{\raisebox{0.8ex}{\!$ \raisebox{1.4ex}[0.0ex][0.0ex]%
 {$\mbox{\raisebox{-2.1ex}[1ex][0.0ex]{\scalebox{0.6}{$\hspace{1.3ex}{{\epsilon}}$} }}
 \atop\rotatebox{180}{\scalebox{0.7}[0.7]{$\curvearrowleft$}}$}\atop
 \mbox{$ \!{{E}}\rlap{$\mathstrut'$} $}$}\!}
}  \!\!\!
\!\raisebox{0.06ex}{
  {\raisebox{0.8ex}{\!$ \raisebox{1.4ex}[0.0ex][0.0ex]%
 {$\mbox{\raisebox{-2.1ex}[1ex][0.0ex]{\scalebox{0.6}{$\hspace{1.3ex}{{\epsilon}}$} }}
 \atop\rotatebox{180}{\scalebox{0.7}[0.7]{$\curvearrowleft$}}$}\atop
 \mbox{$ \,\!\Delta_{{r}}\!\!\rlap{$\mathstrut'$} $}$}\!}
  } \! %
+
\!\!\raisebox{0.00ex}{ 
 {\raisebox{0.8ex}{\!$ \raisebox{1.4ex}[0.0ex][0.0ex]%
 {$\mbox{\raisebox{-2.1ex}[1ex][0.0ex]{\scalebox{0.6}{$\hspace{1.3ex}{{\epsilon}}$} }}
 \atop\rotatebox{180}{\scalebox{0.7}[0.7]{$\curvearrowleft$}}$}\atop
 \mbox{$ \!{{E}} $}$}\!}
 }  \!\!\!\!
  \!\raisebox{0.08ex}{
{\raisebox{0.8ex}{\!$ \raisebox{1.2ex}[0.0ex][0.0ex]%
 {$\mbox{\raisebox{-2.1ex}[1ex][0.0ex]{\scalebox{0.6}{$\hspace{1.3ex}{{\epsilon}}$} }}
 \atop\rotatebox{180}{\scalebox{0.7}[0.7]{$\curvearrowleft$}}$}\atop
 \mbox{$ \,\!\Delta_{{s}}\!\!\rlap{$\mathstrut'$} $}$}\!}
} \! %
-(\lambda+\mu^2)
\bigl(
 \!\!\raisebox{0.06ex}{ 
 {\raisebox{0.8ex}{\!$ \raisebox{1.4ex}[0.0ex][0.0ex]%
 {$\mbox{\raisebox{-2.1ex}[1ex][0.0ex]{\scalebox{0.6}{$\hspace{1.3ex}{{\epsilon}}$} }}
 \atop\rotatebox{180}{\scalebox{0.7}[0.7]{$\curvearrowleft$}}$}\atop
 \mbox{$ \!{{E}}\rlap{$\mathstrut'$} $}$}\!}
 }  %
\!\!\!\!\!\!\!\raisebox{0.00ex}{
 {\raisebox{0.8ex}{\!$ \raisebox{1.4ex}[0.0ex][0.0ex]%
 {$\mbox{\raisebox{-2.1ex}[1ex][0.0ex]{\scalebox{0.6}{$\hspace{1.3ex}{{\epsilon}}$} }}
 \atop\rotatebox{180}{\scalebox{0.7}[0.7]{$\curvearrowleft$}}$}\atop
 \mbox{$ \,\,\,\,\!\Delta_{{p}} $}$}\!}
 } \! %
+
\!\!\raisebox{0.08ex}{ 
{\raisebox{0.8ex}{\!$ \raisebox{1.4ex}[0.0ex][0.0ex]%
 {$\mbox{\raisebox{-2.1ex}[1ex][0.0ex]{\scalebox{0.6}{$\hspace{1.3ex}{{\epsilon}}$} }}
 \atop\rotatebox{180}{\scalebox{0.7}[0.7]{$\curvearrowleft$}}$}\atop
 \mbox{$ \!{{E}} $}$}\!}
}  
\!\!\!\!\!\!\!\raisebox{0.05ex}{
{\raisebox{0.8ex}{\!$ \raisebox{1.4ex}[0.0ex][0.0ex]%
 {$\mbox{\raisebox{-2.1ex}[1ex][0.0ex]{\scalebox{0.6}{$\hspace{1.3ex}{{\epsilon}}$} }}
 \atop\rotatebox{180}{\scalebox{0.7}[0.7]{$\curvearrowleft$}}$}\atop
 \mbox{$ \,\,\,\,\!\Delta_{{q}} $}$}\!}
  } \!\!  %
\bigr)
\\[-2.9ex] &
-\big(
(\mu - ({\ell} - 1)z)
\!\!\raisebox{0.1ex}{ 
{\raisebox{0.8ex}{\!$ \raisebox{1.4ex}[0.0ex][0.0ex]%
 {$\mbox{\raisebox{-2.1ex}[1ex][0.0ex]{\scalebox{0.6}{$\hspace{1.3ex}{{\epsilon}}$} }}
 \atop\rotatebox{180}{\scalebox{0.7}[0.7]{$\curvearrowleft$}}$}\atop
 \mbox{$ \!\!{{E}}\rlap{$\mathstrut'$} $}$}\!}
} \!\!
+
(\lambda+({\ell}+1)\mu{}z)
 \!\!\raisebox{0.05ex}{ 
 {\raisebox{0.8ex}{\!$ \raisebox{1.4ex}[0.0ex][0.0ex]%
 {$\mbox{\raisebox{-2.1ex}[1ex][0.0ex]{\scalebox{0.6}{$\hspace{1.3ex}{{\epsilon}}$} }}
 \atop\rotatebox{180}{\scalebox{0.7}[0.7]{$\curvearrowleft$}}$}\atop
 \mbox{$ \!{{E}} $}$}\!}
 }  \!\!\! %
\big)
\!\!\!\!\raisebox{0.05ex}{
{\raisebox{0.8ex}{\!$ \raisebox{1.4ex}[0.0ex][0.0ex]%
 {$\mbox{\raisebox{-2.1ex}[1ex][0.0ex]{\scalebox{0.6}{$\hspace{1.3ex}{{\epsilon}}$} }}
 \atop\rotatebox{180}{\scalebox{0.7}[0.7]{$\curvearrowleft$}}$}\atop
 \mbox{$ \,\,\,\,\!\Delta_{{r}} $}$}\!}
  }
  \! %
\\[-2.5ex] &
+ (
\!\!\raisebox{0.1ex}{
{\raisebox{0.8ex}{\!$ \raisebox{1.4ex}[0.0ex][0.0ex]%
 {$\mbox{\raisebox{-2.1ex}[1ex][0.0ex]{\scalebox{0.6}{$\hspace{1.3ex}{{\epsilon}}$} }}
 \atop\rotatebox{180}{\scalebox{0.7}[0.7]{$\curvearrowleft$}}$}\atop
 \mbox{$ \!{{E}}\rlap{$\mathstrut'$} $}$}\!}
} \!\!
-\mu
\!\!\raisebox{0.06ex}{ 
{\raisebox{0.8ex}{\!$ \raisebox{1.4ex}[0.0ex][0.0ex]%
 {$\mbox{\raisebox{-2.1ex}[1ex][0.0ex]{\scalebox{0.6}{$\hspace{1.3ex}{{\epsilon}}$} }}
 \atop\rotatebox{180}{\scalebox{0.7}[0.7]{$\curvearrowleft$}}$}\atop
 \mbox{$ \!{{E}} $}$}\!}
}  \!\!\!
    )  %
\!\!\!\!\raisebox{0.05ex}{
  {\raisebox{0.8ex}{\!$ \raisebox{1.4ex}[0.0ex][0.0ex]%
 {$\mbox{\raisebox{-2.1ex}[1ex][0.0ex]{\scalebox{0.6}{$\hspace{1.3ex}{{\epsilon}}$} }}
 \atop\rotatebox{180}{\scalebox{0.7}[0.7]{$\curvearrowleft$}}$}\atop
 \mbox{$ \,\,\,\,\!\Delta_{{s}} $}$}\!}
  } \! %
\\ &
+{{\lfloor0\rfloor}}
    \big(
W_0
\!\!\raisebox{0.06ex}{
{\raisebox{0.8ex}{\!$ \raisebox{1.6ex}[0.0ex][0.0ex]%
 {$\mbox{\raisebox{-2.1ex}[1ex][0.0ex]{\scalebox{0.6}{$\hspace{1.3ex}{{\epsilon}}$} }}
 \atop\rotatebox{180}{\scalebox{0.7}[0.7]{$\curvearrowleft$}}$}\atop
 \mbox{$ \!{{E}} $}$}\!}
}  \!\!
+W_1
\!\raisebox{0.1ex}{ 
{\raisebox{0.8ex}{\!$ \raisebox{1.4ex}[0.0ex][0.0ex]%
 {$\mbox{\raisebox{-2.1ex}[1ex][0.0ex]{\scalebox{0.6}{$\hspace{1.3ex}{{\epsilon}}$} }}
 \atop\rotatebox{180}{\scalebox{0.7}[0.7]{$\curvearrowleft$}}$}\atop
 \mbox{$ \!{{E}}\rlap{$\mathstrut'$} $}$}\!}
} \!\!
-W_2
\!\raisebox{0.1ex}{ 
{\raisebox{0.8ex}{\!$ \raisebox{1.4ex}[0.0ex][0.0ex]%
 {$\mbox{\raisebox{-2.1ex}[1ex][0.0ex]{\scalebox{0.6}{$\hspace{1.3ex}{{\epsilon}}$} }}
 \atop\rotatebox{180}{\scalebox{0.7}[0.7]{$\curvearrowleft$}}$}\atop
 \mbox{$ \!{{E}}\rlap{$\mathstrut''$} $}$}\!}
} \!\! %
-W_3
\!\raisebox{0.1ex}{ 
{\raisebox{0.8ex}{\!$ \raisebox{1.4ex}[0.0ex][0.0ex]%
 {$\mbox{\raisebox{-2.1ex}[1ex][0.0ex]{\scalebox{0.6}{$\hspace{1.3ex}{{\epsilon}}$} }}
 \atop\rotatebox{180}{\scalebox{0.7}[0.7]{$\curvearrowleft$}}$}\atop
 \mbox{$ \!{{E}}\rlap{$\mathstrut'''$} $}$}\!}
}
 \big),
\end{aligned}
\end{equation}
where the following abbreviations are employed
\begin{eqnarray}
W_0 &=&
({\ell} + 1)(\lambda+\mu^2)\mu{}
           z  \!
           {\raisebox{0.8ex}{\!$ \raisebox{1ex}[0.0ex][0.0ex]%
 {$\mbox{\raisebox{-2.1ex}[1ex][0.0ex]{\scalebox{0.6}{$\hspace{1.3ex}{{\epsilon}}$} }}
 \atop\rotatebox{180}{\scalebox{0.7}[0.7]{$\curvearrowleft$}}$}\atop
 \mbox{$ {p} $}$}\!}
\!
\nonumber
           \\&&   
+
z\bigl(-2({\ell} - 2)\lambda + ({\ell}+1)\mu z(4
- 3{{\lfloor0\rfloor}}
+ {{\lfloor\raisebox{0.1ex}{-}1\rfloor}}%
    ^2\mu{}z) \bigr)
\!
{\raisebox{0.8ex}{\!$ \raisebox{1ex}[0.0ex][0.0ex]%
 {$\mbox{\raisebox{-2.1ex}[1ex][0.0ex]{\scalebox{0.6}{$\hspace{1.3ex}{{\epsilon}}$} }}
 \atop\rotatebox{180}{\scalebox{0.7}[0.7]{$\curvearrowleft$}}$}\atop
 \mbox{$ {r} $}$}\!}
\!
+2\mu{}z\!
{\raisebox{0.8ex}{\!$ \raisebox{1ex}[0.0ex][0.0ex]%
 {$\mbox{\raisebox{-2.1ex}[1ex][0.0ex]{\scalebox{0.6}{$\hspace{1.3ex}{{\epsilon}}$} }}
 \atop\rotatebox{180}{\scalebox{0.7}[0.7]{$\curvearrowleft$}}$}\atop
 \mbox{$ {s} $}$}\!}
\!
\nonumber\\&&    
-z^2\big(2{{\lfloor2\rfloor}}
          \lambda
        - {{\lfloor\raisebox{0.1ex}{-}1\rfloor}}
         ({\ell}+1)
         (3-{{\lfloor0\rfloor}}
         )\mu{}z\big)
\!\!\! \raisebox{0.15ex}{
{\raisebox{0.8ex}{\!$ \raisebox{1.2ex}[0.0ex][0.0ex]%
 {$\mbox{\raisebox{-2.1ex}[1ex][0.0ex]{\scalebox{0.6}{$\hspace{1.3ex}{{\epsilon}}$} }}
 \atop\rotatebox{180}{\scalebox{0.7}[0.7]{$\curvearrowleft$}}$}\atop
 \mbox{$ {r}\rlap{$\mathstrut'$} $}$}\!}
}\!\!
-\mu(1 - {{\lfloor\raisebox{0.1ex}{-}1\rfloor}}
           {}z^2)\!\! \raisebox{0.15ex}{  
{\raisebox{0.8ex}{\!$ \raisebox{1.2ex}[0.0ex][0.0ex]%
 {$\mbox{\raisebox{-2.1ex}[1ex][0.0ex]{\scalebox{0.6}{$\hspace{1.3ex}{{\epsilon}}$} }}
 \atop\rotatebox{180}{\scalebox{0.7}[0.7]{$\curvearrowleft$}}$}\atop
 \mbox{$ {s}\rlap{$\mathstrut'$} $}$}\!}
           }\!\!\!,
\nonumber
                                             \\%
W_1 &=&
-({\ell} - 1)(\lambda+\mu^2)
               z \!
{\raisebox{0.8ex}{\!$ \raisebox{1.ex}[0.0ex][0.0ex]%
 {$\mbox{\raisebox{-2.1ex}[1ex][0.0ex]{\scalebox{0.6}{$\hspace{1.3ex}{{\epsilon}}$} }}
 \atop\rotatebox{180}{\scalebox{0.7}[0.7]{$\curvearrowleft$}}$}\atop
 \mbox{$ {p} $}$}\!}
               \!
\nonumber
               \\
               &&              
+z\bigl(
 - 2\mu({\ell} - {{\lfloor3\rfloor}}
                          )
   + ( ({\ell} - 3)(4 - 3{{\lfloor0\rfloor}}
                                  )
                             -{{\lfloor2\rfloor}}
                               \lambda)z
\nonumber\\&&                   \hphantom{  +z\bigl(  }     
+{{\lfloor\raisebox{0.1ex}{-}1\rfloor}}
 {{\lfloor2\rfloor}}
  ({\ell} + 1)\mu{}z^2
\bigr)\!
{\raisebox{0.8ex}{\!$ \raisebox{1.ex}[0.0ex][0.0ex]%
 {$\mbox{\raisebox{-2.1ex}[1ex][0.0ex]{\scalebox{0.6}{$\hspace{1.3ex}{{\epsilon}}$} }}
 \atop\rotatebox{180}{\scalebox{0.7}[0.7]{$\curvearrowleft$}}$}\atop
 \mbox{$ {r} $}$}\!}
\!
\nonumber
                                       \\[0.ex]
                                        &&    
-\big(({\ell} - 1) z  +\mu(1 - {{\lfloor\raisebox{0.1ex}{-}1\rfloor}}
                                {}z^2) \big)
                             \!
{\raisebox{0.8ex}{\!$ \raisebox{1.ex}[0.0ex][0.0ex]%
 {$\mbox{\raisebox{-2.1ex}[1ex][0.0ex]{\scalebox{0.6}{$\hspace{1.3ex}{{\epsilon}}$} }}
 \atop\rotatebox{180}{\scalebox{0.7}[0.7]{$\curvearrowleft$}}$}\atop
 \mbox{$ {s} $}$}\!}
                             \!
\,\,-\,
{{\lfloor2\rfloor}}
  (\lambda + \mu^2)z^2
                            \!\!\! \raisebox{0.15ex}{  
{\raisebox{0.8ex}{\!$ \raisebox{1.2ex}[0.0ex][0.0ex]%
 {$\mbox{\raisebox{-2.1ex}[1ex][0.0ex]{\scalebox{0.6}{$\hspace{1.3ex}{{\epsilon}}$} }}
 \atop\rotatebox{180}{\scalebox{0.7}[0.7]{$\curvearrowleft$}}$}\atop
 \mbox{$ {p}\rlap{$\mathstrut'$} $}$}\!}
                            }\!\!
\nonumber\\&&\hspace{-0ex}
-z^2\big({{\lfloor\raisebox{0.1ex}{-}1\rfloor}}
         (6 + {{\lfloor\raisebox{0.1ex}{-}1\rfloor}}
           ({\ell} - 7))z
+\mu(2{{\lfloor2\rfloor}}
      +{{\lfloor\raisebox{0.1ex}{-}1\rfloor}}%
         ^2(1 - z^2))
    \big) \!\! \raisebox{0.15ex}{  
{\raisebox{0.8ex}{\!$ \raisebox{1.1ex}[0.0ex][0.0ex]%
 {$\mbox{\raisebox{-2.1ex}[1ex][0.0ex]{\scalebox{0.6}{$\hspace{1.3ex}{{\epsilon}}$} }}
 \atop\rotatebox{180}{\scalebox{0.7}[0.7]{$\curvearrowleft$}}$}\atop
 \mbox{$ {r}\rlap{$\mathstrut'$} $}$}\!}
    }\!\!
\nonumber\\&&  \hspace{-0ex}
-{{\lfloor2\rfloor}}
  {}z^2\!\! \raisebox{0.08ex}{  
{\raisebox{0.8ex}{\!$ \raisebox{1.2ex}[0.0ex][0.0ex]%
 {$\mbox{\raisebox{-2.1ex}[1ex][0.0ex]{\scalebox{0.6}{$\hspace{1.3ex}{{\epsilon}}$} }}
 \atop\rotatebox{180}{\scalebox{0.7}[0.7]{$\curvearrowleft$}}$}\atop
 \mbox{$ {s}\mathstrut' $}$}\!}
  }\!\!
- {{\lfloor2\rfloor}}
   {{\lfloor\raisebox{0.1ex}{-}1\rfloor}}%
    ^2{}z^4\!\! \raisebox{0.10ex}{  
{\raisebox{0.8ex}{\!$ \raisebox{1.2ex}[0.0ex][0.0ex]%
 {$\mbox{\raisebox{-2.1ex}[1ex][0.0ex]{\scalebox{0.6}{$\hspace{1.3ex}{{\epsilon}}$} }}
 \atop\rotatebox{180}{\scalebox{0.7}[0.7]{$\curvearrowleft$}}$}\atop
 \mbox{$ {r}\rlap{$\mathstrut''$} $}$}\!}
    }\!\!,
\nonumber
                                              \\
W_2 &=&
2{{\lfloor2\rfloor}}
 {{\lfloor\raisebox{0.1ex}{-}1\rfloor}}%
 ^2{}z^4\!\! \raisebox{0.08ex}{  
 {\raisebox{0.8ex}{\!$ \raisebox{1.2ex}[0.0ex][0.0ex]%
 {$\mbox{\raisebox{-2.1ex}[1ex][0.0ex]{\scalebox{0.6}{$\hspace{1.3ex}{{\epsilon}}$} }}
 \atop\rotatebox{180}{\scalebox{0.7}[0.7]{$\curvearrowleft$}}$}\atop
 \mbox{$ {r}\rlap{$\mathstrut'$} $}$}\!}
 }\!\!
\nonumber
\\&&
+z^2\big({{\lfloor\raisebox{0.1ex}{-}1\rfloor}}
          (10 + {{\lfloor0\rfloor}}
                        ({\ell} - 7))z
+\mu\,({{\lfloor2\rfloor}}
        +{{\lfloor\raisebox{0.1ex}{-}1\rfloor}}%
          ^2(1-z^2)) \big)
\!
{\raisebox{0.8ex}{\!$ \raisebox{1.ex}[0.0ex][0.0ex]%
 {$\mbox{\raisebox{-2.1ex}[1ex][0.0ex]{\scalebox{0.6}{$\hspace{1.3ex}{{\epsilon}}$} }}
 \atop\rotatebox{180}{\scalebox{0.7}[0.7]{$\curvearrowleft$}}$}\atop
 \mbox{$ {r} $}$}\!}
\!,
\nonumber
                                                \\
W_3 &=&
{{\lfloor2\rfloor}}
{{\lfloor\raisebox{0.1ex}{-}1\rfloor}}%
  ^2
z^4
\!
{\raisebox{0.8ex}{\!$ \raisebox{1.ex}[0.0ex][0.0ex]%
 {$\mbox{\raisebox{-2.1ex}[1ex][0.0ex]{\scalebox{0.6}{$\hspace{1.3ex}{{\epsilon}}$} }}
 \atop\rotatebox{180}{\scalebox{0.7}[0.7]{$\curvearrowleft$}}$}\atop
 \mbox{$ {r} $}$}\!}
\!.
\nonumber
\end{eqnarray}
The meaning of the $\Delta$-symbols
and the
abbreviations
$ {{\lfloor\raisebox{0.1ex}{-}1\rfloor}}
  ,
 {{\lfloor0\rfloor}}
 ,
 {{\lfloor2\rfloor}}
   ,
$
was
explained above.
`The diacritic mark'
\raisebox{2.3ex}[2.3ex][0ex]{$
 \mbox{\raisebox{-2.1ex}[1ex][0.0ex]{\scalebox{0.8}{$\hspace{0.9ex}{{\epsilon}}$} }}
 \atop
 \rotatebox{180}{\scalebox{0.9}[0.9]{$\curvearrowleft$}}
 $},
used also in {{Eq.}}~\eqref{eq:Bdelta},
denotes the transformation carrying out the rotation of the function argument
at an angle ${{\epsilon}}\pi$.

If $ {{\epsilon}}=0 $
then
the arguments of all the functions involved in
{{Eq.s}}~\eqref{eq:HoL_A} and \eqref{eq:HoL_B}
are  either $z$ or (somewhere in \eqref{eq:HoL_A}) $1/z$.
Accordingly, for any  
$ z \in {{\mbox{$\mathstrut^\backprime{}\mathbb{C}^*$}}}$
all the constituents of
the both formulas \eqref{eq:HoL_A}, \eqref{eq:HoL_B} are well defined ---
and the equalities they signify hold true
for arbitrary functions ${{E}}, {p}, {q},{r} ,{s}  $
holomorphic in $ {{\mbox{$\mathstrut^\backprime{}\mathbb{C}^*$}}}$.
Further, 
we allow $ {{\epsilon}} $ to vary through the segment
$[0,1]$ and carry out analytic continuation
along the corresponding curves (in fact, the circular arcs)
of all the constituents of the formulas \eqref{eq:HoL_A}, \eqref{eq:HoL_B}.
Fixing the result of this analytic continuation
at the end points corresponding to ${{\epsilon}} =1$,
we obtain the two equalities
in which
some noted
abbreviations acquire the known 
numerical values as follows:
$ {{\lfloor\raisebox{0.1ex}{-}1\rfloor}}
  =
 {{\lfloor\![\mbox{-}1]\!\rfloor}}
   =-1,
 {{\lfloor0\rfloor}}
   =0,
 {{\lfloor2\rfloor}}
   =2$.
 Then it is easy to see
that we obtain in this way the equations (in fact, identities)
\eqref{eq:380} and \eqref{eq:390}.

To prevent egresses of the points of evaluation of our functions from
{\mbox{$\mathstrut^\backprime{}\mathbb{C}^*$}}{}, it is enough
to pick  $z$ from the vicinity of $+{\mathrm{i}}$
in the case  of {{Eq.}}~\eqref{eq:380}
and from the half-plane $\Im z<0$ in the case  of {{Eq.}}~\eqref{eq:390}.
The extending to
greater domains can be carried out
by means of analytic continuation.

\section{Identities utilized in the proof of Theorem \ref{t:060}}\label{a:070}

The identities given below represent the
appropriately adapted
expansions of the
iterated linear operators
$\! 
{\raisebox{1.0ex}{\strut}\raisebox{0.81ex}[0ex][0ex]{${\mbox{\scriptsize
$ {\epsilon} $%
\;\;\;\;}\atop\mbox{$\mathrm{L}_{ {A} }$}}$}}
\!$ and $ \!
{\raisebox{1.0ex}{\strut}\raisebox{0.81ex}[0ex][0ex]{${\mbox{\scriptsize
$ {\epsilon} $%
\;\;\;\;}\atop\mbox{$\mathrm{L}_{ {B} }$}}$}}
\! $
defined by {{Eq.s}}~\eqref{eq:eopA}, \eqref{eq:eopB}.
Specifically, let ${{E}}, {p}, {q}, {r}, {s}$ denote arbitrary holomorphic
functions  and
${{\epsilon}}\in[-1,1]$ be the real parameter.
Then it can be shown by means of straightforward computations show that, at first,
\begin{equation}
    \label{eq:L_A o L_A}
\begin{aligned}
(
{\raisebox{1.0ex}{\strut}\raisebox{0.81ex}[0ex][0ex]{${\mbox{\scriptsize
$ {\epsilon} $%
\;\;\;\;}\atop\mbox{$\mathrm{L}_{ {A} }$}}$}}
\!
\circ
\!
{\raisebox{1.0ex}{\strut}\raisebox{0.81ex}[0ex][0ex]{${\mbox{\scriptsize
$ {\epsilon} $%
\;\;\;\;}\atop\mbox{$\mathrm{L}_{ {A} }$}}$}}\!
)[{{E}}\hspace{0.3ex}](z)
+
e^{{\mathrm{i}}{\ell}\pi}
\lfloor{\hspace{0.08ex}{\mathfrak D}}\rfloor\, 
{{E}} (z)
\equiv&  \;
\\[-01.0ex]   &  \hspace{-22ex}
\;\;\;\;
\big(({s}(z)+\mu{\,}z^{-2}{q}(z)) {{E}}(z)
+(z^2{r}(z)+\mu\,{p}(z)){{E'}}(z)
   \big) 
{\raisebox{0.95ex}{$\raisebox{-0.2ex}[0ex][0ex]{\scriptsize$\epsilon$}
\atop\raisebox{-0.2ex}[0ex][0ex]{$\mathstrut^{A\hspace{-0.1ex}}\hspace{-0.5ex}\Delta_{p}$}$}}   
   (z)
    \\[-0.80ex]  &  \hspace{-20ex}
+\big({q}(z){{E}}(z)+z^2{p}(z){{E'}}(z)\big)
{\raisebox{0.95ex}{$\raisebox{-0.2ex}[0ex][0ex]{\scriptsize$\epsilon$}
\atop\raisebox{-0.2ex}[0ex][0ex]{$\mathstrut^{A\hspace{-0.1ex}}\hspace{-0.5ex}\Delta_{q}$}$}}
          (z) %
                     \\  &\hspace{-20ex}
+
\!
{\raisebox{0.8ex}{\!$ \raisebox{1.ex}[0.0ex][0.0ex]%
 {$\mbox{\raisebox{-2.1ex}[1ex][0.0ex]{\scalebox{0.6}{$\hspace{1.3ex}{{\epsilon}}$} }}
 \atop\rotatebox{180}{\scalebox{0.7}[0.7]{$\curvearrowleft$}}$}\atop
 \mbox{$ {p} $}$}\!}
\!(1/z)
\bigl({p}(z)
{{\mathcal H}}[{{E}}\hspace{0.3ex}](z)
+
{{E'}}(z)\,
\Delta_{{p}}(z)
+
{{E}}(z)\,
 \Delta_{{q}}(z)
\bigr)
                     \\  &\hspace{-20ex}
+
{{\lfloor0\rfloor}}
  \!
  {\raisebox{0.8ex}{\!$ \raisebox{1.ex}[0.0ex][0.0ex]%
 {$\mbox{\raisebox{-2.1ex}[1ex][0.0ex]{\scalebox{0.6}{$\hspace{1.3ex}{{\epsilon}}$} }}
 \atop\rotatebox{180}{\scalebox{0.7}[0.7]{$\curvearrowleft$}}$}\atop
 \mbox{$ {p} $}$}\!}
  \!(1/z) %
W_1
+
\big(e^{\mu({{\lfloor\![0]\!\rfloor}}
             {}z
             +{{\lfloor0\rfloor}}
               {}/z)}-1\big)W_2,
\end{aligned}
\end{equation}
\begin{eqnarray}
\mbox{where }
W_1&=&\big({{\lfloor2\rfloor}}
            \mu{}z^{-2}{q}(z)+{q}'(z)\big){{E}}(z)
+
z^2{p}(z){{E''}}(z)
\nonumber \\ && 
+
\big({q}(z)+({{\lfloor2\rfloor}}
              \mu +2z){p}(z)+z^2{p}'(z) \big){{E'}}(z),
\nonumber
\\
W_2&=&
-
  \!
  {\raisebox{0.8ex}{\!$ \raisebox{1.ex}[0.0ex][0.0ex]%
 {$\mbox{\raisebox{-2.1ex}[1ex][0.0ex]{\scalebox{0.6}{$\hspace{1.3ex}{{\epsilon}}$} }}
 \atop\rotatebox{180}{\scalebox{0.7}[0.7]{$\curvearrowleft$}}$}\atop
 \mbox{$ {p} $}$}\!}
  \!(1/z)   
\Big(
  \big(\mu{}z^{-2}(z^2
                   -{{\lfloor\raisebox{0.1ex}{-}1\rfloor}}%
                                 ^2){q}(z)
                       + {{\lfloor\raisebox{0.1ex}{-}1\rfloor}}
                          {q}'(z)\big){{E}}(z)
\nonumber \\[-0.5ex] &&  \hphantom{  - \! {p} 
{\raisebox{0.8ex}{\!$ \raisebox{1.ex}[0.0ex][0.0ex]%
 {$\mbox{\raisebox{-2.1ex}[1ex][0.0ex]{\scalebox{0.6}{$\hspace{1.3ex}{{\epsilon}}$} }}
 \atop\rotatebox{180}{\scalebox{0.7}[0.7]{$\curvearrowleft$}}$}\atop
 \mbox{$ {p} $}$}\!}
\!(1/z)  \big(  }
+\big(
 (\mu(z^2-{{\lfloor\raisebox{0.1ex}{-}1\rfloor}}%
                      ^2)
         + 2 {{\lfloor\raisebox{0.1ex}{-}1\rfloor}}
           {}z){p}(z)
\nonumber \\   && \hphantom{  -   \!{p} 
\!(1/z)  \big(+\big(  } \;\;\;\;\;\;\;
 + {{\lfloor\raisebox{0.1ex}{-}1\rfloor}}
    {q}(z)
+{{\lfloor\raisebox{0.1ex}{-}1\rfloor}}
  {}z^2{p}'(z)
 \big){{E'}}(z)
\nonumber \\[-0.5ex]    && \hphantom{  -   \! {p} 
\!(1/z)      \big(  }
+
{{\lfloor\raisebox{0.1ex}{-}1\rfloor}}
  {}z^2{p}(z){{E''}}(z)
\Big)
  \nonumber \\[-0.9ex] &&
+
   \!
   {\raisebox{0.8ex}{\!$ \raisebox{1.ex}[0.0ex][0.0ex]%
 {$\mbox{\raisebox{-2.1ex}[1ex][0.0ex]{\scalebox{0.6}{$\hspace{1.3ex}{{\epsilon}}$} }}
 \atop\rotatebox{180}{\scalebox{0.7}[0.7]{$\curvearrowleft$}}$}\atop
 \mbox{$ {q} $}$}\!}
   \!(1/z)  %
\big(
{q}(z)
{{E}}(z)
+
z^2{p}(z)
{{E'}}(z)
\big)
\nonumber
\end{eqnarray}
and, at second,
\\[-2em]
\begin{equation}
   \label{eq:L_B o L_B}
\begin{aligned}
e^{{\mathrm{i}}{\ell}\epsilon\pi}z^{2({\ell}-1)}
\big\lgroup
(
{\raisebox{1.0ex}{\strut}\raisebox{0.81ex}[0ex][0ex]{${\mbox{\scriptsize
$ {\epsilon} $%
\;\;\;\;}\atop\mbox{$\mathrm{L}_{ {B} }$}}$}}
\circ
\!
{\raisebox{1.0ex}{\strut}\raisebox{0.81ex}[0ex][0ex]{${\mbox{\scriptsize
$ {\epsilon} $%
\;\;\;\;}\atop\mbox{$\mathrm{L}_{ {B} }$}}$}}\!
)[{{E}}\hspace{0.3ex}]
+
e^{2{\mathrm{i}}{\ell}\pi}
(\lambda+\mu^2)
\!\raisebox{0.2ex}{ 
 {\raisebox{0.8ex}{\!$ \raisebox{1ex}[0.0ex][0.0ex]%
  {$\mbox{\raisebox{-2.1ex}[1ex][0.0ex]{\scalebox{0.6}{$\hspace{1.3ex}{{\epsilon}}$} }}
  \atop\rotatebox{180}{\scalebox{0.7}[0.7]{$\curvearrowleft$}}$}\atop
  \mbox{$ \lfloor{\hspace{0.08ex}{\mathfrak D}}\rfloor $}$}\!}
}
\!\!\!\!   \raisebox{0.0ex}{  
{\raisebox{0.7ex}{\!$\raisebox{1.8ex}[0.0ex][0.0ex]%
 {$\mbox{\raisebox{-2.1ex}[1ex][0.0ex]{\scalebox{0.6}{$\hspace{1.3ex}{ 2\epsilon }$} }}
 \atop\rotatebox{180}{\scalebox{1}[1]{$\curvearrowleft$}}$}\atop
 \mbox{$ \!{{E}} $}$}\!}
}\!\!
\big\rgroup 
                              \equiv&
 \\[-3.03ex]
& \hspace{-38ex} 
-{{\lfloor\raisebox{0.1ex}{-}1\rfloor}}%
              ^4 z^2
\!\!\!   \raisebox{0.2ex}{ 
{\raisebox{0.7ex}{\!$\raisebox{1.8ex}[0.0ex][0.0ex]%
 {$\mbox{\raisebox{-2.1ex}[1ex][0.0ex]{\scalebox{0.6}{$\hspace{1.3ex}{ 2\epsilon }$} }}
 \atop\rotatebox{180}{\scalebox{1}[1]{$\curvearrowleft$}}$}\atop
 \mbox{$ \!{{E'}} $}$}\!}
}\!\!\!\!
  \cdot
\bigl(
(\lambda+\mu^2)
{\raisebox{0.8ex}{\!$ \raisebox{1.0ex}[0.0ex][0.0ex]%
 {$\mbox{\raisebox{-2.1ex}[1ex][0.0ex]{\scalebox{0.6}{$\hspace{1.3ex}{{\epsilon}}$} }}
 \atop\rotatebox{180}{\scalebox{0.7}[0.7]{$\curvearrowleft$}}$}\atop
 \mbox{$ {r} $}$}\!}
\!   \raisebox{0.2ex}{ 
{\raisebox{0.8ex}{\!$ \raisebox{1.3ex}[0.0ex][0.0ex]%
 {$\mbox{\raisebox{-2.1ex}[1ex][0.0ex]{\scalebox{0.6}{$\hspace{1.3ex}{{\epsilon}}$} }}
 \atop\rotatebox{180}{\scalebox{0.7}[0.7]{$\curvearrowleft$}}$}\atop
 \mbox{$ {\mbox{$\mathstrut^{B\!}\!\Delta_{p}$}} $}$}\!}
}\!\!
+
(\mu
\!\!   \raisebox{0.2ex}{ 
{\raisebox{0.8ex}{\!$ \raisebox{1.4ex}[0.0ex][0.0ex]%
 {$\mbox{\raisebox{-2.1ex}[1ex][0.0ex]{\scalebox{0.6}{$\hspace{1.3ex}{{\epsilon}}$} }}
 \atop\rotatebox{180}{\scalebox{0.7}[0.7]{$\curvearrowleft$}}$}\atop
 \mbox{$ z^2 $}$}\!}
}  \!\!\!
 {\raisebox{0.8ex}{\!$ \raisebox{1ex}[0.0ex][0.0ex]%
 {$\mbox{\raisebox{-2.1ex}[1ex][0.0ex]{\scalebox{0.6}{$\hspace{1.3ex}{{\epsilon}}$} }}
 \atop\rotatebox{180}{\scalebox{0.7}[0.7]{$\curvearrowleft$}}$}\atop
 \mbox{$ {r} $}$}\!}
\!+\!
 {\raisebox{0.8ex}{\!$ \raisebox{1ex}[0.0ex][0.0ex]%
 {$\mbox{\raisebox{-2.1ex}[1ex][0.0ex]{\scalebox{0.6}{$\hspace{1.3ex}{{\epsilon}}$} }}
 \atop\rotatebox{180}{\scalebox{0.7}[0.7]{$\curvearrowleft$}}$}\atop
 \mbox{$ {s} $}$}\!}
 \!
 )
\!   \raisebox{0.2ex}{ 
{\raisebox{0.8ex}{\!$ \raisebox{1.5ex}[0.0ex][0.0ex]%
 {$\mbox{\raisebox{-2.1ex}[1ex][0.0ex]{\scalebox{0.6}{$\hspace{1.3ex}{{\epsilon}}$} }}
 \atop\rotatebox{180}{\scalebox{0.7}[0.7]{$\curvearrowleft$}}$}\atop
 \mbox{$ {\mbox{$\mathstrut^{B\!}\!\Delta_{r}$}} $}$}\!}
}\!\!
\bigr)
                                \\[-2.5ex]
& \hspace{-38ex}
-
\!\!\!   \raisebox{0.0ex}{ 
{\raisebox{0.7ex}{\!$\raisebox{1.8ex}[0.0ex][0.0ex]%
 {$\mbox{\raisebox{-2.1ex}[1ex][0.0ex]{\scalebox{0.6}{$\hspace{1.3ex}{ {2\epsilon} }$} }}
 \atop\rotatebox{180}{\scalebox{1}[1]{$\curvearrowleft$}}$}\atop
 \mbox{$ \!{{E}} $}$}\!}
}\!\!\!\!
\cdot
(
(\lambda+\mu^2)
{\raisebox{0.8ex}{\!$ \raisebox{1.0ex}[0.0ex][0.0ex]%
 {$\mbox{\raisebox{-2.1ex}[1ex][0.0ex]{\scalebox{0.6}{$\hspace{1.3ex}{{\epsilon}}$} }}
 \atop\rotatebox{180}{\scalebox{0.7}[0.7]{$\curvearrowleft$}}$}\atop
 \mbox{$ r $}$}\!}
\!   \raisebox{0.1ex}{ 
{\raisebox{0.8ex}{\!$ \raisebox{1.5ex}[0.0ex][0.0ex]%
 {$\mbox{\raisebox{-2.1ex}[1ex][0.0ex]{\scalebox{0.6}{$\hspace{1.3ex}{{\epsilon}}$} }}
 \atop\rotatebox{180}{\scalebox{0.7}[0.7]{$\curvearrowleft$}}$}\atop
 \mbox{$ {\mbox{$\mathstrut^{B\!}\!\Delta_{q}$}} $}$}\!}
} \!\!
+
(\mu
\!\!   \raisebox{0.2ex}{ 
{\raisebox{0.8ex}{\!$ \raisebox{1.4ex}[0.0ex][0.0ex]%
 {$\mbox{\raisebox{-2.1ex}[1ex][0.0ex]{\scalebox{0.6}{$\hspace{1.3ex}{{\epsilon}}$} }}
 \atop\rotatebox{180}{\scalebox{0.7}[0.7]{$\curvearrowleft$}}$}\atop
 \mbox{$ z^2 $}$}\!}
}  \!\!\!
\!\!  \raisebox{0.07ex}{ 
{\raisebox{0.8ex}{\!$ \raisebox{1.0ex}[0.0ex][0.0ex]%
 {$\mbox{\raisebox{-2.1ex}[1ex][0.0ex]{\scalebox{0.6}{$\hspace{1.3ex}{{\epsilon}}$} }}
 \atop\rotatebox{180}{\scalebox{0.7}[0.7]{$\curvearrowleft$}}$}\atop
 \mbox{$ {r} $}$}\!}
}  \!\!
\!+\!
\!\!  \raisebox{0.15ex}{  
{\raisebox{0.8ex}{\!$ \raisebox{1.0ex}[0.0ex][0.0ex]%
 {$\mbox{\raisebox{-2.1ex}[1ex][0.0ex]{\scalebox{0.6}{$\hspace{1.3ex}{{\epsilon}}$} }}
 \atop\rotatebox{180}{\scalebox{0.7}[0.7]{$\curvearrowleft$}}$}\atop
 \mbox{$ {s} $}$}\!}
}\!\!\!
 )
\,
 \!\!   \raisebox{0.2ex}{   
 {\raisebox{0.8ex}{\!$ \raisebox{1.5ex}[0.0ex][0.0ex]%
 {$\mbox{\raisebox{-2.1ex}[1ex][0.0ex]{\scalebox{0.6}{$\hspace{1.3ex}{{\epsilon}}$} }}
 \atop\rotatebox{180}{\scalebox{0.7}[0.7]{$\curvearrowleft$}}$}\atop
 \mbox{$ {\mbox{$\mathstrut^{B\!}\!\Delta_{s}$}} $}$}\!}
 }\!\!
)
                                   \\[1.ex]
& \hspace{-38ex}
+
{{\lfloor\raisebox{0.1ex}{-}1\rfloor}}%
             ^4 z^2
{\raisebox{0.8ex}{\!$ \raisebox{1ex}[0.0ex][0.0ex]%
 {$\mbox{\raisebox{-2.1ex}[1ex][0.0ex]{\scalebox{0.6}{$\hspace{1.3ex}{{\epsilon}}$} }}
 \atop\rotatebox{180}{\scalebox{0.7}[0.7]{$\curvearrowleft$}}$}\atop
 \mbox{$ {r} $}$}\!}
\!\!\!  \raisebox{0.1ex}{   
{\raisebox{0.7ex}{\!$\raisebox{1.3ex}[0.0ex][0.0ex]%
 {$\mbox{\raisebox{-2.1ex}[1ex][0.0ex]{\scalebox{0.6}{$\hspace{1.3ex}{ {2\epsilon} }$} }}
 \atop\rotatebox{180}{\scalebox{1}[1]{$\curvearrowleft$}}$}\atop
 \mbox{$ \!{r} $}$}\!}
}\!\!
 \!\!   \raisebox{0.2ex}{  
 {\raisebox{0.7ex}{\!$\raisebox{1.7ex}[0.0ex][0.0ex]%
 {$\mbox{\raisebox{-2.1ex}[1ex][0.0ex]{\scalebox{0.6}{$\hspace{1.3ex}{ {2\epsilon} }$} }}
 \atop\rotatebox{180}{\scalebox{1}[1]{$\curvearrowleft$}}$}\atop
 \mbox{$ {{\mathcal H}}\rlap{$[{{E}}\,]$} $}$}\!}
 }\;\;\;\;
+
{{\lfloor\raisebox{0.1ex}{-}1\rfloor}}%
            ^4 z^2
{\raisebox{0.8ex}{\!$ \raisebox{1ex}[0.0ex][0.0ex]%
 {$\mbox{\raisebox{-2.1ex}[1ex][0.0ex]{\scalebox{0.6}{$\hspace{1.3ex}{{\epsilon}}$} }}
 \atop\rotatebox{180}{\scalebox{0.7}[0.7]{$\curvearrowleft$}}$}\atop
 \mbox{$ {r} $}$}\!}
\!\!\!   \raisebox{0.15ex}{ 
{\raisebox{0.7ex}{\!$\raisebox{1.75ex}[0.0ex][0.0ex]%
 {$\mbox{\raisebox{-2.1ex}[1ex][0.0ex]{\scalebox{0.6}{$\hspace{1.3ex}{ {2\epsilon} }$} }}
 \atop\rotatebox{180}{\scalebox{1}[1]{$\curvearrowleft$}}$}\atop
 \mbox{$ \!{{E'}} $}$}\!}
}\!\!
 \!\!\!  \raisebox{0.15ex}{
 {\raisebox{0.7ex}{\!$\raisebox{1.6ex}[0.0ex][0.0ex]%
 {$\mbox{\raisebox{-2.1ex}[1ex][0.0ex]{\scalebox{0.6}{$\hspace{1.3ex}{ {2\epsilon} }$} }}
 \atop\rotatebox{180}{\scalebox{1}[1]{$\curvearrowleft$}}$}\atop
 \mbox{$ \,\Delta_{{r}} $}$}\!}
 }\!\!
+
{\raisebox{0.8ex}{\!$ \raisebox{1ex}[0.0ex][0.0ex]%
 {$\mbox{\raisebox{-2.1ex}[1ex][0.0ex]{\scalebox{0.6}{$\hspace{1.3ex}{{\epsilon}}$} }}
 \atop\rotatebox{180}{\scalebox{0.7}[0.7]{$\curvearrowleft$}}$}\atop
 \mbox{$ {r} $}$}\!}
\!\!\!   \raisebox{0.15ex}{    
{\raisebox{0.7ex}{\!$\raisebox{1.8ex}[0.0ex][0.0ex]%
 {$\mbox{\raisebox{-2.1ex}[1ex][0.0ex]{\scalebox{0.6}{$\hspace{1.3ex}{ {2\epsilon} }$} }}
 \atop\rotatebox{180}{\scalebox{1}[1]{$\curvearrowleft$}}$}\atop
 \mbox{$ \!{{E}} $}$}\!}
}\!\!
\!\!\!  \raisebox{0.15ex}{
{\raisebox{0.7ex}{\!$\raisebox{1.8ex}[0.0ex][0.0ex]%
 {$\mbox{\raisebox{-2.1ex}[1ex][0.0ex]{\scalebox{0.6}{$\hspace{1.3ex}{ {2\epsilon} }$} }}
 \atop\rotatebox{180}{\scalebox{1}[1]{$\curvearrowleft$}}$}\atop
 \mbox{$ \,\Delta_{{s}} $}$}\!}
 }
                                    \\[0.5ex]
&  \hspace{-38ex}
-{{\lfloor0\rfloor}}
   W_1
+(e^{{{\lfloor0\rfloor}}
       \mu(z
            +{{\lfloor\raisebox{0.1ex}{-}1\rfloor}}
              /z)}-1)
              {{\lfloor\raisebox{0.1ex}{-}1\rfloor}}
               W_2.
\end{aligned}
\end{equation}
\begin{eqnarray}
\mbox{where }
W_1&=&
(\mu(1 - \!\!\! \raisebox{0.1ex}{ 
{\raisebox{0.8ex}{\!$ \raisebox{1.4ex}[0.0ex][0.0ex]%
 {$\mbox{\raisebox{-2.1ex}[1ex][0.0ex]{\scalebox{0.6}{$\hspace{1.3ex}{{\epsilon}}$} }}
 \atop\rotatebox{180}{\scalebox{0.7}[0.7]{$\curvearrowleft$}}$}\atop
 \mbox{$ z^2 $}$}\!}
} \!\!
)
{\raisebox{0.8ex}{\!$ \raisebox{1.1ex}[0.0ex][0.0ex]%
 {$\mbox{\raisebox{-2.1ex}[1ex][0.0ex]{\scalebox{0.6}{$\hspace{1.3ex}{{\epsilon}}$} }}
 \atop\rotatebox{180}{\scalebox{0.7}[0.7]{$\curvearrowleft$}}$}\atop
 \mbox{$ {r} $}$}\!}
 -
{\raisebox{0.8ex}{\!$ \raisebox{1.1ex}[0.0ex][0.0ex]%
 {$\mbox{\raisebox{-2.1ex}[1ex][0.0ex]{\scalebox{0.6}{$\hspace{1.3ex}{{\epsilon}}$} }}
 \atop\rotatebox{180}{\scalebox{0.7}[0.7]{$\curvearrowleft$}}$}\atop
 \mbox{$ {s} $}$}\!}
)
\bigl(
{{\lfloor\raisebox{0.1ex}{-}1\rfloor}}%
             ^4 z^2
 \!\!\raisebox{0.0ex}{ 
 {\raisebox{0.7ex}{\!$\raisebox{1.2ex}[0.0ex][0.0ex]%
 {$\mbox{\raisebox{-2.1ex}[1ex][0.0ex]{\scalebox{0.6}{$\hspace{1.3ex}{ {2\epsilon} }$} }}
 \atop\rotatebox{180}{\scalebox{1}[1]{$\curvearrowleft$}}$}\atop
 \mbox{$ \!{r} $}$}\!}
 }\!\!
\!\!\!\raisebox{0.1ex}{ 
{\raisebox{0.7ex}{\!$\raisebox{1.7ex}[0.0ex][0.0ex]%
 {$\mbox{\raisebox{-2.1ex}[1ex][0.0ex]{\scalebox{0.6}{$\hspace{1.3ex}{ {2\epsilon} }$} }}
 \atop\rotatebox{180}{\scalebox{1}[1]{$\curvearrowleft$}}$}\atop
 \mbox{$ \!{{E'}} $}$}\!}
}\!\!\!
+
\!\!\!\raisebox{0.1ex}{ 
{\raisebox{0.7ex}{\!$\raisebox{1.2ex}[0.0ex][0.0ex]%
 {$\mbox{\raisebox{-2.1ex}[1ex][0.0ex]{\scalebox{0.6}{$\hspace{1.3ex}{ {2\epsilon} }$} }}
 \atop\rotatebox{180}{\scalebox{1}[1]{$\curvearrowleft$}}$}\atop
 \mbox{$ \!{s} $}$}\!}
}\!\!
\!\!\!\raisebox{0.1ex}{ 
{\raisebox{0.7ex}{\!$\raisebox{1.8ex}[0.0ex][0.0ex]%
 {$\mbox{\raisebox{-2.1ex}[1ex][0.0ex]{\scalebox{0.6}{$\hspace{1.3ex}{ {2\epsilon} }$} }}
 \atop\rotatebox{180}{\scalebox{1}[1]{$\curvearrowleft$}}$}\atop
 \mbox{$ \!{{E}} $}$}\!}
}\!\!\!
\bigr),
\nonumber
\\
W_2&=&
\bigl(\!
{\raisebox{0.8ex}{\!$ \raisebox{1ex}[0.0ex][0.0ex]%
 {$\mbox{\raisebox{-2.1ex}[1ex][0.0ex]{\scalebox{0.6}{$\hspace{1.3ex}{{\epsilon}}$} }}
 \atop\rotatebox{180}{\scalebox{0.7}[0.7]{$\curvearrowleft$}}$}\atop
 \mbox{$ {s} $}$}\!}
\!\!\!\raisebox{0.1ex}{ 
{\raisebox{0.7ex}{\!$\raisebox{1.3ex}[0.0ex][0.0ex]%
 {$\mbox{\raisebox{-2.1ex}[1ex][0.0ex]{\scalebox{0.6}{$\hspace{1.3ex}{ {2\epsilon} }$} }}
 \atop\rotatebox{180}{\scalebox{1}[1]{$\curvearrowleft$}}$}\atop
 \mbox{$ \!{s} $}$}\!}
} \!\!\!\!
+
{\raisebox{0.8ex}{\!$ \raisebox{1ex}[0.0ex][0.0ex]%
 {$\mbox{\raisebox{-2.1ex}[1ex][0.0ex]{\scalebox{0.6}{$\hspace{1.3ex}{{\epsilon}}$} }}
 \atop\rotatebox{180}{\scalebox{0.7}[0.7]{$\curvearrowleft$}}$}\atop
 \mbox{$ {r} $}$}\!}
\! \cdot\!
\bigl(
-(\mu(1 - {{\lfloor\raisebox{0.1ex}{-}1\rfloor}}%
                       ^2 z^2)
          +{{\lfloor\raisebox{0.1ex}{-}1\rfloor}}
            ({\ell}-1)z)
 \!\!\! \raisebox{0.1ex}{  
 {\raisebox{0.7ex}{\!$\raisebox{1.3ex}[0.0ex][0.0ex]%
 {$\mbox{\raisebox{-2.1ex}[1ex][0.0ex]{\scalebox{0.6}{$\hspace{1.3ex}{ {2\epsilon} }$} }}
 \atop\rotatebox{180}{\scalebox{1}[1]{$\curvearrowleft$}}$}\atop
 \mbox{$ \!{s} $}$}\!}
 } \!\!\!
+{{\lfloor\raisebox{0.1ex}{-}1\rfloor}}%
             ^3 z^2
\!\!\! \raisebox{0.15ex}{    
{\raisebox{0.7ex}{\!$\raisebox{1.8ex}[0.0ex][0.0ex]%
 {$\mbox{\raisebox{-2.1ex}[1ex][0.0ex]{\scalebox{0.6}{$\hspace{1.3ex}{ {2\epsilon} }$} }}
 \atop\rotatebox{180}{\scalebox{1}[1]{$\curvearrowleft$}}$}\atop
 \mbox{$ \!{s}' $}$}\!}
} \!\!\!
\bigr)\!\bigr)
\! 
{\raisebox{0.7ex}{\!$\raisebox{1.8ex}[0.0ex][0.0ex]%
 {$\mbox{\raisebox{-2.1ex}[1ex][0.0ex]{\scalebox{0.6}{$\hspace{1.3ex}{ {2\epsilon} }$} }}
 \atop\rotatebox{180}{\scalebox{1}[1]{$\curvearrowleft$}}$}\atop
 \mbox{$ \!{{E}} $}$}\!}
\nonumber
 \\
&&
+z^2{{\lfloor\raisebox{0.1ex}{-}1\rfloor}}%
                ^3
\Bigl(
{{\lfloor\raisebox{0.1ex}{-}1\rfloor}}
        {\raisebox{0.8ex}{\!$ \raisebox{1ex}[0.0ex][0.0ex]%
 {$\mbox{\raisebox{-2.1ex}[1ex][0.0ex]{\scalebox{0.6}{$\hspace{1.3ex}{{\epsilon}}$} }}
 \atop\rotatebox{180}{\scalebox{0.7}[0.7]{$\curvearrowleft$}}$}\atop
 \mbox{$ {s} $}$}\!}
        \!\!\!
\raisebox{0.1ex}{    
{\raisebox{0.7ex}{\!$\raisebox{1.3ex}[0.0ex][0.0ex]%
 {$\mbox{\raisebox{-2.1ex}[1ex][0.0ex]{\scalebox{0.6}{$\hspace{1.3ex}{ {2\epsilon} }$} }}
 \atop\rotatebox{180}{\scalebox{1}[1]{$\curvearrowleft$}}$}\atop
 \mbox{$ \!{r} $}$}\!}
}\!\!\!\!\!
 \raisebox{0.2ex}{   
 {\raisebox{0.7ex}{\!$\raisebox{1.6ex}[0.0ex][0.0ex]%
 {$\mbox{\raisebox{-2.1ex}[1ex][0.0ex]{\scalebox{0.6}{$\hspace{1.3ex}{ {2\epsilon} }$} }}
 \atop\rotatebox{180}{\scalebox{1}[1]{$\curvearrowleft$}}$}\atop
 \mbox{$ \!{{E'}} $}$}\!}
 }
\nonumber\\
&&                 \hphantom{ + z^2{{\lfloor\raisebox{0.1ex}{-}1\rfloor}}%
                                                ^3  }
+ 
{\raisebox{0.8ex}{\!$ \raisebox{1.ex}[0.0ex][0.0ex]%
 {$\mbox{\raisebox{-2.1ex}[1ex][0.0ex]{\scalebox{0.6}{$\hspace{1.3ex}{{\epsilon}}$} }}
 \atop\rotatebox{180}{\scalebox{0.7}[0.7]{$\curvearrowleft$}}$}\atop
 \mbox{$ {r} $}$}\!}
\!\cdot
\bigl(( \!  
{\raisebox{0.7ex}{\!$\raisebox{1.3ex}[0.0ex][0.0ex]%
 {$\mbox{\raisebox{-2.1ex}[1ex][0.0ex]{\scalebox{0.6}{$\hspace{1.3ex}{ {2\epsilon} }$} }}
 \atop\rotatebox{180}{\scalebox{1}[1]{$\curvearrowleft$}}$}\atop
 \mbox{$ \!{s} $}$}\!}
\!\!
+
{{\lfloor\raisebox{0.1ex}{-}1\rfloor}}%
             ^4  z^2 \!\!
\raisebox{0.2ex}{ 
{\raisebox{0.7ex}{\!$\raisebox{1.8ex}[0.0ex][0.0ex]%
 {$\mbox{\raisebox{-2.1ex}[1ex][0.0ex]{\scalebox{0.6}{$\hspace{1.3ex}{ {2\epsilon} }$} }}
 \atop\rotatebox{180}{\scalebox{1}[1]{$\curvearrowleft$}}$}\atop
 \mbox{$ {r}' $}$}\!}
}\!
)\!\!
 \raisebox{0.2ex}{ 
 {\raisebox{0.7ex}{\!$\raisebox{1.6ex}[0.0ex][0.0ex]%
 {$\mbox{\raisebox{-2.1ex}[1ex][0.0ex]{\scalebox{0.6}{$\hspace{1.3ex}{ {2\epsilon} }$} }}
 \atop\rotatebox{180}{\scalebox{1}[1]{$\curvearrowleft$}}$}\atop
 \mbox{$ \!{{E'}} $}$}\!}
 }
 \nonumber\\[0.7ex]
 &&               \hphantom{ + z^2{{\lfloor\raisebox{0.1ex}{-}1\rfloor}}%
                                              ^3  }\hphantom{ + {r}\;\;\; 
                                                             ( }\;
+
{{\lfloor\raisebox{0.1ex}{-}1\rfloor}}%
            \! 
{\raisebox{0.7ex}{\!$\raisebox{1.3ex}[0.0ex][0.0ex]%
 {$\mbox{\raisebox{-2.1ex}[1ex][0.0ex]{\scalebox{0.6}{$\hspace{1.3ex}{ {2\epsilon} }$} }}
 \atop\rotatebox{180}{\scalebox{1}[1]{$\curvearrowleft$}}$}\atop
 \mbox{$ {r} $}$}\!}
\!\cdot
(-(\mu (1 - {{\lfloor\raisebox{0.1ex}{-}1\rfloor}}%
             ^2 z^2)
           + {{\lfloor\raisebox{0.1ex}{-}1\rfloor}}%
                         ({\ell}-3)z)
\raisebox{0.2ex}{ 
{\raisebox{0.7ex}{\!$\raisebox{1.6ex}[0.0ex][0.0ex]%
 {$\mbox{\raisebox{-2.1ex}[1ex][0.0ex]{\scalebox{0.6}{$\hspace{1.3ex}{ {2\epsilon} }$} }}
 \atop\rotatebox{180}{\scalebox{1}[1]{$\curvearrowleft$}}$}\atop
 \mbox{$ \!{{E'}} $}$}\!}
}
 \nonumber \\[0.7ex]
 &&                \hphantom{ + z^2{{\lfloor\raisebox{0.1ex}{-}1\rfloor}}^3  }
                   \hphantom{ + {r}\;\;\; 
                             ( }
                    \hphantom{+{{\lfloor\raisebox{0.1ex}{-}1\rfloor}}
                                 \hspace{3.6ex} 
                                  \big(
                              }\;\;\,
+{{\lfloor\raisebox{0.1ex}{-}1\rfloor}}%
             ^3 z^2
\raisebox{0.2ex}{ 
{\raisebox{0.7ex}{\!$\raisebox{1.6ex}[0.0ex][0.0ex]%
 {$\mbox{\raisebox{-2.1ex}[1ex][0.0ex]{\scalebox{0.6}{$\hspace{1.3ex}{ {2\epsilon} }$} }}
 \atop\rotatebox{180}{\scalebox{1}[1]{$\curvearrowleft$}}$}\atop
 \mbox{$ \!{{E''}} $}$}\!}
}
)
\bigr)\Bigr).
\nonumber
\end{eqnarray}
Here we employ, in particular, the following abbreviations
\begin{equation}
\label{eq:710}
 {{\lfloor\raisebox{0.1ex}{-}1\rfloor}}
       =e^{{\mathrm{i}}{{\epsilon}}\pi},\;
 {{\lfloor0\rfloor}}
       =1+e^{{\mathrm{i}}{{\epsilon}}\pi},\;
 {{\lfloor\![0]\!\rfloor}}
       =1+e^{-{\mathrm{i}}{{\epsilon}}\pi},\;
 {{\lfloor2\rfloor}}
       =1-e^{{\mathrm{i}}{{\epsilon}}\pi}.
\end{equation}
Further  
abbreviations used for convenience are as follows:
 $ \lfloor{\hspace{0.08ex}{\mathfrak D}}\rfloor $
stands for 
the right-hand side of {{Eq.}}~\eqref{delTa} considered as a function of $z$.
The operator ${{\mathcal H}}$ is defined by {{Eq.}}~\eqref{sDCHop}.
The symbols
$
\;\Delta{{{\mbox{\tiny\ding{74}}}}}(z)
$, where {{\mbox{\small\ding{74}}}}$\;\in\{{p},{q},{r},{s}\}
$,
denote
the differences of the left- and right-hand sides of {{Eq.s}}{~}%
\eqref{eq:ppp'},
\eqref{eq:qqq'},
\eqref{eq:rrr'},
\eqref{eq:sss'}, respectively, which are considered
as the functions of $z$.
Similarly, the symbols
$
{\raisebox{0.95ex}{$\raisebox{-0.2ex}[0ex][0ex]{\scriptsize$\epsilon$}
\atop\raisebox{-0.2ex}[0ex][0ex]{$\mathstrut^{A\hspace{-0.1ex}}\hspace{-0.5ex}
\Delta_{{{\mbox{\tiny\ding{74}}}}}$}$}}
$
stand for the corresponding `${{\epsilon}}$-deformed' differences
${\mbox{$\mathstrut^{A}\!\Delta_{{{\mbox{\tiny\ding{74}}}}}$}}$ 
of the left- and right-hand sides of {{Eq.s}}{~}\eqref{eq:pppA}-\eqref{eq:sssA}
which are defined by {{Eq.s}}~\eqref{eq:490}.
`The diacritic mark'
{\!\raisebox{2.3ex}[0ex][0ex]{$ \mbox{\raisebox{-2.1ex}[1ex][0.0ex]{
 \scalebox{0.8}{$\hspace{0.1ex}{{\epsilon}}$} }} \atop
 \rotatebox{180}{\scalebox{0.9}[0.9]{$\curvearrowleft$}}$}}%
indicates 
the transformation carrying out the rotation of the function argument
at an angle ${{\epsilon}}\pi$, i.e.\
$\vphantom{{I^I}^I}
 {\raisebox{0.8ex}{\!$ \raisebox{1.3ex}[0.0ex][0.0ex]%
  {$\mbox{\raisebox{-2.1ex}[1ex][0.0ex]{\scalebox{0.6}{$\hspace{1.3ex}{{\epsilon}}$} }}
  \atop\rotatebox{180}{\scalebox{0.7}[0.7]{$\curvearrowleft$}}$}\atop
  \mbox{$ {{\mbox{\small\ding{74}}}} $}$}\!}
(z)= {{\mbox{\small\ding{74}}}}(e^{{\mathrm{i}}{{\epsilon}}\pi} z)$.
Similar `accent'
\!
\!\!{\raisebox{2.3ex}[0ex][0ex]{$ \mbox{\raisebox{-2.1ex}[1ex][0.0ex]{
\scalebox{0.8}{$\hspace{0.1ex}2{{\epsilon}}$} }} \atop
 \rotatebox{180}{\scalebox{0.9}[0.9]{$\curvearrowleft$}} $}}\!\!\!
carries the concordant meaning:
as compared  to 
{\!\raisebox{2.3ex}[0ex][0ex]{$ \mbox{\raisebox{-2.1ex}[1ex][0.0ex]{
 \scalebox{0.8}{$\hspace{0.1ex}{{\epsilon}}$} }} \atop
 \rotatebox{180}{\scalebox{0.9}[0.9]{$\curvearrowleft$}}$}}%
\!\!,
the rotation angle for it amounts to $2{{\epsilon}}\pi  $.

Some abuse of notations is
related
 to usage of 
the symbols
$
\!\!   \raisebox{0.5ex}{
             {\raisebox{0.8ex}{\!$ \raisebox{1.1ex}[0.0ex][0.0ex]%
 {$\mbox{\raisebox{-2.1ex}[1ex][0.0ex]{\scalebox{0.6}{$\hspace{1.3ex}{{\epsilon}}$} }}
 \atop\rotatebox{180}{\scalebox{0.7}[0.7]{$\curvearrowleft$}}$}\atop
 \mbox{$ {\mbox{$\mathstrut^{B\!}\!\Delta_{\rlap{{{\mbox{\tiny\ding{74}}}}}\,\,}$}} $}$}\!}
                   }
$.
Similarly to the symbols
$
{\raisebox{0.95ex}{$\raisebox{-0.2ex}[0ex][0ex]{\scriptsize$\epsilon$}
\atop\raisebox{-0.2ex}[0ex][0ex]{$\mathstrut^{A\hspace{-0.1ex}}\hspace{-0.5ex}
\Delta_{{{\mbox{\tiny\ding{74}}}}}$}$}}
$,
they refer to the
 `${{\epsilon}}$-deformed' differences
of the left- and right-sides of (this time) {{Eq.s}}~\eqref{eq:pppB}-\eqref{eq:sssB},
which are defined   by {{Eq.s}}~\eqref{eq:580}.
However, now additionally ``the
{\!\raisebox{2.3ex}[0ex][0ex]{$ \mbox{\raisebox{-2.1ex}[1ex][0.0ex]{
 \scalebox{0.8}{$\hspace{0.9ex}{{\epsilon}}$} }} \atop
 \rotatebox{180}{\scalebox{0.9}[0.9]{$\curvearrowleft$}}$}}%
\!-rotation'' of the argument  of the ${{\epsilon}}$-deformation result
considered as a function of $z$
has to be carried out afterwards. Thus in this case one deals, in a sense,  with the
``${{\epsilon}}$-deformed and 
{\!\raisebox{2.3ex}[0ex][0ex]{$ \mbox{\raisebox{-2.1ex}[1ex][0.0ex]{
 \scalebox{0.8}{$\hspace{0.1ex}{{\epsilon}}$} }} \atop
 \rotatebox{180}{\scalebox{0.9}[0.9]{$\curvearrowleft$}}$}}%
\!-rotated'' differences of
the left- and right-hand sides of {{Eq.s}}~\eqref{eq:pppB}-\eqref{eq:sssB}.

If ${{\epsilon}}=0$ then all the instances of the functions
${{E}}, {p}, {q}, {r}, {s}$ and their derivatives  involved in
{{Eq.s}}~\eqref{eq:L_A o L_A}, \eqref{eq:L_B o L_B}
are evaluated either at $z$ or at $1/z$.
As a consequence, 
all the constituents of the these formulas
are well defined  for arbitrary $z\in   {{\mbox{$\mathstrut^\backprime{}\mathbb{C}^*$}}}$
--- and the equalities they represent hold true.

Next, 
we allow 
the parameter
$ {{\epsilon}} $ to vary through the segment
$[0,1]$ and carry out analytic continuation
along the corresponding curves (in fact, circular arcs) in the function domains.
At their end points corresponding to ${{\epsilon}} =1$
the
coefficients represented by the abbreviations
 ${{\lfloor\raisebox{0.1ex}{-}1\rfloor}},
 {{\lfloor\![\mbox{-}1]\!\rfloor}},
 {{\lfloor0\rfloor}},
 {{\lfloor2\rfloor}}$
acquire the values $-1,-1,0,2$, respectively, see {{Eq.s}}~\eqref{eq:710}.
This allows us,
in particular, to ignore
the last lines in the both
formulas
\eqref{eq:L_A o L_A} and \eqref{eq:L_B o L_B}. 
Simultaneously,
the effect of
the rotation  of the function arguments
tagged by `the accent' 
{\!\raisebox{2.3ex}[0ex][0ex]{$ \mbox{\raisebox{-2.1ex}[1ex][0.0ex]{
 \scalebox{0.8}{$\hspace{0.1ex}{{\epsilon}}$} }} \atop
 \rotatebox{180}{\scalebox{0.9}[0.9]{$\curvearrowleft$}}$}}
turns into  the 
action of the semi-monodromy operator ${{{\mathcal{M}}}^{1/2}}$
(see Theorem \ref{t:030})
which we indicate also by `the accent'
$ \raisebox{1.5ex} {\rotatebox{185}{\scalebox{0.99}[0.99]{$\curvearrowleft$}}} $ %
over the function symbol,
 see, e.g., the proof of   Theorem \ref{t:060}.
As to the exceptional symbols
\raisebox{0ex}[2.8ex]{
$
 \!\!\!   \raisebox{0.2ex}{
                          {\raisebox{0.8ex}{\!$ \raisebox{1.1ex}[0.0ex][0.0ex]%
 {$\mbox{\raisebox{-2.1ex}[1ex][0.0ex]{\scalebox{0.6}{$\hspace{1.3ex}{{\epsilon}}$} }}
 \atop\rotatebox{180}{\scalebox{0.7}[0.7]{$\curvearrowleft$}}$}\atop
 \mbox{$ {\mbox{$\mathstrut^{B\!}\!\Delta_{\rlap{{{\mbox{\tiny\ding{74}}}}}\,\,}$}} $}$}\!}
                         }
$},
{{\mbox{\small\ding{74}}}}{}$\;\in\{{p},{q},{r},{s}\} $,
it is easy to see that
at the end point of the curve of analytic continuation
they become equal to the expressions
denoted in {{Eq.}}~\eqref{eq:440}
by the symbols
\raisebox{0ex}[2.5ex]{$\,
{\hspace{-0.2em}\raisebox{1.3ex}{\!$\raisebox{0.4ex}[0.0ex][0.0ex]%
 {$\rotatebox{185}{\scalebox{0.7}[0.9]{$\curvearrowleft$}} $} \atop
 \raisebox{-0.5ex}[0.0ex][0.0ex]{$
 {\mbox{$\mathstrut^{{B}\!}\!\Delta_{ {{{\mbox{\tiny\ding{74}}}}} }$}}
 $}$}\hspace{-0.2em}}
\,$}.

A separate note  on the effect of `the accent'
\!
\!\!{\raisebox{2.3ex}[0ex][0ex]{$ \mbox{\raisebox{-2.1ex}[1ex][0.0ex]{
\scalebox{0.8}{$\hspace{0.1ex}2{{\epsilon}}$} }} \atop
 \rotatebox{180}{\scalebox{0.9}[0.9]{$\curvearrowleft$}} $}}\!\!\!
is necessary. In the limit as ${{{\epsilon}}\nearrow 1}$
it also `rotates' the argument of the function
to be transformed
but now the rotation 
angle amounts to $2\pi$
meaning, in a sense, a full revolution. 
It had been noticed that
such transformations are termed monodromy.
We denoted the operator carrying out the monodromy transformation
by the symbol ${\mathcal{M}}$ but in some formulas (e.g.\ in {{Eq.}}~\eqref{eq:440}) it is also
indicated by `the diacritic mark'
{\raisebox{0.7ex}{\scalebox{0.7}[0.7]{$\circlearrowleft$}}}%
. 

It has also to be noted that in case of monodromy transformation
some precaution on structure of the domain of the function
to which it acts needs to be taken. This point is briefly discussed in the proof
of Theorem \ref{t:060}.

Now, collecting all the modifications of the formulas
{{Eq.s}}~\eqref{eq:L_A o L_A} and \eqref{eq:L_B o L_B},
arising when the analytic continuation corresponding to
$\lim_{{{\epsilon}}\nearrow 1}$ has been carried out,
one finds 
that they finally convert to {{Eq.s}}~\eqref{eq:430} and \eqref{eq:440},
respectively.

%
%
%
%
%
%

\end{document}